\documentclass[11pt,a4paper]{article}
\usepackage[english]{babel}
\usepackage{amsmath, empheq}
\usepackage{graphicx}
\usepackage{parskip} 
\usepackage{subcaption} 
\usepackage{multirow}	
\usepackage{float} 
\usepackage{siunitx}
\usepackage{color}
\usepackage{marginnote}
\usepackage{amssymb} 
\usepackage{mathtools} 
\usepackage{verbatim}
\usepackage{tikz}
\usepackage{Bert_ii_variables}
\usepackage{Bert_ii_variables_afleidingen}
\usepackage{Bert_ii_bestestimator}
\usepackage{Bert_ii_gridsketch}
\usepackage{Bert_ii_zeroD_oneD_paths}
\usepackage[makeroom]{cancel}
\usepackage{adjustbox}
\usepackage{calc}

\mathtoolsset{showonlyrefs}
\usepackage[a4paper,text={6.5in,10.0in},centering,
            includefoot,foot=0.6in]{geometry}
\usepackage{authblk}

\usepackage{amsthm} 
\usepackage{thmtools}

\tikzstyle{neutraledge}=[gray]
\tikzstyle{deadedge}=[red,thick,dotted]
\tikzstyle{winneredge}=[blue,very thick]
\newcommand{\clippicture}[5]{%
\node[inner sep=0pt] (#5) at (#2,#3) {\includegraphics[width=#4,height=#4,trim=2.52cm 2.25cm 2.02cm 2.42cm,clip]{#1}};
}%
\newcommand{\triplebestfig}[4]{
\def\figloc{#1}%
\begin{tikzpicture}[x=1cm, y=1cm]%
\def\ticklength{.1cm}%
\def\height{5cm}%
\def\width{5cm}%
\def\xspacing{\width/4}%
\clippicture{\figloc_1D0D_anal_paper.pdf}{0}{0}{\height}{1}%
\clippicture{\figloc_1D0D_num_paper.pdf}{\width+\xspacing}{0}{\height}{2}%
\clippicture{\figloc_1D1D_num_paper.pdf}{2*\width+2*\xspacing}{0}{\height}{3}%
\foreach \i in {0,1,2} {%
 \draw (\width*\i+\i*\xspacing-\width/2,-\height/2) -- (\width*\i+\i*\xspacing-\width/2,\height/2);%
 \draw (\width*\i+\i*\xspacing-\width/2,-\height/2) -- (\width*\i+\i*\xspacing-\width/2-\ticklength,-\height/2) node[left] () {\Large 0};%
 \draw (\width*\i+\i*\xspacing-\width/2,\height/2) -- (\width*\i+\i*\xspacing-\width/2-\ticklength,\height/2) node[left] () {\Large 10};%
 \draw (\width*\i+\i*\xspacing-\width/2,-\height/2) -- (\width*\i+\i*\xspacing-\width/2,-\height/2-\ticklength) node[below] () {\Large 0};%
 \draw (\width*\i+\i*\xspacing+\width/2,-\height/2) -- (\width*\i+\i*\xspacing+\width/2,-\height/2-\ticklength) node[below] () {\Large 1};%
 \draw (\width*\i+\i*\xspacing-\width/2,-\height/2) -- (\width*\i+\i*\xspacing+\width/2,-\height/2);%
}%
\def\dummyopacitynaam{1}%
\ifdim .5 pt > #2 pt \relax \def\dummyopacitynaam{0}\fi%
\ifdim .5 pt > #4 pt \relax {\draw (0,\height/2+\xspacing/3) node[text opacity=\dummyopacitynaam] (naam1) {\LARGE analytical};
\draw (\width+\xspacing,\height/2+\xspacing/3) node[text opacity=\dummyopacitynaam] (naam2) {\LARGE numerical};
\draw (2*\width+2*\xspacing,\height/2+\xspacing/3) node[text opacity=\dummyopacitynaam] (naam3) {\LARGE numerical};
\draw (-\width/2-2*\ticklength, 0) node[above,rotate=90] () {\LARGE $\VAratet\VAdomlength/|v|$};}\fi%
\ifdim .5 pt < #4 pt \relax {\draw (0,\height/2+\xspacing/3) node[text opacity=\dummyopacitynaam] (naam1) {\LARGE mass source};
\draw (\width+\xspacing,\height/2+\xspacing/3) node[text opacity=\dummyopacitynaam] (naam2) {\LARGE momentum source};
\draw (2*\width+2*\xspacing,\height/2+\xspacing/3) node[text opacity=\dummyopacitynaam] (naam3) {\LARGE momentum source};
\draw (-\width/2-2*\ticklength, 0) node[above,rotate=90] () {\LARGE $\VAratet\VAdomlength/|v|$};}\fi%
\draw (naam1) node[yshift=2*\ticklength,above,text opacity=\dummyopacitynaam] () {\LARGE 1D0D};%
\draw (naam2) node[yshift=2*\ticklength,above,text opacity=\dummyopacitynaam] () {\LARGE 1D0D};%
\draw (naam3) node[yshift=2*\ticklength,above,text opacity=\dummyopacitynaam] () {\LARGE 1D1D};%
\foreach \i in {0,1,2} {%
 \def\dummyopacitybot{1}%
 \ifdim .5 pt > #3 pt \relax \def\dummyopacitybot{0}\fi%
 \draw (\width*\i+\xspacing*\i,-\height/2-2*\ticklength) node[below,text opacity=\dummyopacitybot] () {\LARGE $\VAPr$};%
}%
\end{tikzpicture}%
}%
\newcommand{\maakmooieticks}[4]{
\def\figfullloc{#1}%
\begin{tikzpicture}[x=1cm, y=1cm]%
\def\ticklength{.1cm}%
\def\height{5cm}%
\def\width{5cm}%
\def\xspacing{\width/4}%
\clippicture{\figfullloc}{0}{0}{\height}{1}%
\foreach \i in {0} {%
 \draw (\width*\i+\i*\xspacing-\width/2,-\height/2) -- (\width*\i+\i*\xspacing-\width/2,\height/2);%
 \draw (\width*\i+\i*\xspacing-\width/2,-\height/2) -- (\width*\i+\i*\xspacing-\width/2-\ticklength,-\height/2) node[left] () {\Large 0};%
 \draw (\width*\i+\i*\xspacing-\width/2,\height/2) -- (\width*\i+\i*\xspacing-\width/2-\ticklength,\height/2) node[left] () {\Large 10};%
 \draw (\width*\i+\i*\xspacing-\width/2,-\height/2) -- (\width*\i+\i*\xspacing-\width/2,-\height/2-\ticklength) node[below] () {\Large 0};%
 \draw (\width*\i+\i*\xspacing+\width/2,-\height/2) -- (\width*\i+\i*\xspacing+\width/2,-\height/2-\ticklength) node[below] () {\Large 1};%
 \draw (\width*\i+\i*\xspacing-\width/2,-\height/2) -- (\width*\i+\i*\xspacing+\width/2,-\height/2);%
}%
\def\dummyopacityRtLv{1}%
\ifdim .5 pt > #2 pt \relax \def\dummyopacityRtLv{0}\fi%
\draw (0,\height/2+\xspacing/1.5) node[] (naam1) {\LARGE #4};
\draw (-\width/2-2*\ticklength, 0) node[above,rotate=90,text opacity=\dummyopacityRtLv] () {\LARGE $\VAcst\VAdomlength$};%
\foreach \i in {0} {%
 \def\dummyopacitybot{1}%
 \ifdim .5 pt > #3 pt \relax \def\dummyopacitybot{0}\fi%
 \draw (\width*\i+\xspacing*\i,-\height/2-2*\ticklength) node[below,text opacity=\dummyopacitybot] () {\LARGE $\VAPr$};%
}%
\end{tikzpicture}%
}%
\newcommand{\maakmooietickscb}[2]{
\def\figfullloc{#1}%
\begin{tikzpicture}[x=1cm, y=1cm]%
\def\ticklength{.1cm}%
\def\height{5cm}%
\def\width{1cm}%
\def\xspacing{\width/4}%
\node[inner sep=0pt] (0) at (0,0) {\includegraphics[width=\width,height=\height,trim=1cm 4.25cm 1cm 4.53cm,clip]{\figfullloc}};
\foreach \i in {0} {%
 \draw (\width*\i+\i*\xspacing-\width/2,-\height/2) -- (\width*\i+\i*\xspacing-\width/2,\height/2);
 \draw (\width*\i+\i*\xspacing+\width/2,-\height/2) -- (\width*\i+\i*\xspacing+\width/2,\height/2);
 \draw (\width*\i+\i*\xspacing+\width/2-\ticklength,-\height/2) -- (\width*\i+\i*\xspacing+\width/2+\ticklength,-\height/2) node[right] () {\Large $10^0$};
 \draw (\width*\i+\i*\xspacing+\width/2-\ticklength,0) -- (\width*\i+\i*\xspacing+\width/2+\ticklength,0) node[right] () {\Large $10^1$};%
 \draw (\width*\i+\i*\xspacing+\width/2-\ticklength,\height/2) -- (\width*\i+\i*\xspacing+\width/2+\ticklength,\height/2) node[right] () {\Large $10^2$};%
\foreach \smalltickpos in {.1505,.2386,.301,.3495,.3891,.4225,.4515,.4771}{%
 \draw (\width*\i+\i*\xspacing+\width/2+\ticklength/1.5,-\height/2+\smalltickpos*\height) -- (\width*\i+\i*\xspacing+\width/2,-\height/2+\smalltickpos*\height) node[left] () {};%
 \draw (\width*\i+\i*\xspacing+\width/2+\ticklength/1.5,\smalltickpos*\height) -- (\width*\i+\i*\xspacing+\width/2,\smalltickpos*\height) node[left] () {};%
}%
 \draw (\width*\i+\i*\xspacing-\width/2,-\height/2) -- (\width*\i+\i*\xspacing+\width/2,-\height/2);%
  \draw (\width*\i+\i*\xspacing-\width/2,\height/2) -- (\width*\i+\i*\xspacing+\width/2,\height/2);%
  \draw (\width*\i+\xspacing*\i,-\height/2-2*\ticklength) node[below,text opacity=0] () {\LARGE $\VAPr$};%
  \draw (0,\height/2+\xspacing/1.5) node[text opacity=0] (naam1) {\LARGE #2};
}%
\end{tikzpicture}%
}%
\usepackage{soul}

\title{A comparison of source term estimators in coupled finite-volume/Monte-Carlo methods with applications to plasma edge simulations in nuclear fusion}
\author[1]{Bert Mortier}
\author[2]{Martine Baelmans}
\author[1]{Giovanni Samaey} 

\affil[1]{NUMA, Dept. Computer Science, KU Leuven}
\affil[2]{Dept. Mechanical Engineering, KU Leuven}

\renewcommand\vec{\mathbf}

\newcommand{\localcolorlegendoneA}{
\begin{tabular}{cl}
{{\adjustbox{trim={.45\width} {.45\height} {.45\width} {.45\height}, clip}{\includegraphics[width=2cm]{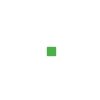}}}} & \texttt{nac\_ne}\\
{{\adjustbox{trim={.45\width} {.45\height} {.45\width} {.45\height}, clip}{\includegraphics[width=2cm]{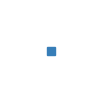}}}} & \texttt{nac\_tl}\\
{{\adjustbox{trim={.45\width} {.45\height} {.45\width} {.45\height}, clip}{\includegraphics[width=2cm]{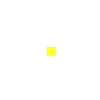}}}} & \texttt{natl\_ne}\\
{{\adjustbox{trim={.45\width} {.45\height} {.45\width} {.45\height}, clip}{\includegraphics[width=2cm]{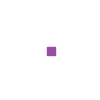}}}} & \texttt{natl\_tl}\\
\end{tabular}}

\newcommand{\localcolorlegendoneB}{
\begin{tabular}{cl}
{{\adjustbox{trim={.45\width} {.45\height} {.45\width} {.45\height}, clip}{\includegraphics[width=2cm]{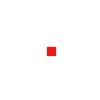}}}} & \texttt{a\_ne}\\
{{\adjustbox{trim={.45\width} {.45\height} {.45\width} {.45\height}, clip}{\includegraphics[width=2cm]{figuren/nac_ex.png}}}} & \texttt{nac\_ne}\\
{{\adjustbox{trim={.45\width} {.45\height} {.45\width} {.45\height}, clip}{\includegraphics[width=2cm]{figuren/natl_ex.png}}}} & \texttt{natl\_ne}\\
{{\adjustbox{trim={.45\width} {.45\height} {.45\width} {.45\height}, clip}{\includegraphics[width=2cm]{figuren/natl_tl.png}}}} & \texttt{natl\_tl}\\
\end{tabular}}

\newcommand{\localcolorlegendoneC}{
\begin{tabular}{cl}
{{\adjustbox{trim={.45\width} {.45\height} {.45\width} {.45\height}, clip}{\includegraphics[width=2cm]{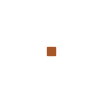}}}} & \texttt{nac\_c}\\
{{\adjustbox{trim={.45\width} {.45\height} {.45\width} {.45\height}, clip}{\includegraphics[width=2cm]{figuren/nac_ex.png}}}} & \texttt{nac\_ne}\\
{{\adjustbox{trim={.45\width} {.45\height} {.45\width} {.45\height}, clip}{\includegraphics[width=2cm]{figuren/natl_ex.png}}}} & \texttt{natl\_ne}\\
{{\adjustbox{trim={.45\width} {.45\height} {.45\width} {.45\height}, clip}{\includegraphics[width=2cm]{figuren/natl_tl.png}}}} & \texttt{natl\_tl}\\
\end{tabular}}

\newcommand{\localcolorlegendoneD}{
\begin{tabular}{cl}
{{\adjustbox{trim={.45\width} {.45\height} {.45\width} {.45\height}, clip}{\includegraphics[width=2cm]{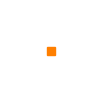}}}} & \texttt{a\_c}\\
{{\adjustbox{trim={.45\width} {.45\height} {.45\width} {.45\height}, clip}{\includegraphics[width=2cm]{figuren/a_ex.png}}}} & \texttt{a\_ne}\\
{{\adjustbox{trim={.45\width} {.45\height} {.45\width} {.45\height}, clip}{\includegraphics[width=2cm]{figuren/nac_ex.png}}}} & \texttt{nac\_ne}\\
{{\adjustbox{trim={.45\width} {.45\height} {.45\width} {.45\height}, clip}{\includegraphics[width=2cm]{figuren/natl_ex.png}}}} & \texttt{natl\_ne}\\
{{\adjustbox{trim={.45\width} {.45\height} {.45\width} {.45\height}, clip}{\includegraphics[width=2cm]{figuren/natl_tl.png}}}} & \texttt{natl\_tl}\\
\end{tabular}}

\newcommand{\localcolorlegendoneE}{
\begin{tabular}{cl}
{{\adjustbox{trim={.45\width} {.45\height} {.45\width} {.45\height}, clip}{\includegraphics[width=2cm]{figuren/nac_ex.png}}}} & \texttt{nac\_ne}\\
{{\adjustbox{trim={.45\width} {.45\height} {.45\width} {.45\height}, clip}{\includegraphics[width=2cm]{figuren/natl_ex.png}}}} & \texttt{natl\_ne}\\
{{\adjustbox{trim={.45\width} {.45\height} {.45\width} {.45\height}, clip}{\includegraphics[width=2cm]{figuren/natl_tl.png}}}} & \texttt{natl\_tl}\\
\end{tabular}}

\begin{document}

\def \varmapmetfiguren {figuren/}
\def \varmapmetfigurenresults {figuren/}

\maketitle

\begin{abstract}
In many applications, such as plasma edge simulation of a nuclear fusion reactor, a coupled PDE/kinetic description is required, which is usually solved with a coupled finite-volume/Monte-Carlo method. Different procedures have been proposed to estimate the source terms in the finite volume part that appear from the Monte Carlo part of the simulation. In this paper, we present a systematic comparison of the variance and computational cost of a coherent set of such estimation procedures. We compare the different estimation procedures for mass in a simplified forward-backward scattering model problem, where an analytical comparison is possible, and for mass and momentum in a model problem with realistic scattering. Our results reveal a non-trivial dependence of the optimal choice of estimator on the model parameters and show that different estimation procedures prevail for different quantities of interest.
\end{abstract}


\section{Introduction}

In many applications, the relevant processes are accurately described by a partial differential equation (PDE) of reaction-advection-diffusion type coupled to a Boltzmann-type kinetic equation that models a distribution of particles in position-velocity space. Examples of such applications include bacterial chemotaxis~\cite{rousset2013bacterialchemotaxis}, rarified gas dynamics~\cite{sun2004hybridrarifiedgas}, and plasma physics~\cite{reiter2005eirene,stangeby2000plasmaboundary,viola2014snowflakevsconv}. The simulation of such coupled models is computationally challenging due to the different dimensionality of both parts of the model.

This work focuses on particle-tracing methods for the Boltzmann-BGK kinetic equation, which are implemented in the neutral particle transport codes EIRENE~\cite{reiter2005eirene} and DEGAS2~\cite{stotler1994degas2} for nuclear fusion applications. These codes are coupled to PDE solvers for the plasma, such as the deterministic B2~\cite{reiter2005eirene} or UEDGE~\cite{fenstermacher1995uedgedegas} codes or the stochastic EMC3 code~\cite{feng2002emc3}. Such coupled codes are used to perform plasma edge simulations to evaluate and design fusion reactors and operational conditions. Due to evolutions in fusion reactor design, the Monte Carlo simulations are often the computational bottleneck due to highly heterogeneous and highly reactive conditions~\cite{krasheninnikov2017physicsdetachement}. These conditions therefore require careful selection of the most suitable Monte Carlo estimation procedure, or (potentially) the use of different estimation procedures in different parts of the space-time domain. This forms the motivation for this paper, in which we study these particle tracing methods at length.

In this paper, we consider a prototypical model of this type that appears in plasma edge simulations in nuclear fusion reactors, such as ITER \cite{wiesen2015SOLPSITER}. The neutral model contains the interactions with the plasma and during these interactions, mass, momentum, and energy are exchanged between plasma and neutral particles. The exchanges between neutrals and plasma are modeled as source terms in the plasma equations and estimating them is the aim of the neutral model simulation. We will discuss this neutral model, together with three different unbiased simulation strategies and four different source term estimators. Combined, this yields eleven different relevant procedures to estimate the source terms for the plasma equation.

Different source term estimation procedures result in different statistical behaviour and a different computational cost. Selecting the best estimator forms a fundamental way of reducing the variance and simulation cost of the Monte Carlo procedure~\cite{kahn1956MCapplications}. Currently, only a few works are available that compare the performance of these estimation procedures, usually in a very restrictive setting, leaving the choice of estimation procedure to the preference or experience of the user. A first important comparison of the estimation procedures was conducted in~\cite{macmillan1966comparisonneutrMC}, where an invariant imbedding methodology~\cite{bellman1975ii} is used to derive ODEs for the statistical error of a limited set of estimation procedures in a one-dimensional forward scattering scenario. In this degenerate scenario, neutrals always have the same velocity and do not change direction, leading to very simple particle paths. Indira~\cite{indira1988analyticalleakageest} performed a similar study that included forward-backward scattering, but limited the study to estimation procedures for leakage, the number of particles that leave the domain. Both the setting in \cite{macmillan1966comparisonneutrMC} and \cite{indira1988analyticalleakageest} allowed for significant simplifications, resulting in ODEs for the statistical error with a comprehensive analytical solution. We also refer to the work of Lux for approximate formulas for the variance of the most commonly used estimators \cite{lux1979vareff_techrep}, and for sufficient conditions for one estimator to outperform an other one \cite{lux1978standardvarred}. While useful, Lux' results do not capture the highly non-trivial behaviour at high scattering rates and low absorption rates, where the paths are generally the most complex. For completeness, we also refer to \cite{spanier1970analyticvarred, sarkar1979scoremomentumIS, indira1986weightmomentseqexpandscattangl, indira1989optimization, solomon2011scoremomentumWW}, in which analytical calculations of the variance are presented to optimize importance sampling.

This paper forms a next step in the systematic study of source term estimators in coupled finite-volume/Monte-Carlo methods. In particular, we discuss the different procedures that can be used for the Monte Carlo estimation of the mass and momentum source terms, and find which estimation procedure attains the lowest statistical error and computational cost as a function of the plasma background parameters. We first consider mass source estimation in a simplified backward-forward scattering model problem, for which analytical expressions of the statistical error and computational error can be found via the invariant imbedding methodology. This problem, although simplified, represents a significant complication over previous work. We extend our analytical work with experiments that cover both mass and momentum source estimation in a multi-speed setting where an invariant imbedding procedure would become infeasible. This numerical extension shows the relevance of the highly simplified forward-backward scattering model for selection of the mass source estimation procedure, but also the large difference between optimal choice of estimation procedure for mass source estimation compared to momentum source estimation. In Chapters~2 and~3 of~\cite{mortier2020phdthesis}, a more elaborate version of the here presented results are given.

The remainder of this paper is organized as follows. In Section~\ref{sec:iif_modelsim}, we briefly describe the model problem and the three simulation types we consider. Then, in Section~\ref{sec:iif_est}, we present the four estimator types and present the eleven relevant source term estimation procedures we arrive at. In Section~\ref{sec:iif_bestest_1D0D}, we present the simplified problem and the invariant imbedding procedure that enables us to calculate the statistical properties for each of these estimation procedures. With these calculated statistical properties, we identify the best estimation procedure throughout the parameter domain and we present the potential gain in doing so. In Section~\ref{sec:ii_res_1D1D}, we expand on this by considering mass and momentum estimation in a more relevant but more complex setting, for which we conduct and present the numerical experiments. In Section~\ref{sec:iif_conclusion}, we finally present our conclusions and hint at future work.

\section{Kinetic neutral model and its Monte Carlo simulation\label{sec:iif_modelsim}}

In this section, we briefly review the kinetic neutral model and the three Monte Carlo discretizations we consider. A detailed description and derivation of these methods can be found in~\cite[Chapter 2]{mortier2020phdthesis}.

The neutrals in a fusion reactor originate according to a source term $\VAnsource(x,v)$ and then move with their velocity, until they undergo a collision. This collision can either be with the wall or with plasma particles. Generally, a collision with a plasma particle can either be an ionization (absorption) at which the neutral disappears from the neutral simulation, or a charge-exchange interaction (scattering), which is modeled by the particle receiving a new velocity from the post-collisional velocity distribution $\VApostcolveldistrvr$. With absorption rate $\VAratea(x)$ and scattering rate $\VArates(x)$, this behaviour can be written as the following stationary kinetic equation for the neutral particle distribution $\phi_\text{n}(x,v)$:
\begin{equation}
\underbrace{\vphantom{\int}{v}\nabla\phi_\text{n}(x,v)}_\text{transport}=\underbrace{\VAnsource(x,v)}_{\substack{\text{source from}\\ \text{the plasma}}}-\underbrace{\vphantom{\int}\VAratea({x})\phi_\text{n}(x,v)}_{\substack{\text{sink due to}\\ \text{absorption}}}
-\underbrace{\vphantom{\int}\VArates({x})\phi_\text{n}(x,v)+\int\!\!\VArates({x})\phi_\text{n}(x,v')\VApostcolveldistrvr\text{d}{v}'}_{\substack{\text{velocity redistribution}\\ \text{due to scattering}}}\,,\label{eq:iif_modelsim_model_kinetic}
\end{equation}
where we ignored boundary hits.

We consider three types of simulation:
\begin{itemize}
\item \textbf{Analog simulation.} A first method to obtain samples according to Equation~\eqref{eq:iif_modelsim_model_kinetic}, is by modeling particle paths $(x(t),v(t))$ that adhere to the underlying particle model that gave rise to the kinetic equation. For that reason, this simulation is called \emph{analog} and referred to by \texttt{a}. We simulate particles by letting them undergo collisions with the plasma, which can be either an absorption or a scattering. Upon scattering, particles change velocity. Particles disappear from the simulation as soon as they undergo an absorption event. A potential weakness in this method, is that with the finite number of particles, potentially zero of them could penetrate highly absorbing regions in the domain. This caveat is resolved by the two subsequent simulation methods.
\item \textbf{Non-analog collision type simulation.} A second method to obtain samples according to Equation~\eqref{eq:iif_modelsim_model_kinetic} also simulates particles by sampling from $\VAnsource(x,v)$ and letting them undergo collisions with the plasma at the total rate $R_\text{a}(x)+R_\text{s}(x)$. However, instead of executing either an absorption or a scattering effect, the collision-type simulation \emph{always} executes a scattering event, and takes care of absorption by introducing a particle weight that is updated at each collision. This no longer follows the underlying particle model of absorption exactly, but instead modifies the collisions, hence the name \emph{non-analog collision type simulation}, which we refer to by \texttt{nac}.
\item \textbf{Non-analog track-length type simulation.} Finally, a third method takes care of absorption by continuously adjusting the particle weight along the trajectory. The only collision events that remain are then scattering events, which occur with the rate $R_\text{s}(x)$. This method is called the \emph{non-analog track-length type simulation} and referred to by \texttt{natl}.
\end{itemize}

On top of collisions, particles can also undergo boundary hits. A usual boundary condition is partial reflection, which can also be implemented as a probability of disappearing from the simulation or as a weight update. The rates $R_\text{a}(x)$ and $R_\text{s}(x)$ depend on the plasma state, which is usually given as a piecewise constant function on a finite volume grid. To keep track of the local rates as the particle moves, we include grid cell crossings as a third type of event. These will be of practical importance for some estimators presented in Section~\ref{sec:iif_est}.

\paragraph{Equations of motion.} The particle state is represented as a function $t\in[0,T_\text{end}]\mapsto(x(t),v(t),w(t))$, with $x(t)$ the time-dependent position, $v(t)$ the velocity and $w(t)$ the particle weight. At the initial time $T_0=0$ the initial state, $(x_0,v_0)$ samples $\VAnsource(x,v)$, and $w(0)=1$. From there, the particle moves with a constant velocity, that potentially changes at every event.  
\begin{align}
(x_0,v_0)&\sim\VAnsource(x,v)\,,\\
x(0)&=x_0\,,\\
\frac{\text{d}x(t)}{\text{d}t}&=v(t)\,,\\
v(t)&=v_\VAeventno\text{  for  }t\in[T_\VAeventno,T_{\VAeventno+1}),\ \ k\in\{0,\dots,\VAnofevents-1\}\,.
\end{align}
The $T_\VAeventno$ form the event times and $T_\VAnofevents=T_\text{end}$ is the time instance at which the particle disappears from the simulation. This final event can either be an absorption collision or the particle leaving the domain through the boundary. The velocity at this final time does not change, ${v}(T_\VAnofevents)=\VAdiscreteveleend_{\VAnofevents-1}$. We now first look into the calculation of the event times $T_k$, after which we discuss the generation of a new velocity and the updating of the particle weight.

\paragraph{Event times.} If we ignore the possibility of an other event taking place first, the next time at which an event takes place of each of the three kinds separately: collision with the plasma, boundary hit, or grid cell crossing, the next time, can be readily found from $(x(T_{\VAeventno}),v(T_{\VAeventno}))$. We will discuss each of the three events separately now and then find the actual next event as the one that occurs first. 

The next collision time would be sampled by inverting
\begin{equation}
\int_{T_\VAeventno}^{T_{\VAeventno+1}^c}\VArate_{\ast}(x(T_\VAeventno+(t-T_\VAeventno)v(T_\VAeventno))\text{d}t\,,\text{  with  }\epsilon\sim\mathcal{E}(1)\,,\label{eq:single_particle_scatabsT}
\end{equation}
where $\VArate_{\ast}(x)$ represents the total rate of collisions with the plasma occurring in the simulation. In the analog and non-analog collision type simulations this is $\VAratet(x)=\VAratea(x)+\VArates(x)$ and in the non-analog track-length type simulation, absorption is taken care of along the particle trajectory, and not during collisions, so the total rate of collisions is only $\VArates(x)$.

A boundary hit occurs at the first instance the boundary is reached, hence
\begin{equation}
T_{\VAeventno+1}^b=T_\VAeventno+\text{min}(\tau|\tau>0,x(T_{\VAeventno})+\tau v(T_{\VAeventno})\in\partial\mathcal{D})\,,
\end{equation}
with $\partial\mathcal{D}$ the boundary of the domain.

Similarly to a boundary hit, a grid cell crossing occurs at the first instance a grid cell is reached. With the grid cells denoted by the set $\left\{\VAgridcelldefgridcellno\right\}_{\VAgridcellno=1}^\VAnofgridcells$, the grid cell boundary can be written as $\cup_\VAgridcellno\partial\VAgridcelldefgridcellno$ and the grid cell crossing time, not considering other events, is
\begin{equation}
T_{\VAeventno+1}^g=T_\VAeventno+\text{min}(\tau|\tau>0,x(T_{\VAeventno})+\tau v(T_{\VAeventno})\in\cup_\VAgridcellno\partial\VAgridcelldefgridcellno\setminus\partial\mathcal{D})\,.
\end{equation}

The next event that actually occurs is that for which $T_{\VAeventno+1}^\star$, $\star\in\{c,b,g\}$, is the lowest. To encode the nature of the $({\VAeventno+1})$-th event, we use the variables $c_{\VAeventno+1}$ for a collision with the plasma, $b_{\VAeventno+1}$ for a boundary hit, and $g_{\VAeventno+1}$ for a grid cell crossing, such that
\begin{align}
(T_{\VAeventno+1},c_{\VAeventno+1},b_{\VAeventno+1},g_{\VAeventno+1})=\left\{\begin{array}{lll}
(T_{\VAeventno+1}^c,1,0,0) & \text{ if }&T_{\VAeventno+1}^c=\text{min}(T_{\VAeventno+1}^c,T_{\VAeventno+1}^b,T_{\VAeventno+1}^g)\,,\\
(T_{\VAeventno+1}^b,0,1,0) & \text{ else if }&T_{\VAeventno+1}^b=\text{min}(T_{\VAeventno+1}^c,T_{\VAeventno+1}^b,T_{\VAeventno+1}^g)\,,\\
(T_{\VAeventno+1}^g,0,0,1) & \text{ else if }&T_{\VAeventno+1}^g=\text{min}(T_{\VAeventno+1}^c,T_{\VAeventno+1}^b,T_{\VAeventno+1}^g)\,.
\end{array}\right.
\end{align}
Now the next event time, and the nature of the next event is established, we discuss how each event is simulated.

\paragraph{Collision events in analog simulation.} In the non-analog collision type and non-analog track-length type simulations, a collision event with the plasma is always executed as a scattering event. In the analog simulation, a collision is either an absorption (with probability $\frac{\VAratea(x(T_\VAeventno))}{\VAratet(x(T_\VAeventno))}$) or a scattering event (with the complementary probability $\frac{\VArates(x(T_\VAeventno))}{\VAratet(x(T_\VAeventno))}$). 

In an analog simulation, the nature of this collision determined by a Bernouilli distributed random number
\begin{equation}
a_\VAeventno\sim\mathcal{B}\left(\frac{\VAratea(x(T_\VAeventno))}{\VAratet(x(T_\VAeventno))}\right)\,.\label{eq:iif_modelsim_absorption}
\end{equation}
If $a_\VAeventno$ equals one at a plasma collision, this signifies absorption, if it is zero, a scattering occurs. If the event is an absorption, it is the last event, with event index
\begin{equation}
\VAa{\VAnofevents}=\text{min}(\VAeventno|a_\VAeventno c_\VAeventno=1)\,.\label{eq:iif_modelsim_Kabsorption}
\end{equation}
If it is a scattering event, a new velocity is sampled as
\begin{equation}
{v}_{\VAeventno}^s\sim \VApostcolveldistr({v}|{x}(T_{\VAeventno}))\,.
\end{equation}

\paragraph{Collision events in a non-analog collision type simulation.}In the non-analog collision type simulation, the absorption events are not executed as such, but instead all collisions with the plasma are executed as scattering collisions. To remove the resulting bias, a weight loss is introduced, where the particles lose a fraction $\frac{\VAratea(x(T_\VAeventno))}{\VAratet(x(T_\VAeventno))}$ at every collision with the plasma. One way of interpreting this different simulation type is by considering every particle as representing an infinite number of particles moving together. At every collision a fraction $\frac{\VAratea(x(T_\VAeventno))}{\VAratet(x(T_\VAeventno))}$ of this infinite amount is absorbed and the remaining fraction scatters together and moves further. For non-analog collision type simulations, we therefore update the particle weight as 
\begin{align}
{w}(t)&={\VAdiscreteweight}_\VAeventno \text{ for } t\in[{T}_\VAeventno,{T}_{\VAeventno+1}[,\text{ }\VAeventno\in\{0,\dots,{\VAnofevents}-1\}\,.\label{eq:single_particle_w1}\\
{\VAdiscreteweight}_{\VAeventno+1}&={\VAdiscreteweight}_\VAeventno\left(1-{c}_{\VAeventno+1}\dfrac{\VAratea({{x}}({T}_{\VAeventno+1}))}{\VAratet({{x}}({T}_{\VAeventno+1}))}\right)\,,\label{eq:single_particle_reweigh_nac}\\	
\end{align}

\paragraph{Absorption in non-analog track-length type simulation.} The track-length simulation strategy removes the absorption collisions with the plasma entirely and instead performs absorption continuously along the particle path. This means during every $\text{d}t$, the expected absorbed fraction, $\VAratea(x(t))\text{d}t$ is removed. By adopting the viewpoint of every particle representing an infinite set of particles, this can be interpreted as that during every interval $\text{d}t$ the expected amount of particles is absorbed.
The particle behaviour is described as a function $t\in[0,{T}_\text{end}]\mapsto ({{x}}(t),{{v}}(t),{w}(t))$, where the weight is now described as a continuous function, determined by Equation~\eqref{eq:single_particle_w_natl}: 
\begin{align}
\dfrac{\text{d}{w}(t)}{\text{d}t}=-\VAratea({{x}}(t)){w}(t)\,,\label{eq:single_particle_w_natl}\qquad {\VAdiscreteweight}(0)=1\,.
\end{align}

In practical fusion simulations, the neutral particles of the Monte Carlo simulation usually move against a piecewise constant plasma background. This facilitates sampling of the next event times by transforming the integral in Equation~\eqref{eq:single_particle_scatabsT} to a sum. For both the analog simulation (\texttt{a}) as the non-analog collision type simulation (\texttt{nac}), this is the only simplification with respect to the general plasma model. For the non-analog track-length type simulation, the weight function ${w}(t)$ also simplifies. Since we take grid cell crossings to be events, the plasma background is always constant in between two events. Consequently, with $\VAratea^{{\VAgridcellno}(\VAeventno)}$ the absorption rate in between events $\VAeventno$ and $\VAeventno+1$, it follows from Equation~\eqref{eq:single_particle_w_natl} that the weight function becomes
\begin{equation}
{w}(t)={w}_\VAeventno e^{-\VAratea^{{\VAgridcellno}(\VAeventno)}(t-{T}_{\VAeventno})}\text{ for }t\in[{T}_\VAeventno,{T}_{\VAeventno+1}], k\in[0,\dots,{K}-1].\label{eq:iif_modelsim_natl_pwconst_weight}
\end{equation}
Therefore, the non-analog track-length type simulation can also be performed by a simple weight update at events. With ${d}_\VAeventno=|{\VAdiscreteposeend}_{\VAeventno+1}-{\VAdiscreteposeend}_\VAeventno|$ the distance between events $\VAeventno$ and $\VAeventno+1$, the weight update can be written as
\begin{equation}
{\VAdiscreteweight}_{\VAeventno+1}=e^{-\VAratea\frac{{d}_\VAeventno}{{\VAdiscreteveleend}_\VAeventno}}{\VAdiscreteweight}_\VAeventno\,.
\end{equation}

Table~\ref{tab:overview_weightchanges} provides an overview of how the different simulation types reweigh at the $\VAeventno$-th event when the reaction rates are constant in between events.
\begin{table}[H]
\centering
\begin{tabular}{cccc}
Simulation 						&	Scattering event	&	Boundary hit & Grid cell crossing\\ \hline
\texttt{a} 						&	1				& 1 &1\\
\texttt{nac}	&	$\VAratea/\VAratet$ & 1 &1\\
\texttt{natl} &	$e^{-\frac{\VAratea}{|{\VAdiscreteveleend}_{\VAeventno}|}{d}_\VAeventno}$ & $e^{-\frac{\VAratea}{|{\VAdiscreteveleend}_{\VAeventno}|}{d}_\VAeventno}$ & $e^{-\frac{\VAratea}{|{\VAdiscreteveleend}_{\VAeventno}|}{d}_\VAeventno}$
\end{tabular}
\caption{$\VAdiscreteweight_{\VAeventno+1}/\VAdiscreteweight_{\VAeventno}$, the post-event reweighing factor at event $\VAeventno$ in the different simulations, with $d_\VAeventno=|\VAdiscreteposeend_{\VAeventno+1}-\VAdiscreteposeend_{\VAeventno}|$, and $\VAratea$, respectively $\VAratet$, the constant absorption rate, respectively total reaction rate, between events $\VAeventno$ and $\VAeventno+1$.}
\label{tab:overview_weightchanges}
\end{table}

\paragraph{Boundary hit.} At the boundary, the particle can either be absorbed (with probability $\alpha(x)$ or reflected (with probability $1-\alpha(x)$). To select the type, we again generate a Bernouilli distributed variable:
\begin{equation}
\beta_{\VAeventno}\sim\mathcal{B}(\alpha(\VAdiscreteposeend_{\VAeventno}^b))\,.\label{eq:single_particle_refl}\\
\end{equation}
If $b_{\VAeventno}=1$ and $\beta_{\VAeventno}=1$, the particle leaves the domain, as is expressed by
\begin{align}
\VAnofevents_\text{out}&=\min\left(\VAeventno|b_{\VAeventno}\beta_{\VAeventno}=1\right)\,.\label{eq:single_particle_nout}
\end{align}
Otherwise, if $b_{\VAeventno}=1$ and $\beta_{\VAeventno}=0$, the particle is reflected, which leads to a new deterministic velocity based on its current velocity and the normal to the boundary at the reflection location, via the equation
\begin{equation}
{v}^r_{\VAeventno+1}=\VAreflectionfunction(\VAdiscreteveleend_{\VAeventno})\,,\label{eq:single_particle_reflV}
\end{equation}
with $\VAreflectionfunction({v})\equiv -{v}$, a deterministic function that represents perfect reflection at the wall. 

\paragraph{Grid cell crossing.} At a grid cell crossing, essentially nothing happens, i.e. the particle continues to move with its current velocity.

\paragraph{New velocity.} Depending on the event type, the new velocity differs. This can be described as
\begin{equation}
v_{\VAeventno+1}=\underbrace{(1-a_{\VAeventno+1})c_{\VAeventno+1}v_{\VAeventno+1}^s}_\text{scattering collision}
+\underbrace{(1-\beta_{\VAeventno+1})b_{\VAeventno+1}v_{\VAeventno+1}^r}_\text{boundary reflection}
+\underbrace{(a_{\VAeventno+1}c_{\VAeventno+1}+\beta_{\VAeventno+1}b_{\VAeventno+1}+g_{\VAeventno+1})v_\VAeventno}_\text{absorption, exit at boundary, or grid cell crossing}\,,
\end{equation}
where the velocity at the final event does not change.

\paragraph{End time.}
The particle path ends when it leaves the domain (or at the first absorption event for an analog simulation), hence
\begin{align}
\VAnofevents&=\min(\VAnofevents_\text{a},\VAnofevents_\text{out})\,,\label{eq:iif_modelsim_a_K}\\
T_\text{end}&=T_\VAnofevents=\min\left(T_{\VAnofevents_\text{a}},T_{\VAnofevents_\text{out}}\right)\,.\label{eq:iif_modelsim_a_Tend}
\end{align}

\paragraph{Particle trajectory.} The entire particle trajectory can be determined by the set of events with the variables that encode the event types:
\begin{equation}
\VAparticlepath=\{\VAdiscreteposeend_\VAeventno,\VAdiscreteveleend_\VAeventno,c_\VAeventno,b_\VAeventno,g_\VAeventno,a_\VAeventno,\beta_\VAeventno\}_{\VAeventno=0}^\VAnofevents\,.
\end{equation}
For non-analog collision type and non-analog track-length type simulations the $a_\VAeventno$ have no impact of course.

\section{Source term estimation procedures\label{sec:iif_est}}

The simulation techniques of Section~\ref{sec:iif_modelsim} provide particle trajectories sampling the Boltzmann-BGK equation, Equation~\eqref{eq:iif_modelsim_model_kinetic}. The quantities of interest we are after are the stationary rates of mass, momentum, and energy transfer from the neutral particles to the plasma in each of the $\VAnofgridcells$ grid cells $\{\VAgridcelldefgridcellno\}_{\VAgridcellno=1}^\VAnofgridcells$. These quantities feature as source terms in the plasma equations. To find these source terms, the simulation methods need to be combined with source term estimators, which we introduce in this section.

In their general form, the stationary rates of mass, momentum, and energy transfer can be written as
\begin{equation}
\VAgeneriek{S}^{{\VAgridcellno}}=
\int \left(\VAratea({x})\VAageneriek{
\VAsourceatev}({v})+\VArates({x})\int \VAsgeneriek{\VAsourceatev}({v}\rightarrow{v}')\VApostcolveldistrvpr\text{d}{v}'\right)\VAgridcellinddefgridcellno({x}) \phi_\text{n}({x},{v})\text{d}{v}\text{d}{x}\,,
\label{eq:sources_fromcols_gen_FV}
\end{equation}
with $\VAgridcellinddefgridcellno({x})$ the characteristic function of the grid cell $\VAgridcelldefgridcellno$ and with $\VAageneriek{\VAsourceatev}({v})$ and $\VAsgeneriek{\VAsourceatev}({v}\rightarrow {v}')$ the exchanges due to a particle with velocity $v$ that is absorbed, respectively a particle with velocity $v$ that is scattered into a particle with velocity $v'$. In this series of papers, we study the mass and momentum source extensively, the source contributions for these are
\begin{align}
\VAamass{\VAsourceatev}({v})&=1\,,
& &\VAsmass{\VAsourceatev}({v}\rightarrow {v}')=0\,,\label{eq:sources_sourceatev_mass}\\
\VAamom{\VAsourceatev}({v})&={v}\,,
& &\VAsmom{\VAsourceatev}({v}\rightarrow {v}')={v}-{v}'\,.\label{eq:sources_sourceatev_mom}
\end{align}

A Monte Carlo approximation replaces Equation~\eqref{eq:sources_fromcols_gen_FV} by contributions at events, averaged over $\VAnofparticles$ particle paths $\{\VAparticlepath_\VAparticleno\}_{\VAparticleno=1}^\VAnofparticles$. In general, each of the $\VAnofevents^\VAparticleno+1$ events of a particle path can result in a \emph{score} for the source term estimator in the grid cell $\VAgridcellno$. Depending on the event type and the estimator type, the scoring can differ. We denote the score at an event as $\VAsourceatev^{\VAgridcellno,\VAplaceholderest}_{\text{e},\VAplaceholderquant}$, with the superscript $\VAplaceholderest$ serving as a place-holder for the estimator type, the $\VAplaceholderev$ for the event type, and the subscript $\VAplaceholderquant$ for the estimated quantity. The resulting Monte Carlo estimator then reads
\begin{multline}
\VAMCapprox{S}^{{\VAgridcellno},\VAplaceholderest}_\VAplaceholderquant=\frac{1}{\VAnofparticles}\sum_{\VAparticleno=1}^\VAnofparticles\sum_{\VAeventno=0}^{\VAnofevents^\VAparticleno} \left(\vphantom{\VAsourceatev^{\VAgridcellno,\VAplaceholderest}_k}\right.\underbrace{\VAageneriek{\VAsourceatev}^{\VAgridcellno,\VAplaceholderest}(\VAparticlepath^\VAparticleno,\VAeventno)b_\VAeventno^\VAparticleno a_\VAeventno^\VAparticleno}_\text{absorption}+\underbrace{\VAsgeneriek{\VAsourceatev}^{\VAgridcellno,\VAplaceholderest}(\VAparticlepath^\VAparticleno,\VAeventno)c_\VAeventno^\VAparticleno(1-a_\VAeventno^\VAparticleno)}_\text{scattering}\left.\vphantom{\VAsourceatev^{\VAgridcellno,\VAplaceholderest_k}}\right.\\
\left.+\underbrace{\VAegeneriek{\VAsourceatev}^{\VAgridcellno,\VAplaceholderest}(\VAparticlepath^\VAparticleno,\VAeventno)b_\VAeventno^\VAparticleno(1-\beta_\VAeventno^\VAparticleno)}_\text{boundary absorption}+\underbrace{\VArgeneriek{\VAsourceatev}^{\VAgridcellno,\VAplaceholderest}(\VAparticlepath^\VAparticleno,\VAeventno)b_\VAeventno^\VAparticleno\beta_\VAeventno^\VAparticleno}_\text{boundary reflection}+\underbrace{\VAggeneriek{\VAsourceatev}^{\VAgridcellno,\VAplaceholderest}(\VAparticlepath^\VAparticleno,\VAeventno)g_\VAeventno^\VAparticleno}_\text{grid cell crossing}\right)\,,
\label{eq:sources_frompaths_general_an}
\end{multline}
where all the symbols are introduced in Section~\ref{sec:iif_modelsim}.

The Monte Carlo approximation in Equation~\eqref{eq:sources_frompaths_general_an} is for analog particle paths. For non-analog particle paths, the scores should be multiplied by the appropriate weights. Most of the different estimators can be used for the three types of simulations we provided in Section~\ref{sec:iif_modelsim}. Combined, a simulation type and source term estimator type form a \emph{source term estimation procedure}.

In the remainder of this section, we will present four types of estimators: the analog estimator (Section~\ref{subsec:iif_est_a}), collision estimator (Section~\ref{subsec:ii_est_c}), next-event estimator (Section~\ref{subsec:ii_est_ne}), and track-length estimator (Section~\ref{subsec:ii_est_tl}). For each of these estimators, we will present the source contributions at every event for the analog simulation, discuss how to apply it to the two non-analog simulation types, and present the fundamental cases of estimating the expected number of absorption events, with
\begin{equation}
\VAaabs{\VAsourceatev}({v})=1\quad \text{and}\quad\VAsabs{\VAsourceatev}({v}\rightarrow {v}')=0\,,\label{eq:iif_est_nofabssources}
\end{equation} and the expected number of scattering events, with
\begin{equation}
\VAascat{\VAsourceatev}({v})=0\quad\text{and}\quad\VAsscat{\VAsourceatev}({v}\rightarrow {v}')=1\,,\label{eq:iif_est_nofscatsources}
\end{equation}
for the practically relevant case of a piecewise constant plasma. All other potential quantities of interest can be easily found by multiplying with the actual $\VAageneriek{\VAsourceatev}({v})$ and $\VAsgeneriek{\VAsourceatev}({v}\rightarrow {v}')$ given in Equation~\eqref{eq:sources_sourceatev_mass} for mass and~\eqref{eq:sources_sourceatev_mom} for momentum. Then, in Section~\ref{subsec:iif_est_overview}, we present an overview of how to execute the resulting source term estimation procedures in the practically relevant case of a piecewise constant plasma.

As for the simulation techniques, we will restrict ourselves to a short review of the estimators and not derive them in full, which has been done in~\cite[Chapter 2]{mortier2020phdthesis}.

\subsection{Analog estimator\label{subsec:iif_est_a}}

The most straightforward estimator for $\VAgeneriek{S}^{{\VAgridcellno}}$ is the analog estimator. An analog estimator scores according to the physical event that is being simulated. In an analog simulation, an analog estimator scores 
\begin{equation}
\VAageneriek{\VAsourceatev}^{\VAgridcellno,\texttt{a}}(\VAparticlepath^\VAparticleno,\VAeventno)=\VAageneriek{\VAsourceatev}(v_\VAeventno)\VAgridcellinddefgridcellno({x}_\VAeventno^\VAparticleno)
\label{eq:ii_est_a_scoreatev_abs}
\end{equation}
at each absorption collision and 
\begin{equation}
\VAsgeneriek{\VAsourceatev}^{\VAgridcellno,\texttt{a}}(\VAparticlepath^\VAparticleno,\VAeventno)=\VAsgeneriek{\VAsourceatev}(v_\VAeventno\rightarrow v_{\VAeventno+1})\VAgridcellinddefgridcellno({x}_\VAeventno^\VAparticleno)
\label{eq:ii_est_a_scoreatev_scat}
\end{equation}
at each scattering collision. The superscript \texttt{a} denotes the analog estimator. At other events, there is no physical exchange with the plasma, hence $\VAegeneriek{\VAsourceatev}^{\VAgridcellno,\texttt{a}}(\VAparticlepath^\VAparticleno,\VAeventno)=\VArgeneriek{\VAsourceatev}^{\VAgridcellno,\texttt{a}}(\VAparticlepath^\VAparticleno,\VAeventno)=\VAggeneriek{\VAsourceatev}^{\VAgridcellno,\texttt{a}}(\VAparticlepath^\VAparticleno,\VAeventno)=0$.

\paragraph{Extension to non-analog simulations.} Non-analog simulations are not relevant for an analog estimator, since the events lose their physical meaning.

\paragraph{Fundamental cases.} The source term estimator for the expected number of absorption events is readily found by using Equation~\eqref{eq:iif_est_nofabssources} in Equations~\eqref{eq:ii_est_a_scoreatev_abs} and~\eqref{eq:ii_est_a_scoreatev_scat}, giving $\VAaabs{\VAsourceatev}(\VAparticlepath,\VAeventno)=1$ and $\VAsabs{\VAsourceatev}(\VAparticlepath,\VAeventno)=0$.
Similarly for the expected number of scattering events with Equation~\eqref{eq:iif_est_nofscatsources}, we find $\VAascat{\VAsourceatev}^{\VAgridcellno,\texttt{a}}(\VAparticlepath,\VAeventno)=0$ and $\VAsscat{\VAsourceatev}^{\VAgridcellno,\texttt{a}}(\VAparticlepath,\VAeventno)=\VAgridcellinddefgridcellno(\VAdiscreteposeend_\VAeventno)$. The plasma being piecewise-constant has no impact on what has to be scored here. We denote the analog estimator for the number of absorption events in an analog simulation by \texttt{a\_a\_abs} and for the number of scattering events by \texttt{a\_a\_sc}.

\subsection{Collision estimator\label{subsec:ii_est_c}}

A first alternative estimator does not distinguish between absorption and scattering collisions, but samples the expected contribution due to a collision,
\begin{equation}
\VAcgeneriek{\VAsourceatev}({x},{v})=\frac{\VAratea({x})}{\VAratet({x})}\VAageneriek{\VAsourceatev}({v})+\frac{\VArates({x})}{\VAratet({x})}\int \VAsgeneriek{\VAsourceatev}({v}\rightarrow{v}')\VApostcolveldistrvpr\text{d}{v}'
\label{eq:ii_est_c_sourceatcol}
\end{equation}
at every collision with the plasma in the considered grid cell, or $\VAcgeneriek{\VAsourceatev}({x},{v})\VAgridcellinddefgridcellno(\VAdiscreteposeend_\VAeventno^\VAparticleno)$ in general. This estimator is called the collision estimator and we use \texttt{c} to refer to it.

Note that the integral over the post-collisional velocity ${v}'$ in Equation~\eqref{eq:ii_est_c_sourceatcol} is not difficult to compute: since $\VAsgeneriek{\VAsourceatev}({v}\rightarrow{v}')$ is a simple polynomial in ${v}'$, the solution of the inner integral is a combination of moments of $\VApostcolveldistr^{\,{\VAgridcellno}}({v}')$, which are known from the plasma simulation. For instance, if the momentum is the estimated quantity, $\VAamom{\VAsourceatev}({v})={v}$ and $\VAsmom{\VAsourceatev}({v}\rightarrow {v}')={v}-{v}'$ and Equation~\eqref{eq:ii_est_c_sourceatcol} becomes
\begin{equation}
\VAcmom{\VAsourceatev}({x},{v})=\frac{\VAratea({x})}{\VAratet({x})}{v}+\frac{\VArates({x})}{\VAratet({x})}({v}-\VAplasmaspeed({x}))\,,\label{eq:ii_est_c_sourceatcol_mom}
\end{equation} 
with $\VAplasmaspeed({x})$ the expected post-collisional velocity at position ${x}$, a quantity that is computed by the plasma simulation.

With $\VAcgeneriek{\VAsourceatev}({x},{v})$ as defined in Equation~\eqref{eq:ii_est_c_sourceatcol}, the scores for a collision estimator can be written as
\begin{equation}
\VAageneriek{\VAsourceatev}^{\VAgridcellno,\texttt{c}}(\VAparticlepath,\VAeventno)=\VAsgeneriek{\VAsourceatev}^{\VAgridcellno,\texttt{c}}(\VAparticlepath,\VAeventno)
=\VAcgeneriek{\VAsourceatev}(\VAdiscreteposeend_\VAeventno,\VAdiscreteveleend_\VAeventno)\VAgridcellinddefgridcellno(\VAdiscreteposeend_\VAeventno^\VAparticleno)\,,\label{eq:ii_est_c_scoreatev}
\end{equation}
and $\VAegeneriek{\VAsourceatev}^{\VAgridcellno,\texttt{c}}(\VAparticlepath,\VAeventno)=\VArgeneriek{\VAsourceatev}^{\VAgridcellno,\texttt{c}}(\VAparticlepath,\VAeventno)=\VAggeneriek{\VAsourceatev}^{\VAgridcellno,\texttt{c}}(\VAparticlepath,\VAeventno)=0$.

\paragraph{Extension to non-analog simulations.} For the non-analog collision type simulation, the only change consists in multiplying the score at each collision event $k$ with the weight at the collision, being ${\VAdiscreteweight}_{\VAeventno-1}$. In the non-analog track-length type simulation, the score at each collision event $k$ has to be multiplied by ${\VAdiscreteweight}_\VAeventno\frac{\VAratet({x})}{\VArates({x})}$, where the second factor compensates for the fact that the absorption collisions with the plasma do not occur.

\paragraph{Fundamental cases.} When the plasma is piecewise constant, the reaction rates are constants $\VAratea^\VAgridcellno$, $\VArates^\VAgridcellno$, and $\VAratet^\VAgridcellno$ in each grid cell $\VAgridcellno$ and the post-collisional velocity distribution is a position-independent function $\VApostcolveldistr^{\,{\VAgridcellno}}({v}')$. Then, we can rewrite the score at a collision as
\begin{equation}
\VAcgeneriek{\VAsourceatev}({x},{v})\VAgridcellinddefgridcellno({x})=\VAcgeneriek{\VAsourceatev}^\VAgridcellno({v})\VAgridcellinddefgridcellno({x})=\left(\frac{\VAratea^\VAgridcellno}{\VAratet^\VAgridcellno}\VAageneriek{\VAsourceatev}({v})+\frac{\VArates^\VAgridcellno}{\VAratet^\VAgridcellno}\int \VAsgeneriek{\VAsourceatev}({v}\rightarrow{v}')\VApostcolveldistr^{\,\VAgridcellno}({v}')\text{d}{v}'\right)\VAgridcellinddefgridcellno({x})\,.
\label{eq:ii_est_c_sourceatcol_piecewiseconst}
\end{equation}
When we want to estimate the number of absorption and scattering events, this means $\frac{\VAratea^\VAgridcellno}{\VAratet^\VAgridcellno}$, respectively $\frac{\VArates^\VAgridcellno}{\VAratet^\VAgridcellno}$ has to be scored at every collision. When applying these two estimators to the same simulation, the result only differs by a constant. Consequently, the statistical properties of a collision estimator for the expected number of absorption events and for the expected number of scattering events are identical, and shared by a collision estimator for the expected total number of collisions. We thus restrict ourselves to only considering one of these fundamental cases, namely a collision estimator for the total number of collisions. Since we can apply the collision estimator to all three simulation types, we obtain three estimation procedures: \texttt{a\_c}, \texttt{nac\_c}, and \texttt{natl\_c}.

\subsection{Next-event estimator\label{subsec:ii_est_ne}}

The next-event estimator looks one step further into the future. As the name indicates, the estimation looks as far as the next event and counts at the beginning of every flight path the expected contribution of the event (collision with the plasma, boundary hit, or grid cell crossing). In the general case, this type of estimator requires the computation of an integral at every scoring instant, rendering this estimator very expensive. In the relevant case of piecewise constant plasma, however, this integral has a simple analytical solution.

In general, the expected contribution at $x$ for a particle with velocity $v$ during the infinitesimal time interval $\text{d}\tau$ is
\begin{equation}
\left(\VAratea({x})\VAageneriek{\VAsourceatev}({v})+\VArates({x})\int \VAsgeneriek{\VAsourceatev}({v}\rightarrow{v}')\VApostcolveldistrvpr\text{d}{v}'\right)\text{d}t=\VAratet(x)\VAcgeneriek{\VAsourceatev}({x},{v})\text{d}\tau\,.
\end{equation}
If a particle starts a flight at $(x_\VAeventno,v_\VAeventno)$, it will retain its velocity for the entirety of the flight path and it can thus only reach positions that are on the line
\begin{equation}
\VAdiscreteposeend_\VAeventno+\frac{\VAdiscreteveleend_\VAeventno}{|\VAdiscreteveleend_\VAeventno|}d\,,\qquad 0\leq d\,.
\end{equation}
Since a grid cell crossing represents an event, the maximal distance the particle could potentially travel is the distance to the edge of the current grid,
\begin{equation}
D_\VAeventno^\VAgridcellno=D^\VAgridcellno(\VAdiscreteposeend_\VAeventno,\VAdiscreteveleend_\VAeventno)=\max\left(d|\VAdiscreteposeend_\VAeventno+\frac{\VAdiscreteveleend_\VAeventno}{|\VAdiscreteveleend_\VAeventno|}d\in\VAgridcelldefgridcellno\right)\,,
\label{eq:est_ne_bigD}
\end{equation}
when $\VAgridcellinddefgridcellno(\VAdiscreteposeend_\VAeventno)=1$. For convenience, we assume here that the domain boundaries are also grid cell boundaries.

Hence, all the positions from $\VAdiscreteposeend_\VAeventno$ to $\VAdiscreteposeend_\VAeventno+\frac{\VAdiscreteveleend_\VAeventno}{|\VAdiscreteveleend_\VAeventno|}D_\VAeventno^\VAgridcellno$ can be reached in this flight. However, each position $\VAdiscreteposeend_\VAeventno+\frac{\VAdiscreteveleend_\VAeventno}{|\VAdiscreteveleend_\VAeventno|}d$, $0\leq d\leq D_\VAeventno^\VAgridcellno$, is only reached by the fraction of particles that did not collide with the plasma beforehand, which is an expected fraction
\begin{equation}
e^{-\int_0^d\VAratet\left(\VAdiscreteposeend_\VAeventno+\frac{\VAdiscreteveleend_\VAeventno}{|\VAdiscreteveleend_\VAeventno|}\ell\right)\frac{\text{d}\ell}{|\VAdiscreteveleend_\VAeventno|}}\,
\end{equation}
in analog simulations.

Combining the above information leads to the following expected contribution in the grid cell $\VAgridcelldefgridcellno$ from $(\VAdiscreteposeend_\VAeventno,\VAdiscreteveleend_\VAeventno)$ until the next event:
\begin{multline}
\VAsgeneriek{\VAsourceatev}^{\VAgridcellno,\texttt{ne}}(\VAparticlepath,\VAeventno)=\VArgeneriek{\VAsourceatev}^{\VAgridcellno,\texttt{ne}}(\VAparticlepath,\VAeventno)=\VAggeneriek{\VAsourceatev}^{\VAgridcellno,\texttt{ne}}(\VAparticlepath,\VAeventno)\\
=\VAgridcellinddefgridcellno(\VAdiscreteposeend_\VAeventno)\int_0^{D_\VAeventno^\VAgridcellno}\VAratet\left(\VAdiscreteposeend_\VAeventno+\frac{\VAdiscreteveleend_\VAeventno}{|\VAdiscreteveleend_\VAeventno|}d\right)\VAcgeneriek{\VAsourceatev}\left(\VAdiscreteposeend_\VAeventno+\frac{\VAdiscreteveleend_\VAeventno}{|\VAdiscreteveleend_\VAeventno|}d,\VAdiscreteveleend_\VAeventno\right)e^{-\int_0^d\VAratet\left(\VAdiscreteposeend_\VAeventno+\frac{\VAdiscreteveleend_\VAeventno}{|\VAdiscreteveleend_\VAeventno|}\ell\right)\frac{\text{d}\ell}{|\VAdiscreteveleend_\VAeventno|}}\frac{\text{d}d}{|\VAdiscreteveleend_\VAeventno|}\,\label{eq:iif_est_ne_scorewithint}
\end{multline}
and $\VAageneriek{\VAsourceatev}^{\VAgridcellno,\texttt{ne}}(\VAparticlepath,\VAeventno)=\VAegeneriek{\VAsourceatev}^{\VAgridcellno,\texttt{ne}}(\VAparticlepath,\VAeventno)=0$.

The remaining integral in Equation~\eqref{eq:iif_est_ne_scorewithint} becomes tractable when the plasma background is of a simple form. For instance, in the practically relevant case with piecewise constant plasma, Equation~\eqref{eq:iif_est_ne_scorewithint} simplifies to
\begin{equation}
\VAsgeneriek{\VAsourceatev}^{\VAgridcellno,\texttt{ne}}(\VAparticlepath,\VAeventno)=\VArgeneriek{\VAsourceatev}^{\VAgridcellno,\texttt{ne}}(\VAparticlepath,\VAeventno)=\VAggeneriek{\VAsourceatev}^{\VAgridcellno,\texttt{ne}}(\VAparticlepath,\VAeventno)=
\VAgridcellinddefgridcellno(\VAdiscreteposeend_\VAeventno)\VAcgeneriek{\VAsourceatev}^\VAgridcellno(\VAdiscreteveleend_\VAeventno)\left(1-e^{-\frac{\VAratet^\VAgridcellno}{|\VAdiscreteveleend_\VAeventno|}D_\VAeventno^\VAgridcellno}\right)\,,\label{eq:iif_est_ne_score_const}
\end{equation}
with $\VAcgeneriek{\VAsourceatev}^\VAgridcellno(\VAdiscreteveleend_\VAeventno)$ as in Equation~\eqref{eq:ii_est_c_sourceatcol_piecewiseconst}.

\paragraph{Extension to non-analog simulations.} In an \texttt{nac} or \texttt{natl} simulation, the only difference lies in applying the correct weight at all the scattering collisions, reflections, and grid cell crossings, namely ${\VAdiscreteweight}_\VAeventno$.

\paragraph{Fundamental cases.} In the two fundamental cases of estimating the expected number of absorption and scattering events for a piecewise-constant plasma, we can again restrict ourselves to only considering estimating the expected total number of collisions (with $\VAcgeneriek{\VAsourceatev}^\VAgridcellno(\VAdiscreteveleend_\VAeventno)=1$), since they only differ by a constant factor. Applying a next-event estimator for the total number of collisions to the three simulation types, leads to three different estimation procedures: \texttt{a\_ne}, \texttt{nac\_ne}, and \texttt{natl\_ne}.

\subsection{Track-length estimator\label{subsec:ii_est_tl}}

In contrast to the three previous estimator types, the track-length estimator only works if the plasma background has a simple shape. This is true for the piecewise constant plasma background, we will consider here.

The track-length estimator then scores
\begin{equation}
\VAsgeneriek{\VAsourceatev}^{\VAgridcellno,\texttt{tl}}(\VAparticlepath,\VAeventno)=
\VArgeneriek{\VAsourceatev}^{\VAgridcellno,\texttt{tl}}(\VAparticlepath,\VAeventno)=
\VAggeneriek{\VAsourceatev}^{\VAgridcellno,\texttt{tl}}(\VAparticlepath,\VAeventno)\\
=\VAgridcellinddefgridcellno(\VAdiscreteposeend_\VAeventno)\VAgridcellinddefgridcellno(\VAdiscreteposeend_{\VAeventno+1})\VAcgeneriek{\VAsourceatev}^\VAgridcellno(\VAdiscreteveleend_\VAeventno)\frac{\VAratet^\VAgridcellno}{|\VAdiscreteveleend_\VAeventno|}|\VAdiscreteposeend_{\VAeventno+1}-\VAdiscreteposeend_\VAeventno|\,,\label{eq:ii_est_tl_sourceatev}
\end{equation}
and
\begin{equation}
\VAageneriek{\VAsourceatev}^{\VAgridcellno,\texttt{tl}}(\VAparticlepath,\VAeventno)=\VAegeneriek{\VAsourceatev}^{\VAgridcellno,\texttt{tl}}(\VAparticlepath,\VAeventno)=0\,.
\end{equation}
The name track-length arises from the fact that the score is related to the travelled length $|\VAdiscreteposeend_{\VAeventno+1}-\VAdiscreteposeend_{\VAeventno}|$. The factor by which the travelled length is multiplied,
\begin{equation}
\VAcgeneriek{\VAsourceatev}^\VAgridcellno(\VAdiscreteveleend_\VAeventno)\frac{\VAratet^\VAgridcellno}{|\VAdiscreteveleend_\VAeventno|}\,,\label{eq:iif_est_tl_expectedscoreperlength}
\end{equation}
thus expresses the expected contributed to the source by an infinitesimal travelled length.

Using the probability density function for the travelled length $d_\VAeventno=|\VAdiscreteposeend_{\VAeventno+1}-\VAdiscreteposeend_{\VAeventno}|$,
\begin{equation}
\text{P}(d_\VAeventno|x_\VAeventno,v_\VAeventno)=\left\{\begin{array}{ll}
\VAratet^\VAgridcellno e^{-\VAratet^\VAgridcellno \frac{d_\VAeventno}{|\VAdiscreteveleend_\VAeventno|}}\text{d}d_\VAeventno & \text{if }0\leq d_\VAeventno<D_\VAeventno^\VAgridcellno\,, \\
e^{-\VAratet\frac{D_\VAeventno^\VAgridcellno}{|v_\VAeventno|}} & \text{if }d_\VAeventno=D_\VAeventno^\VAgridcellno\,,
\end{array}\right.
\end{equation}
the estimator from Equation~\eqref{eq:ii_est_tl_sourceatev} can be seen to be unbiased by taking the expected value over all possible distances $d_\VAeventno$:
\begin{align}
\mathbb{E}\left[\frac{\VAratet^\VAgridcellno}{|\VAdiscreteveleend_\VAeventno|}d_\VAeventno\right]&=\int_0^{D_\VAeventno^\VAgridcellno}\frac{\VAratet^\VAgridcellno}{|\VAdiscreteveleend_\VAeventno|}d \VAratet^\VAgridcellno e^{-\VAratet^\VAgridcellno\frac{d}{|\VAdiscreteveleend_\VAeventno|}}\text{d}d+\frac{\VAratet^\VAgridcellno}{|\VAdiscreteveleend_\VAeventno|}D_\VAeventno^\VAgridcellno e^{-\VAratet\frac{D_\VAeventno^\VAgridcellno}{|\VAdiscreteveleend_\VAeventno|}}\,,\\
&=1-\left(1+\VAratet^\VAgridcellno\frac{D_\VAeventno^\VAgridcellno}{|\VAdiscreteveleend_\VAeventno|}\right)e^{-\VAratet^\VAgridcellno\frac{D_\VAeventno^\VAgridcellno}{|\VAdiscreteveleend_\VAeventno|}}+\frac{\VAratet^\VAgridcellno}{|\VAdiscreteveleend_\VAeventno|}D_\VAeventno^\VAgridcellno e^{-\VAratet\frac{D_\VAeventno^\VAgridcellno}{|\VAdiscreteveleend_\VAeventno|}}\,,\\
&=1-e^{-\VAratet^\VAgridcellno\frac{D_\VAeventno^\VAgridcellno}{|\VAdiscreteveleend_\VAeventno|}}\,,
\end{align}
which equals the corresponding factor in Equation~\eqref{eq:iif_est_ne_score_const}, proving unbiasedness.

\paragraph{Extension to non-analog simulations.} For a non-analog estimation procedure, the time spent in a grid cell should be multiplied by its corresponding weight. In an \texttt{nac} simulation, this weight is constant during a single flight path and equals the weight at the beginning. In an \texttt{natl} simulation, the weight changes continuously along the path. When the plasma is constant during that flight, the weight during the time interval $[{T}_\VAeventno,{T}_{\VAeventno+1}]$ is expressed, according to Equation~\eqref{eq:iif_modelsim_natl_pwconst_weight}, as
\begin{equation}
{w}_\VAeventno e^{-\VAratea^\VAgridcellno\frac{\ell}{|{\VAdiscreteveleend}_\VAeventno|}}\,,
\end{equation}
with $\ell\in[0,{d}_\VAeventno]$ and ${d}_\VAeventno=|{\VAdiscreteposeend}_{\VAeventno+1}-{\VAdiscreteposeend}_\VAeventno|$. Using this correct weight during every infinitesimal travelled length $\text{d}\ell$, gives
\begin{equation}
\int_0^{{d}_\VAeventno}{w}_\VAeventno e^{-\VAratea^\VAgridcellno\frac{\ell}{|{\VAdiscreteveleend}_\VAeventno|}}\text{d}\ell={w}_\VAeventno\frac{|{\VAdiscreteveleend}_\VAeventno|}{\VAratea^\VAgridcellno}\left(1-e^{-\VAratea^\VAgridcellno\frac{\ell}{|{\VAdiscreteveleend}_\VAeventno|}}\right)
\end{equation}
as a weighted travelled distance. By multiplying with the expected score per travelled length from Equation~\eqref{eq:iif_est_tl_expectedscoreperlength}, we find the correct score for an \texttt{natl\_tl} estimation procedure to be
\begin{equation}
\VAsgeneriek{{\VAsourceatev}}^{\VAgridcellno,\texttt{tl}}(\VAparticlepath,\VAeventno)=\VArgeneriek{{\VAsourceatev}}^{\VAgridcellno,\texttt{tl}}(\VAparticlepath,\VAeventno)=\VAggeneriek{{\VAsourceatev}}^{\VAgridcellno,\texttt{tl}}(\VAparticlepath,\VAeventno)=
\VAgridcellinddefgridcellno({\VAdiscreteposeend}_\VAeventno)\VAgridcellinddefgridcellno({\VAdiscreteposeend}_{\VAeventno+1})\VAcgeneriek{\VAsourceatev}^\VAgridcellno(\VAdiscreteveleend_\VAeventno)\VAnatl{\VAdiscreteweight}_\VAeventno\frac{\VAratet^\VAgridcellno}{\VAratea^\VAgridcellno}\!\!\left(\!1-e^{-\frac{\VAratea^\VAgridcellno|{\VAdiscreteposeend}_{\VAeventno+1}-{\VAdiscreteposeend}_\VAeventno|}{|{\VAdiscreteveleend}_{\VAeventno-1}|}}\right)
\label{eq:ii_est_tl_sourceatev_natl}\,.
\end{equation}
and
\begin{equation}
\VAageneriek{{\VAsourceatev}}^{\VAgridcellno,\texttt{tl}}(\VAparticlepath,\VAeventno)=\VAegeneriek{{\VAsourceatev}}^{\VAgridcellno,\texttt{tl}}(\VAparticlepath,\VAeventno)=0\,.
\end{equation}

\paragraph{Fundamental cases.} As for the collision and next-event estimators, changing the estimated quantity namely only modifies the factor $\VAcgeneriek{\VAsourceatev}^\VAgridcellno(\VAdiscreteveleend_\VAeventno)$ in Equation~\eqref{eq:ii_est_tl_sourceatev} (or similar equations for \texttt{nac} and \texttt{natl}). The track-length estimators for the expected number of absorption and scattering events thus have identical statistical properties as a track-length estimator for the expected total number of collisions. We consequently only have to consider three different estimation procedures for the expected total number of collisions: \texttt{a\_tl}, \texttt{nac\_tl}, and \texttt{natl\_tl}.

\subsection{Overview of the estimation procedures applied to the fundamental cases\label{subsec:iif_est_overview}}

In the previous four sections we have discussed the estimation procedures that are studied in this paper series. We have spent specific attention to the practically relevant case of a piecewise constant plasma and the estimation of the number of absorption and the number of scattering collisions, since these fundamental cases lead to the quantities of interest by combining them appropriately with $\VAageneriek{
\VAsourceatev}({v})$ and $\VAsgeneriek{\VAsourceatev}({v}\rightarrow{v}')$. We found in the previous sections that these fundamental cases result in eleven estimation procedures with different statistical behaviour.

As we have discussed in Sections~\ref{subsec:ii_est_c}--\ref{subsec:ii_est_tl}, the statistical properties of estimating the expected number of absorption and scattering events for collision, next-event, and track-length estimators are identical to each other and to estimators for the expected total number of collisions. For that reason, we only consider the latter. For the analog estimator however, we do consider the estimator for the expected number of absorption events and for the expected number of scattering events separately, since they have different statistical properties.

Table~\ref{tab:overview_scores} provides an overview of what should be scored at every event by each of these fundamental estimation procedures. The estimation procedures are applied to estimate the expected total number of collisions, except for the analog estimators, which estimate the total number of absorption collisions (\texttt{a\_a\_abs}), respectively the total number of scattering collisions (\texttt{a\_a\_sc}). Together with Table~\ref{tab:overview_weightchanges}, the table of this section provides a full overview of how to practically implement these estimation procedures.

\begin{table}[H]
\centering
\begin{tabular}{ccccc}
Estimator	& $\VAa{\VAsourceatev}^\VAgridcellno$	& $\VAs{\VAsourceatev}^\VAgridcellno$ & $\VAr{\VAsourceatev}^\VAgridcellno=\VAg{\VAsourceatev}^\VAgridcellno$ & $\VAe{\VAsourceatev}^\VAgridcellno$ \\
\hline
\texttt{a\_a\_abs}	& 1				& 0				& 0 & 0 \\
\texttt{a\_a\_sc}	& 0				& 1				& 0 & 0 \\
\texttt{a\_c}		& 1 & 1	& 0 & 0  \\
\texttt{nac\_c}		& / & ${\VAdiscreteweight}_{\VAeventno-1}$	& 0 & 0  \\
\texttt{natl\_c}	& / & ${\VAdiscreteweight}_\VAeventno\frac{\VAratet^\VAgridcellno}{\VArates^\VAgridcellno}$	& 0& 0  \\
\texttt{a\_tl}		& 0 & $\frac{\VAratet^\VAgridcellno}{|\VAdiscreteveleend_{\VAeventno}|} d_\VAeventno$	& $\frac{\VAratet^\VAgridcellno}{|\VAdiscreteveleend_{\VAeventno}|} d_\VAeventno$	&0\\
\texttt{nac\_tl}	& /	& ${\VAdiscreteweight}_{\VAeventno}\frac{\VAratet^\VAgridcellno}{|{\VAdiscreteveleend}_{\VAeventno}|} {d}_\VAeventno$	& ${\VAdiscreteweight}_{\VAeventno}\frac{\VAratet^\VAgridcellno}{|{\VAdiscreteveleend}_{\VAeventno}|} {d}_\VAeventno$	& 0\\
\texttt{natl\_tl}	& / & ${\VAdiscreteweight}_{\VAeventno}\frac{\VAratet^\VAgridcellno}{\VAratea^\VAgridcellno}\left(1-e^{-\frac{\VAratea^\VAgridcellno}{|{\VAdiscreteveleend}_{\VAeventno}|}{d}_\VAeventno}\right)$	& ${\VAdiscreteweight}_{\VAeventno}\frac{\VAratet^\VAgridcellno}{\VAratea^\VAgridcellno}\left(1-e^{-\frac{\VAratea^\VAgridcellno}{|{\VAdiscreteveleend}_{\VAeventno}|}{d}_\VAeventno}\right)$& 0	\\
\texttt{a\_ex}		& 0 & $1-e^{-\frac{\VAratet^\VAgridcellno}{|\VAdiscreteveleend_\VAeventno|}D_\VAeventno^\VAgridcellno}$ & $1-e^{-\frac{\VAratet^\VAgridcellno}{|\VAdiscreteveleend_\VAeventno|}D_\VAeventno^\VAgridcellno}$ & 0 \\
\texttt{nac\_ex}	& /	& ${\VAdiscreteweight}_\VAeventno\left(1-e^{-\frac{\VAratet^\VAgridcellno}{|{\VAdiscreteveleend}_\VAeventno|}{D}_\VAeventno^\VAgridcellno}\right)$ & ${\VAdiscreteweight}_\VAeventno\left(1-e^{-\frac{\VAratet^\VAgridcellno}{|{\VAdiscreteveleend}_\VAeventno|}{D}_\VAeventno^\VAgridcellno}\right)$& 0 \\
\texttt{natl\_ex}	& / &${\VAdiscreteweight}_\VAeventno\left(1-e^{-\frac{\VAratet^\VAgridcellno}{|{\VAdiscreteveleend}_\VAeventno|}{D}_\VAeventno^\VAgridcellno}\right)$ & ${\VAdiscreteweight}_\VAeventno\left(1-e^{-\frac{\VAratet^\VAgridcellno}{|{\VAdiscreteveleend}_\VAeventno|}{D}_\VAeventno^\VAgridcellno}\right)$& 0 
\end{tabular}
\caption{An overview of how the different estimation procedures score at events in grid cell $\VAgridcellno$ with $d_\VAeventno=|\VAdiscreteposeend_{\VAeventno+1}-\VAdiscreteposeend_{\VAeventno}|$ and $D_\VAeventno^\VAgridcellno=\max(d|\VAdiscreteposeend_{\VAeventno}+d\frac{\VAdiscreteveleend}{|\VAdiscreteveleend|}\in\VAgridcelldefgridcellno)$. The actual scores are the factors in this table multiplied by the factor $\VAgridcellinddefgridcellno(\VAdiscreteposeend_\VAeventno)$, and by an additional factor $\VAgridcellinddefgridcellno(\VAdiscreteposeend_{\VAeventno+1})$ for the track-length estimators. These were left out for clarity.}
\label{tab:overview_scores}
\end{table}

\section{Best estimation procedure in a simplified setting\label{sec:iif_bestest_1D0D}}

To estimate the different source terms, each of the eleven different source term estimation procedures introduced in the previous sections can be used. The main question this papers tries to answer is which of these should be selected, given knowledge of the background plasma parameters such as the event rates and post-collisional velocity distribution. To do so, we study the statistical and computational properties of each of the estimation procedures and determine based on these which estimation procedures performs best. This will aid in a substantiated estimation procedure selection given the problem setting, and even enables the selection of different estimation procedures for different regions of the problem domain. This can be beneficial when the background is very heterogeneous or when the background changes during the simulation, for instance due to iterations of the coupled FV/MC system or due to grid refinement.

In this section, we first consider a simplified forward-backward scattering setting, which admits closed sets of ordinary differential equations (ODEs) for the statistical error and the computational cost for the eleven estimation procedures under study. With these results, we will compare the different estimation procedures and indicate the best estimation procedure throughout the parameter space. In Section~\ref{sec:ii_res_1D1D} we will provide several numerical extensions to mass and momentum source term estimation in a more realistic setting.

In Section~\ref{subsec:iif_bestest_1D0D_1D0D}, we first present the simplified setting for which we have analytical results. Then, in Section~\ref{subsec:iif_bestest_1D0D_ii}, we present the invariant imbedding methodology that will yield closed sets of ODEs. By evaluating these ODEs, we find the best estimation procedure depending on the model parameters. The resulting partition of the parameter space based on the best mass source term estimation procedure is presented in Section~\ref{subsec:iif_bestest_1D0D_bestest} both when considering statistical error and computational cost.

\subsection{Simplified 1D0D simulation\label{subsec:iif_bestest_1D0D_1D0D}}

To facilitate an analytical study of the performance of the different estimation procedures, we simplify the one-dimensional model to be spatially homogeneous and to only have forward-backward scattering. Concretely, we restrict the velocity to being $\pm 1$, making the model zero-dimensional in velocity. Since it is still one-dimensional in space, we refer to this model as \emph{1D0D}.

The post-collisional velocity in the {1D0D} setting is completely determined by one parameter value: the probability of going right after a collision, $\VAPr$. This reduces the post-collision velocity distribution to
\begin{equation}
\VApostcolveldistr^\text{1D0D}(v)=\VAPr\delta(v-1)+(1-\VAPr)\delta(v+1)\,,
\end{equation}
with $\delta$ the Dirac-delta. The constant size of the velocity and the space-independence of the rates, result in constant values $\VAcsa=\frac{\VAratea}{|v|}$, $\VAcss=\frac{\VArates}{|v|}$, and $\VAcst=\frac{\VAratet}{|v|}$. We call these quantities the cross-sections and use them as parameters in the sequel.

We use as an initial condition that all particles enter from the left,
\begin{equation}
\VAnsourceveldistr^\text{1D0D}(v)=\delta(v-1)\,.
\end{equation}
Furthermore, we take the probability of being reflected equal to be 0 at each end of the domain, hence $\alpha(r)\equiv 1$.

This reduced problem, has three remaining dimensions to its parameter space: the survival probability at a collision $\VAcss/\VAcst=\VArates/\VAratet$, the dedimensionalized total collision rate $\VAcst\VAdomlength=\frac{\VAratet}{|v|}\VAdomlength$ and the post-collisional parameter $\VAPr$.

\subsection{Invariant imbedding procedure\label{subsec:iif_bestest_1D0D_ii}}

For the simplified 1D0D model described in Section~\ref{subsec:iif_bestest_1D0D_1D0D}, we can derive ODEs for the statistical properties of the different estimation procedures with an invariant imbedding procedure. The invariant imbedding procedure~\cite{chandrasekhar1950radtransii, bellman1975ii} consists of expressing a moment of a quantity, the score for example, for a slab of length $\VAiidomainlengthvar+\Delta\VAiidomainlengthvar$ as a function of moments of quantities for a slab of length $\VAiidomainlengthvar$ and its extension $\Delta \VAiidomainlengthvar$. By taking the limit of $\Delta \VAiidomainlengthvar\rightarrow 0$, an ODE is formed with the domain length as an integration variable. For the other quantities that arise, a similar procedure can be followed until the set of ODEs is closed.

For each of the eleven mass source estimation procedures, this invariant imbedding procedure leads to a system of ODEs with sizes up to 22 for the statistical properties. In each of these ODE systems the parameters $(\VAcsa,\VAcss,\VAPr)$ feature as parameters that can be fixed and the domain length $\VAiidomainlengthvar$ is the integration variable. We note that in this 1D0D problem, the constant size of the velocity has as an additional effect that energy of the particles always remains proportional to their mass. This means that the performance results for mass source estimation that we attain by this method hold for energy source estimation as well.

To illustrate this procedure, we include part of the invariant imbedding procedure for the non-analog collision type track-length estimation procedure (\texttt{nac\_tl}). The full derivation for each mass source estimation procedure can be found in~\cite{mortier2020iiappendix}.

\subsection{Invariant imbedding applied to the track-length estimator on a non-analog collision type simulation\label{subsec:ii2_ii_example}}

In a non-analog collision type simulation (\texttt{nac}), every collision is executed as a scattering event, at which the particle weight is multiplied by the factor $\VAcss/\VAcst$ to keep the simulation unbiased, see Section~\ref{sec:iif_modelsim}. The total factor by which the particle weight changed after the particle passed through a domain of length $\VAiidomainlengthvar$ in an \texttt{nac} simulation, is denoted by $\VAiiWnac(\VAiidomainlengthvar)$. During this passage, multiple collisions might have taken place. A track-length estimator for the expected number of collisions, scores $\VAcst d$, with $d$ the distance travelled and $\VAcst$ the total cross-section, see Section~\ref{subsec:ii_est_tl}. To estimate the expected number of a certain type of collisions, $\VAcst$ is replaced by $\VAcsa$ for absorption, respectively by $\VAcss$ for scattering events. To clarify the difference between the cross-section as the probability of colliding per travelled length and its role in the score, we denote the score for a travelled distance $d$ as $\VAcs d$. The total score by an \texttt{nac\_tl} estimation procedure by a single particle path through a domain of length $\VAiidomainlengthvar$ is denoted by $\VAiitlnac(\VAiidomainlengthvar)$.

Here, we will focus on the second moment of an \texttt{nac\_tl} estimation procedure sample, an indispensable quantity to compute the statistical properties of the estimation procedure. We only consider the contribution by paths that leave and enter the domain from the left. Other outcomes are treated similarly. We denote the outcome by a subscript of two letters, of which the first denotes the place of entry and the second the place of exit. The probability of a path under the condition that it enters and leaves from the left is denoted by $\VAiiPnacll(\VAiidomainlengthvar)$.

\begin{figure}[H]
\centering
\resizebox{.6\columnwidth}{!}{\def\iioutcomefigdist{.3}%
\def\iioutcomefigdepthbeforeturn{3}%
\newcommand{\domain}[2]{
\draw (#1,#2) -- (#1+1,#2) -- (#1+7, #2) -- (#1+7, #2+1) -- (#1,#2+1) -- (#1,#2);
\draw (#1+1,#2) -- (#1+1, #2+1);
}
\newcommand{\turn}[2]{
 \pgfmathparse{#1}
 \edef\la{\pgfmathresult}
 \pgfmathparse{#2}
 \edef\lala{\pgfmathresult}
 \pgfmathparse{\lala-\la}
 \edef\lalala{\pgfmathresult}
 \draw[blue, densely dashed] (1,#1) -- (\iioutcomefigdepthbeforeturn,#1)
   to[out=0, in=0, distance=100] (\iioutcomefigdepthbeforeturn,#2) -- (1,#2);
}%
\newcommand{\scatter}[2]{
 \draw (#1,#2) node[circle,fill=cyan!70!blue!70,scale=.3,draw=blue!60!black!70] {\huge \textbf{s}};
}%
\begin{tikzpicture}[x=1cm, y=-.7cm]
 \foreach\yy in {0,...,4} {
  \draw (-.5,\yy+\iioutcomefigdist*\yy+.5) node {\pgfmathparse{\yy+1}\pgfmathprintnumber{\pgfmathresult}};
  \domain{0}{\yy+\iioutcomefigdist*\yy}
 }%
 \def\iioutcomefigentryheight{.3}%
 \def\iioutcomefigexitheight{.7}%
 \foreach\yy in {0,1,3} {
  \draw[blue] (0, \yy+\iioutcomefigdist*\yy+\iioutcomefigentryheight) -- (1, \yy+\iioutcomefigdist*\yy+\iioutcomefigentryheight);
  \turn{\yy+\iioutcomefigdist*\yy+\iioutcomefigentryheight}{\yy+\iioutcomefigdist*\yy+\iioutcomefigexitheight}
  \draw[blue, -latex] (1, \yy+\iioutcomefigdist*\yy+\iioutcomefigexitheight) -- (0, \yy+\iioutcomefigdist*\yy+\iioutcomefigexitheight);
 }
 \def\iioutcomefigentryheight{.2}%
 \def\iioutcomefigexitheight{.8}%
 \def\yy{4}
 \draw[blue] (0, \yy+\iioutcomefigdist*\yy+\iioutcomefigentryheight) -- (1, \yy+\iioutcomefigdist*\yy+\iioutcomefigentryheight);
 \turn{\yy+\iioutcomefigdist*\yy+\iioutcomefigentryheight}{\yy+\iioutcomefigdist*\yy+\iioutcomefigentryheight*.6+\iioutcomefigexitheight*.4}
 \turn{\yy+\iioutcomefigdist*\yy+\iioutcomefigentryheight*.4+\iioutcomefigexitheight*.6}{\yy+\iioutcomefigdist*\yy+\iioutcomefigexitheight}
 \draw[blue, -latex] (1, \yy+\iioutcomefigdist*\yy+\iioutcomefigexitheight) -- (0, \yy+\iioutcomefigdist*\yy+\iioutcomefigexitheight);
 \draw[blue] (.5,\yy+\iioutcomefigdist*\yy+\iioutcomefigentryheight*.6+\iioutcomefigexitheight*.4) -- (1,\yy+\iioutcomefigdist*\yy+\iioutcomefigentryheight*.6+\iioutcomefigexitheight*.4);
 \draw[blue] (.5,\yy+\iioutcomefigdist*\yy+\iioutcomefigentryheight*.4+\iioutcomefigexitheight*.6) -- (1,\yy+\iioutcomefigdist*\yy+\iioutcomefigentryheight*.4+\iioutcomefigexitheight*.6);
 \scatter{.5}{\yy+\iioutcomefigdist*\yy+\iioutcomefigentryheight*.5+\iioutcomefigexitheight*.5}
 \def\iioutcomefigentryheight{.3}%
 \def\iioutcomefigexitheight{.7}%
 \def\yy{2}
 \draw[blue] (0, \yy+\iioutcomefigdist*\yy+\iioutcomefigentryheight) -- (.5, \yy+\iioutcomefigdist*\yy+\iioutcomefigentryheight);
 \draw[blue, -latex] (.5, \yy+\iioutcomefigdist*\yy+\iioutcomefigentryheight*.5+\iioutcomefigexitheight*.5) -- (0, \yy+\iioutcomefigdist*\yy+\iioutcomefigentryheight*.5+\iioutcomefigexitheight*.5);
 \scatter{.5}{\yy+\iioutcomefigdist*\yy+\iioutcomefigentryheight*.75+\iioutcomefigexitheight*.25}
 \def\yy{1}
 \scatter{.5}{\yy+\iioutcomefigdist*\yy+\iioutcomefigentryheight}
 \def\yy{3}
 \scatter{.5}{\yy+\iioutcomefigdist*\yy+\iioutcomefigexitheight}
 \def\yy{5}
 \path[latex-latex, draw=black] (0,\yy+\iioutcomefigdist*\yy-\iioutcomefigdist+.15) -- (1,\yy+\iioutcomefigdist*\yy-\iioutcomefigdist+.15);
 \path[latex-latex, draw=black] (1,\yy+\iioutcomefigdist*\yy-\iioutcomefigdist+.15) -- (7,\yy+\iioutcomefigdist*\yy-\iioutcomefigdist+.15);
 \draw[] (.5,\yy+\iioutcomefigdist*\yy-\iioutcomefigdist+.5) node {$\Delta x$};
 \draw[] (4,\yy+\iioutcomefigdist*\yy-\iioutcomefigdist+.5) node {$x$};
 \draw[] (-.5,\yy+\iioutcomefigdist*\yy-\iioutcomefigdist+.5) node {$j$};
\end{tikzpicture}%
}
\caption{The possible paths in a domain of length $\VAiidomainlengthvar+\Delta\VAiidomainlengthvar$ that start and end on the left side in a non-analog simulation and have a probability of order at most 1 in $\Delta\VAiidomainlengthvar$. The symbol \protect\newcommand{\scatter}[2]{
 \draw (#1,#2) node[circle,fill=cyan!70!blue!70,scale=.3,draw=blue!60!black!70] {\huge \textbf{s}};
}%
\begin{tikzpicture}
 \scatter{0}{0}
\end{tikzpicture}%
 in the part of length $\Delta\VAiidomainlengthvar$ identifies that a collision took place there. The dashed lines in the $\VAiidomainlengthvar$ part of the domain signify that the behaviour there is irrelevant, as long as the outcome is correct: entering and leaving the $\VAiidomainlengthvar$ part from the left. The figure is copied from~\cite{mortier2020iiappendix}.}
\label{fig:ii_illustration}
\end{figure}
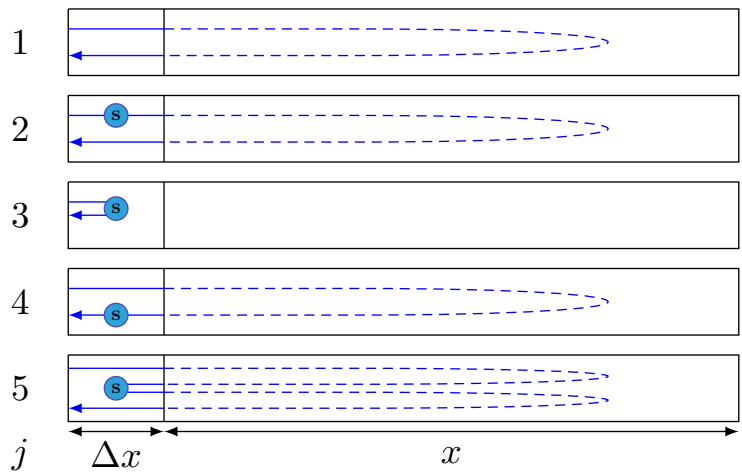

We focus on the quantity $\VAiiPnacll(\VAiidomainlengthvar)\mathbb{E}[\VAiitlnacll^2(\VAiidomainlengthvar)]$ in this example, which is the contribution to the second moment of a track-length estimator in a non-analog collision type simulation in a domain of length $\VAiidomainlengthvar$ by the particle paths that enter and leave from the left. We will express $\VAiiPnacll(\VAiidomainlengthvar+\Delta\VAiidomainlengthvar)\mathbb{E}[\VAiitlnacll^2(\VAiidomainlengthvar+\Delta\VAiidomainlengthvar)]$, as a function of quantities over $\VAiidomainlengthvar$. To do so, we condition the paths that enter and leave the domain of length $\VAiidomainlengthvar+\Delta\VAiidomainlengthvar$ from the left by their behaviour in the part of length $\Delta\VAiidomainlengthvar$. Since our aim is to arrive at an ODE for $\VAiiPnacll(\VAiidomainlengthvar)\mathbb{E}[\VAiitlnacll^2(\VAiidomainlengthvar)]$, we can neglect contributions to $\VAiiPnacll(\VAiidomainlengthvar+\Delta\VAiidomainlengthvar)\mathbb{E}[\VAiitlnacll^2(\VAiidomainlengthvar+\Delta\VAiidomainlengthvar)]$ of order $o(\Delta\VAiidomainlengthvar),\ \Delta\VAiidomainlengthvar\rightarrow0$. This means that we can neglect paths that collide more than once in the $\Delta\VAiidomainlengthvar$ part of the domain, since these have a probability of order $\mathcal{O}(\Delta\VAiidomainlengthvar^2),\ \Delta\VAiidomainlengthvar\rightarrow0$ to occur, since, when travelling a length of $\Delta\VAiidomainlengthvar\rightarrow0$, the probability of colliding is $\VAcst\Delta\VAiidomainlengthvar,\ \Delta\VAiidomainlengthvar\rightarrow0$. When restricting the outcome of the paths to leaving and entering from the left, and to having at most one collision in the $\Delta\VAiidomainlengthvar$ part of the domain, there are five options left of how the particle can act within the domain of length $\Delta\VAiidomainlengthvar$, which are shown in Figure~\ref{fig:ii_illustration}.

\newcommand{\VAcasenumber}{j}

With $\VAiicondition{\VAiiPnac}{ll,\VAcasenumber}(\VAiidomainlengthvar+\Delta\VAiidomainlengthvar)$ and $\mathbb{E}[\VAiicondition{\VAiitlnac}{ll,\VAcasenumber}^2(\VAiidomainlengthvar+\Delta\VAiidomainlengthvar)]$ the probability of the $\VAcasenumber$-th case of Figure~\ref{fig:ii_illustration}, respectively the second moment of $\VAiitlnacll(\VAiidomainlengthvar+\Delta\VAiidomainlengthvar)$ conditioned on the $\VAcasenumber$-th case of Figure~\ref{fig:ii_illustration}, we can write
\begin{equation}
\VAiiPnacll(\VAiidomainlengthvar+\Delta\VAiidomainlengthvar)\mathbb{E}\left[\VAiitlnacll^2(\VAiidomainlengthvar+\Delta\VAiidomainlengthvar)\right]=\sum_{\VAcasenumber=1}^5\VAiicondition{\VAiiPnac}{ll,\VAcasenumber}(\VAiidomainlengthvar+\Delta\VAiidomainlengthvar)\mathbb{E}\left[\VAiicondition{\VAiitlnac}{ll,\VAcasenumber}^2(\VAiidomainlengthvar+\Delta\VAiidomainlengthvar)\right]\,,\label{eq:ii_ii_exa_PllTll_conditioned}
\end{equation}
by Taylor expansion of the exponential. We will now explicitly elaborate $\VAiicondition{\VAiiPnac}{ll,1}(\VAiidomainlengthvar)$ and $\mathbb{E}[\VAiicondition{\VAiitlnac}{ll,1}(\VAiidomainlengthvar)]$, after which we will provide all five terms and show how this leads to an ODE.

In the first case of Figure~\ref{fig:ii_illustration}, the particle does not collide in the $\Delta\VAiidomainlengthvar$ part of the domain after its entry, returns to the $\Delta\VAiidomainlengthvar$ part as the outcome of its passage through the $x$ part of the domain and, again does not collide in the $\Delta\VAiidomainlengthvar$ part, after which it exits the domain of length $\VAiidomainlengthvar+\Delta\VAiidomainlengthvar$. The probability of this case is thus
\begin{equation}
\VAiicondition{\VAiiPnac}{ll,1}(\VAiidomainlengthvar+\Delta\VAiidomainlengthvar)=\left(1-e^{-\VAcst\Delta\VAiidomainlengthvar}\right)\VAiiPnacll(\VAiidomainlengthvar)\left(1-e^{-\VAcst\Delta\VAiidomainlengthvar}\right)\,,
\end{equation}
the probability of not colliding when passing through the $\Delta\VAiidomainlengthvar$ part for the first time, times the probability of returning to the $\Delta\VAiidomainlengthvar$ part, times the probability of not colliding in the second passage through the $\Delta\VAiidomainlengthvar$ part. Since we can neglect all contributions of $o(\Delta\VAiidomainlengthvar),\ \Delta\VAiidomainlengthvar\rightarrow0$, we write
\begin{equation}
\VAiicondition{\VAiiPnac}{ll,1}(\VAiidomainlengthvar+\Delta\VAiidomainlengthvar)=\VAiiPnacll(\VAiidomainlengthvar)(1-2\VAcst\Delta\VAiidomainlengthvar)+\mathcal{O}(\Delta\VAiidomainlengthvar^2),\ \Delta\VAiidomainlengthvar\rightarrow0\,.\label{eq:ii_ii_exa_Pnacll1}
\end{equation}
The score by a particle conditioned on the first case equals
\begin{equation}
\VAiicondition{\VAiitlnac}{ll,1}(\VAiidomainlengthvar+\Delta\VAiidomainlengthvar)=\VAcs\Delta\VAiidomainlengthvar+\VAiitlnacll(\VAiidomainlengthvar)+\VAiiWnacll(\VAiidomainlengthvar)\VAcs\Delta\VAiidomainlengthvar\,,\label{eq:ii2_ii_exa_Etlnacll1_-1}
\end{equation}
with $\Delta\VAiidomainlengthvar$ the travelled length during the first passage through the $\Delta\VAiidomainlengthvar$ part of the domain, and consequently $\VAcs\Delta\VAiidomainlengthvar$ is the corresponding track-length score, where a general cross-section is used which is to be replaced by $\VAcst$ for estimating the total number of collisions. Then, the particle moves through the $\VAiidomainlengthvar$ part of the domain, the point of entry and exit for that passage is the left side, meaning that the score contribution is represented by $\VAiitlnacll(\VAiidomainlengthvar)$ and the weight change during this passage is the factor $\VAiiWnacll(\VAiidomainlengthvar)$. Finally, during the second passage through the $\Delta\VAiidomainlengthvar$ part of the domain, with a reduced weight $\VAiiWnacll(\VAiidomainlengthvar)$, the score contribution is $\VAiiWnacll(\VAiidomainlengthvar)\VAcs\Delta\VAiidomainlengthvar$. From Equation~\eqref{eq:ii2_ii_exa_Etlnacll1_-1}, we readily find the second moment of $\VAiicondition{\VAiitlnac}{ll,1}(\VAiidomainlengthvar+\VAcs\Delta\VAiidomainlengthvar)$, conditioned on the first case, to equal
\begin{multline}
\mathbb{E}\left[(\VAiitlnacll(\VAiidomainlengthvar)+\VAcs\Delta\VAiidomainlengthvar+\VAiiWnacll(\VAiidomainlengthvar)\VAcs\Delta\VAiidomainlengthvar)^2\right]=\mathbb{E}\left[\VAiitlnacll^2(\VAiidomainlengthvar)+2\VAiitlnacll(\VAiidomainlengthvar)\VAcs\Delta\VAiidomainlengthvar+2\VAiitlnacll(\VAiidomainlengthvar)\VAiiWnacll(\VAiidomainlengthvar)\VAcs\Delta\VAiidomainlengthvar\right]\\
+\mathcal{O}(\Delta\VAiidomainlengthvar^2),\ \Delta\VAiidomainlengthvar\rightarrow0\,,\label{eq:ii_ii_Etlnacll1_0}
\end{multline}
\begin{multline}
\hphantom{\mathbb{E}\left[(\VAiitlnacll(\VAiidomainlengthvar)+\VAcs\Delta\VAiidomainlengthvar+\VAiiWnacll(\VAiidomainlengthvar)\VAcs\Delta\VAiidomainlengthvar)^2\right]}=\mathbb{E}\left[\VAiitlnacll^2(\VAiidomainlengthvar)\right]+2\VAcs\Delta\VAiidomainlengthvar\mathbb{E}\left[\VAiitlnacll(\VAiidomainlengthvar)\right]+2\VAcs\Delta\VAiidomainlengthvar\mathbb{E}\left[\VAiitlnacll(\VAiidomainlengthvar)\VAiiWnacll(\VAiidomainlengthvar)\right]\\
+\mathcal{O}(\Delta\VAiidomainlengthvar^2),\ \Delta\VAiidomainlengthvar\rightarrow0\,,\label{eq:ii_ii_Etlnacll1}
\end{multline}

The new quantities on the right hand side of Equation~\eqref{eq:ii_ii_Etlnacll1} ($\mathbb{E}[\VAiitlnacll(\VAiidomainlengthvar)]$, $\mathbb{E}[\VAiitlnacll^2(\VAiidomainlengthvar)]$, and $\mathbb{E}[\VAiitlnacll(\VAiidomainlengthvar)\VAiiWnacll(\VAiidomainlengthvar)]$) are all at most second-order in the path variables ($\VAiitlnacll(\VAiidomainlengthvar)$ and $\VAiiWnacll(\VAiidomainlengthvar)$), as is the quantity, $\mathbb{E}[\VAiitlnacll^2(\VAiidomainlengthvar+\Delta \VAiidomainlengthvar)]$, we are after. This feature of the invariant imbedding procedure is true for all our derivations, and allows to find a closed set of ODEs, since the number of second-order moments we can take is finite and each of them depends on moments of at most second order.

With Equations~\eqref{eq:ii_ii_exa_Pnacll1} and~\eqref{eq:ii_ii_Etlnacll1} we can find the first term of the right hand side of Equation~\eqref{eq:ii_ii_exa_PllTll_conditioned} as
\begin{multline}
\VAiicondition{\VAiiPnac}{ll,1}(\VAiidomainlengthvar+\Delta\VAiidomainlengthvar)\mathbb{E}\left[\VAiicondition{\VAiitlnac}{ll,1}^2(\VAiidomainlengthvar+\Delta\VAiidomainlengthvar)\right]=\VAiiPnacll(\VAiidomainlengthvar)\mathbb{E}\left[\VAiitlnacll^2(\VAiidomainlengthvar)\right]-2\VAcst\VAiiPnacll(\VAiidomainlengthvar)\mathbb{E}\left[\VAiitlnacll^2(\VAiidomainlengthvar)\right]\Delta\VAiidomainlengthvar\\
+2\VAcs\VAiiPnacll(\VAiidomainlengthvar)\mathbb{E}\left[\VAiitlnacll(\VAiidomainlengthvar)\right]\Delta\VAiidomainlengthvar+2\VAcs\VAiiPnacll(\VAiidomainlengthvar)\mathbb{E}\left[\VAiitlnacll(\VAiidomainlengthvar)\VAiiWnacll(\VAiidomainlengthvar)\right]\Delta\VAiidomainlengthvar+\mathcal{O}(\Delta\VAiidomainlengthvar^2),\ \Delta\VAiidomainlengthvar\rightarrow0\,,
\end{multline}

For the other four cases of Figure~\ref{fig:ii_illustration}, a similar method yields
\begin{align}
\VAiicondition{\VAiiPnac}{ll,2}(\VAiidomainlengthvar+\Delta\VAiidomainlengthvar)\mathbb{E}\left[\VAiicondition{\VAiitlnac}{ll,2}^2(\VAiidomainlengthvar+\Delta\VAiidomainlengthvar)\right]&=\VAPr\frac{\VAcss^2}{\VAcst}\VAiiPnacll(\VAiidomainlengthvar)\mathbb{E}\left[\VAiitlnacll^2(\VAiidomainlengthvar)\right]\Delta\VAiidomainlengthvar+\mathcal{O}\left(\Delta\VAiidomainlengthvar^2\right),\ \Delta\VAiidomainlengthvar\rightarrow0\\
\VAiicondition{\VAiiPnac}{ll,3}(\VAiidomainlengthvar+\Delta\VAiidomainlengthvar)\mathbb{E}\left[\VAiicondition{\VAiitlnac}{ll,3}^2(\VAiidomainlengthvar+\Delta\VAiidomainlengthvar)\right]&=0+\mathcal{O}\left(\Delta\VAiidomainlengthvar^2\right),\ \Delta\VAiidomainlengthvar\rightarrow0\\
\VAiicondition{\VAiiPnac}{ll,4}(\VAiidomainlengthvar+\Delta\VAiidomainlengthvar)\mathbb{E}\left[\VAiicondition{\VAiitlnac}{ll,4}^2(\VAiidomainlengthvar+\Delta\VAiidomainlengthvar)\right]&=(1-\VAPr)\VAcst\VAiiPnacll(\VAiidomainlengthvar)\mathbb{E}\left[\VAiitlnacll^2(\VAiidomainlengthvar)\right]\Delta\VAiidomainlengthvar+\mathcal{O}\left(\Delta\VAiidomainlengthvar^2\right),\ \Delta\VAiidomainlengthvar\rightarrow0\\
\begin{split}
\VAiicondition{\VAiiPnac}{ll,5}(\VAiidomainlengthvar+\Delta\VAiidomainlengthvar)\mathbb{E}\left[\VAiicondition{\VAiitlnac}{ll,5}^2(\VAiidomainlengthvar+\Delta\VAiidomainlengthvar)\right]&=\VAPr\VAcst\VAiiPnacll(\VAiidomainlengthvar)\mathbb{E}\left[\VAiitlnacll^2(\VAiidomainlengthvar)\right]\!\Delta\VAiidomainlengthvar\\
&\quad+2\VAPr\VAcss\VAiiPnacll(\VAiidomainlengthvar)\mathbb{E}\left[\VAiitlnacll(\VAiidomainlengthvar)\right]\VAiiPnacll(\VAiidomainlengthvar)\mathbb{E}\left[\VAiiWnacll(\VAiidomainlengthvar)\VAiitlnacll(\VAiidomainlengthvar)\right]\!\Delta\VAiidomainlengthvar\\
&\quad+\VAPr\frac{\VAcss^2}{\VAcst}\VAiiPnacll(\VAiidomainlengthvar)\mathbb{E}\left[\VAiiWnacll^2(\VAiidomainlengthvar)\right]\VAiiPnacll(\VAiidomainlengthvar)\mathbb{E}\left[\VAiitlnacll^2(\VAiidomainlengthvar)\right]\!\Delta\VAiidomainlengthvar+\mathcal{O}\left(\Delta\VAiidomainlengthvar^2\right),\ \Delta\VAiidomainlengthvar\rightarrow0
\end{split}
\end{align}

Substituting these five terms in Equation~\eqref{eq:ii_ii_exa_PllTll_conditioned} gives an expression for $\VAiiPnacll(\VAiidomainlengthvar+\Delta\VAiidomainlengthvar)\mathbb{E}[\VAiitlnacll^2(\VAiidomainlengthvar+\Delta\VAiidomainlengthvar)]$. Dividing that result by $\Delta\VAiidomainlengthvar$ and taking the limit $\Delta\VAiidomainlengthvar\rightarrow0$, results in the ordinary differential equation for $\VAiiPnacll\mathbb{E}[\VAiiEtlnacll^2]$,
\begin{multline}
\frac{\text{d}\left(\VAiiPnacll(\VAiidomainlengthvar)\mathbb{E}[\VAiitlnacll^2(\VAiidomainlengthvar)]\right)}{\text{d}\VAiidomainlengthvar}=-2\VAcst\VAiiPnacll(\VAiidomainlengthvar)\mathbb{E}\left[\VAiitlnacll^2(\VAiidomainlengthvar)\right]
+2\VAcs\VAiiPnacll(\VAiidomainlengthvar)\mathbb{E}\left[\VAiitlnacll(\VAiidomainlengthvar)\right]+2\VAcs\VAiiPnacll(\VAiidomainlengthvar)\mathbb{E}\left[\VAiitlnacll(\VAiidomainlengthvar)\VAiiWnacll(\VAiidomainlengthvar)\right]\\
+\VAPr\frac{\VAcss^2}{\VAcst}\VAiiPnacll(\VAiidomainlengthvar)\mathbb{E}\left[\VAiitlnacll^2(\VAiidomainlengthvar)\right]
+(1-\VAPr)\VAcst\VAiiPnacll(\VAiidomainlengthvar)\mathbb{E}\left[\VAiitlnacll^2(\VAiidomainlengthvar)\right]\\
+\VAPr\VAcst\VAiiPnacll(\VAiidomainlengthvar)\mathbb{E}\left[\VAiitlnacll^2(\VAiidomainlengthvar)\right]
+2\VAPr\VAcss\VAiiPnacll(\VAiidomainlengthvar)\mathbb{E}\left[\VAiitlnacll(\VAiidomainlengthvar)\right]\VAiiPnacll(\VAiidomainlengthvar)\mathbb{E}\left[\VAiiWnacll(\VAiidomainlengthvar)\VAiitlnacll(\VAiidomainlengthvar)\right]\\
+\VAPr\frac{\VAcss^2}{\VAcst}\VAiiPnacll(\VAiidomainlengthvar)\mathbb{E}\left[\VAiiWnacll^2(\VAiidomainlengthvar)\right]\VAiiPnacll(\VAiidomainlengthvar)\mathbb{E}\left[\VAiitlnacll^2(\VAiidomainlengthvar)\right]\,,\label{eq:ii2_exa_result}
\end{multline}
which equals Equation~(66) in~\cite{mortier2020iiappendix}. The initial value to solve this ODE is 0, because the probability to turn in a slab of length zero is zero, hence so is $\VAiiPnacll(0)$ and the travelled length in a slab of length zero also becomes zero, so $\mathbb{E}[\VAiitlnacll^2(0)]=0$.

The same procedure can be applied for each of the new terms in Equation~\eqref{eq:ii2_exa_result}, leading to a closed system of ODEs. The full details of this, and all other derivations for the other estimation procedures, are included in~\cite{mortier2020iiappendix}.

To evaluate the performance of the different estimation procedures for the space spanned by the non-dimensional parameters introduced in Section~\ref{subsec:iif_bestest_1D0D_1D0D} ($\VAcss/\VAcst$, $\VAcst\VAdomlength$, $\VAPr$), the ODE systems are to be evaluated until $\VAdomlength$. In this integration process, the results for smaller values of $\VAdomlength$ are also found. Hence, to find the measure of performance of the estimation procedures on a fine mesh of parameter values, we have to integrate these ODEs for the desired values of $\VAcss/\VAcst$ and $\VAPr$ and the step length and final value of the ODE integration are determined by the desired values of $\VAcst\VAdomlength$.

\subsection{Best estimation procedure for the mass source in the 1wawaD0D model\label{subsec:iif_bestest_1D0D_bestest}}

For each of the estimation procedures described in Sections~\ref{sec:iif_modelsim} and~\ref{sec:iif_est}, the method described in Sections~\ref{subsec:iif_bestest_1D0D_ii} and~\ref{subsec:ii2_ii_example} provides us with the statistical properties for the entire ($\VAcss/\VAcst$,$\VAcst\VAdomlength$,$\VAPr$) parameter space, the execution of which can be found in~\cite{mortier2020iiappendix}. These statistical properties are discussed at length in \cite[Chapter 3]{mortier2020phdthesis}, and some partial results are presented in~\cite{mortier2017invimb}. In this section, we will use the obtained statistical properties to obtain a partition of the ($\VAcss/\VAcst$,$\VAcst\VAdomlength$,$\VAPr$) parameter space based on which estimation procedure performs best.

In Section~\ref{subsubsec:iif_bestest_1D0D_bestest_var}, we will use the statistical error for a fixed number of simulated particles as a measure of performance and in Section~\ref{subsubsec:iif_bestest_1D0D_bestest_cost}, we will consider the computational cost for a given statistical error as a measure of performance. Both measures of performance apply to a different setting: when the number of particles is considered fixed, the goal of the estimation procedure selection is to minimize the statistical error, which is achieved via the partition of Section~\ref{subsubsec:iif_bestest_1D0D_bestest_var}. Alternatively, when the number of simulated particles is not fixed, but a certain statistical error is aimed for, the computational cost for that error is to be minimized, as is achieved in Section~\ref{subsubsec:iif_bestest_1D0D_bestest_cost}.

In Section~\ref{sec:ii_res_1D1D}, we will expand on these results by also considering momentum sources and a more realistic 1D1D setting. In that next section, we only have numerical results.

\subsubsection{Best estimation procedure for the mass source in the 1D0D model based on statistical error\label{subsubsec:iif_bestest_1D0D_bestest_var}}

\begin{figure}[H]\centering
\begin{tabular}{m{.25\textwidth}m{.25\textwidth}m{.25\textwidth}m{.2\textwidth}}
\resizebox{.24\textwidth}{!}{\maakmooieticks{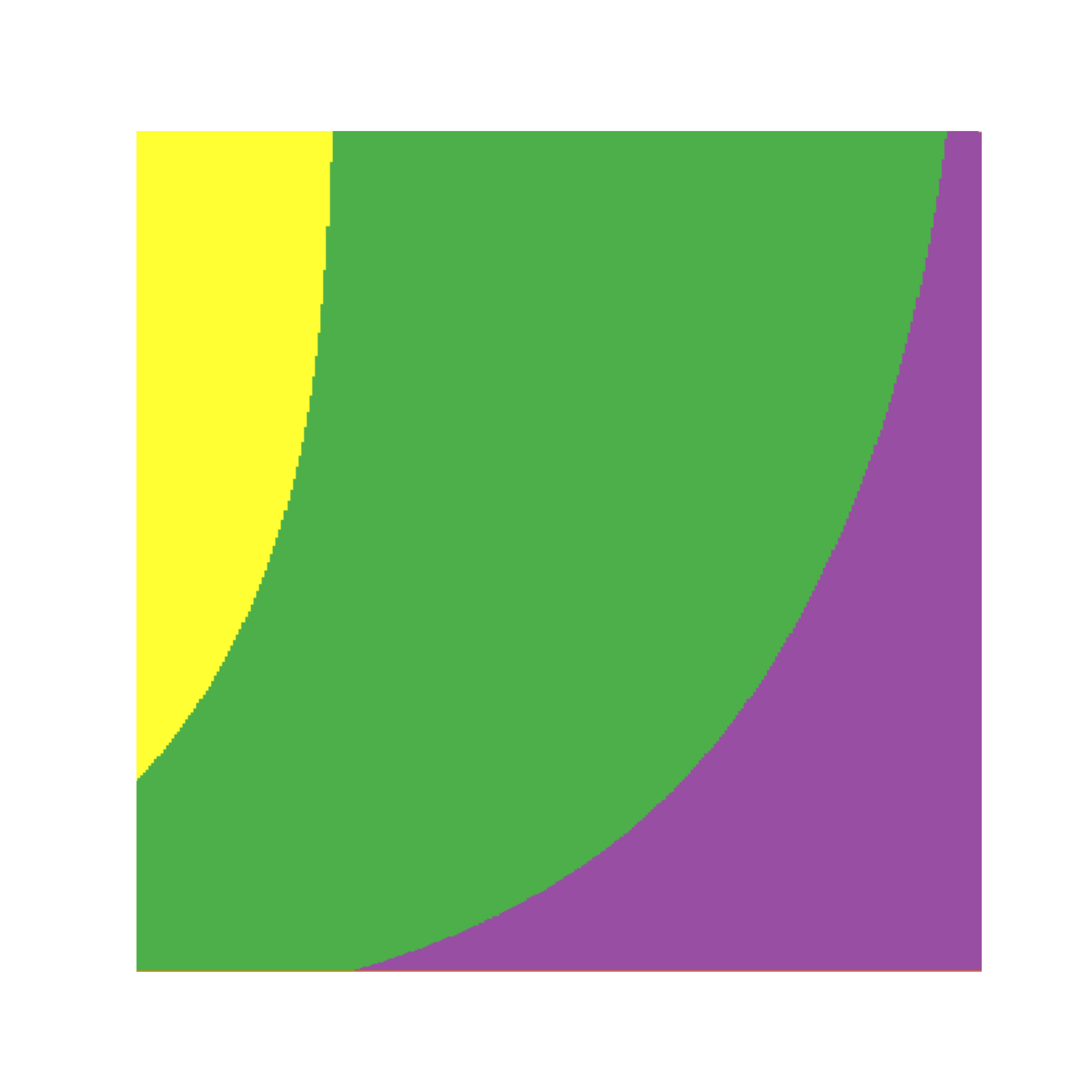}{1}{1}{$\dfrac{\VAcss}{\VAcst}=0.75$}} &
\resizebox{.24\textwidth}{!}{\maakmooieticks{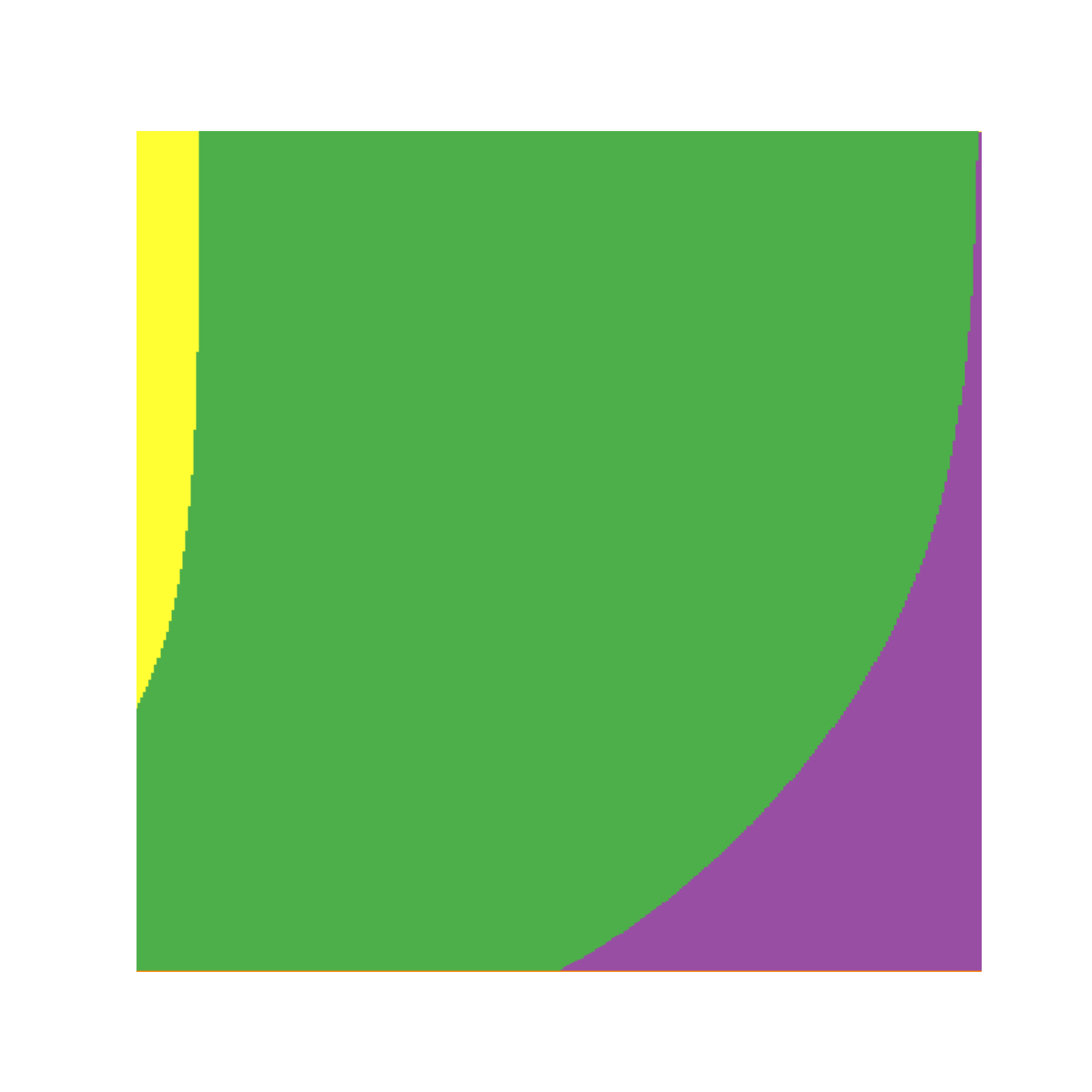}{0}{1}{$\dfrac{\VAcss}{\VAcst}=0.5$}} &
\resizebox{.24\textwidth}{!}{\maakmooieticks{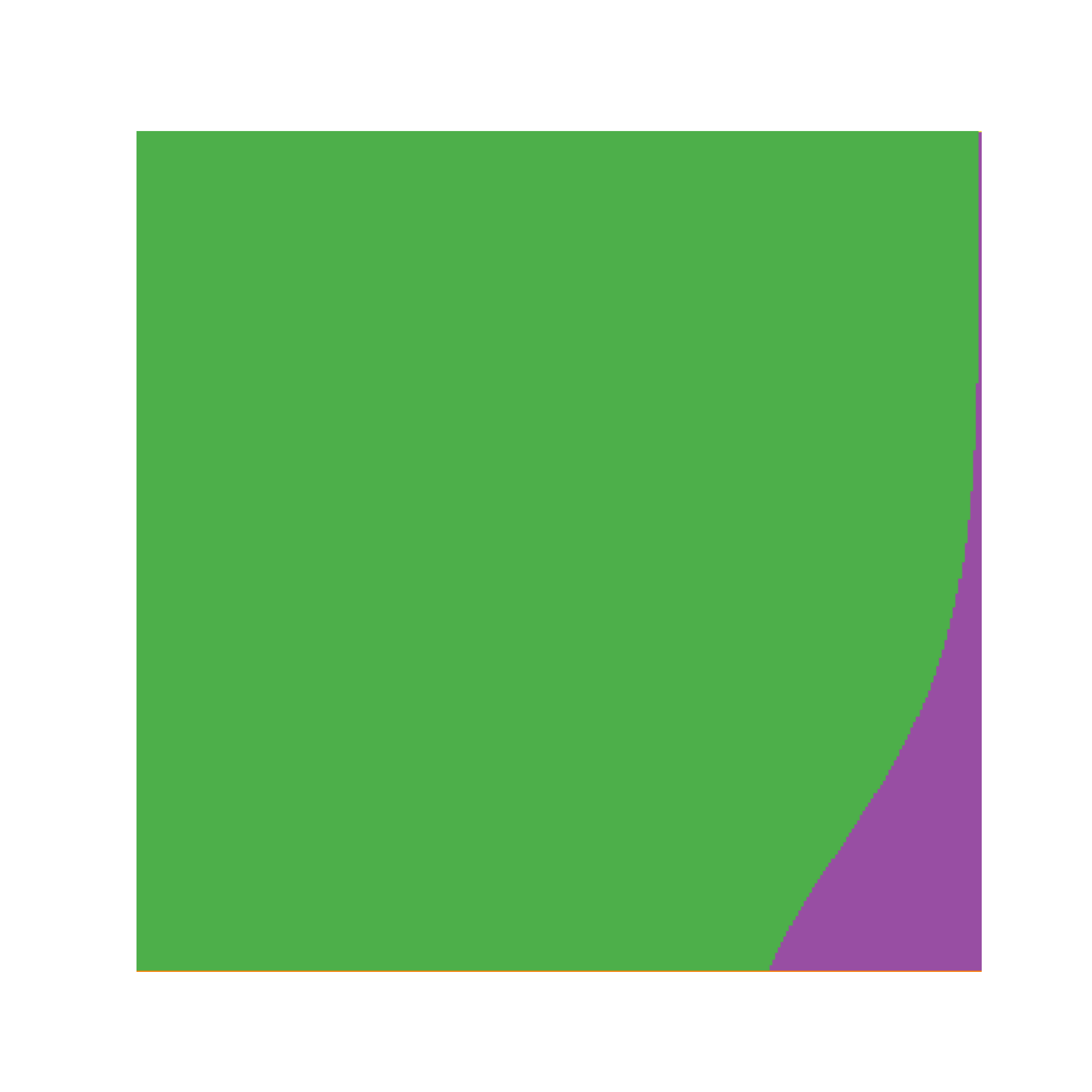}{0}{1}{$\dfrac{\VAcss}{\VAcst}=0.25$}} &
\localcolorlegendoneE
\end{tabular}
\caption{A partition of the parameter domain based on the 1D0D mass source estimation procedure with the lowest statistical error.}\label{fig:best_1D0D_mass_stdv_bestest}
\end{figure}

The first measure of performance we consider is the expected statistical error on the result, given the number of simulated particles. This error is proportional to the standard deviation of the contribution due to a single particle and ODEs for this quantity are derived in~\cite{mortier2020iiappendix}. Figure~\ref{fig:best_1D0D_mass_stdv_bestest} shows the best 1D0D mass source estimation procedure for a large part of the parameter domain. This figure shows there are three competitive estimation procedures: \texttt{nac\_ne}, \texttt{natl\_ne}, and \texttt{natl\_tl}.

At this point, it is relevant to point to two analytical results by Lux: \cite[Theorem 5.22, p. 267]{lux1991MCPT} and~\cite[Theorem 5.24, p. 276]{lux1991MCPT}, which are initially provided in~\cite{lux1978standardvarred}. These state that the next-event estimator for the number of collisions (and thus for the mass source) never has a variance that is larger than the collision estimator for the same simulation, respectively that the collision type survival biasing never has a variance larger than the corresponding analog simulation, granted that the estimator used remains the same. These results are visible in Figure~\ref{fig:best_1D0D_mass_stdv_bestest} by the absence of any collision estimator and of any estimation procedure with an analog simulation. 

The dominating estimation procedure is a next-event estimator with a non-analog collision type simulation, since at low survival probability ($\VAcss/\VAcst$) most of the score is established at the first event in such a simulation --- which is characteristic for a next-event estimator and at higher survival probabilities and not too low $\VAPr$, the higher number of scoring events in a non-analog collision type estimator is beneficial. At low $\VAPr$ and significant collisionality, the neutral that enters from the left, will nearly always leave the domain to the left. The main variation arises from when exactly this occurs. Shorter paths have a small \texttt{ne} contribution, but balanced by the larger weight and vice versa for longer paths. This balance of the resulting score reduces the variance.

The \texttt{natl\_tl} estimation procedure becomes optimal when scattering is trivial ($\VAPr=1$), since its variance is zero then, since the \texttt{natl\_tl} score only depends on the travelled length, which will always be exactly $\VAdomlength$ for an \texttt{natl} simulation with $\VAPr=1$. 
Track-length aspects also become beneficial when the particle trajectories become more complex (high survival probability $\VAcss/\VAcst$) because it takes the more global `path length' into account and when the total collisionality is low, the situation for which it was originally designed. The precise cut-offs are non-trivial in each of the cases.

In the rest of the 1D0D parameter space --- being the part with $\VAcst\VAdomlength=\frac{\VAratet}{|v|}\VAdomlength>10$ --- the picture remains intact. There remain significant portions with the \texttt{natl\_ne} estimator having the upper hand and the \texttt{natl\_tl} estimation procedure remains optimal for a small fraction near $\VAcss/\VAcst=1$ and $\VAPr=1$.

\begin{figure}[H]\centering
\begin{tabular}{b{.25\textwidth}b{.25\textwidth}b{.25\textwidth}b{.2\textwidth}}
\resizebox{.24\textwidth}{!}{\maakmooieticks{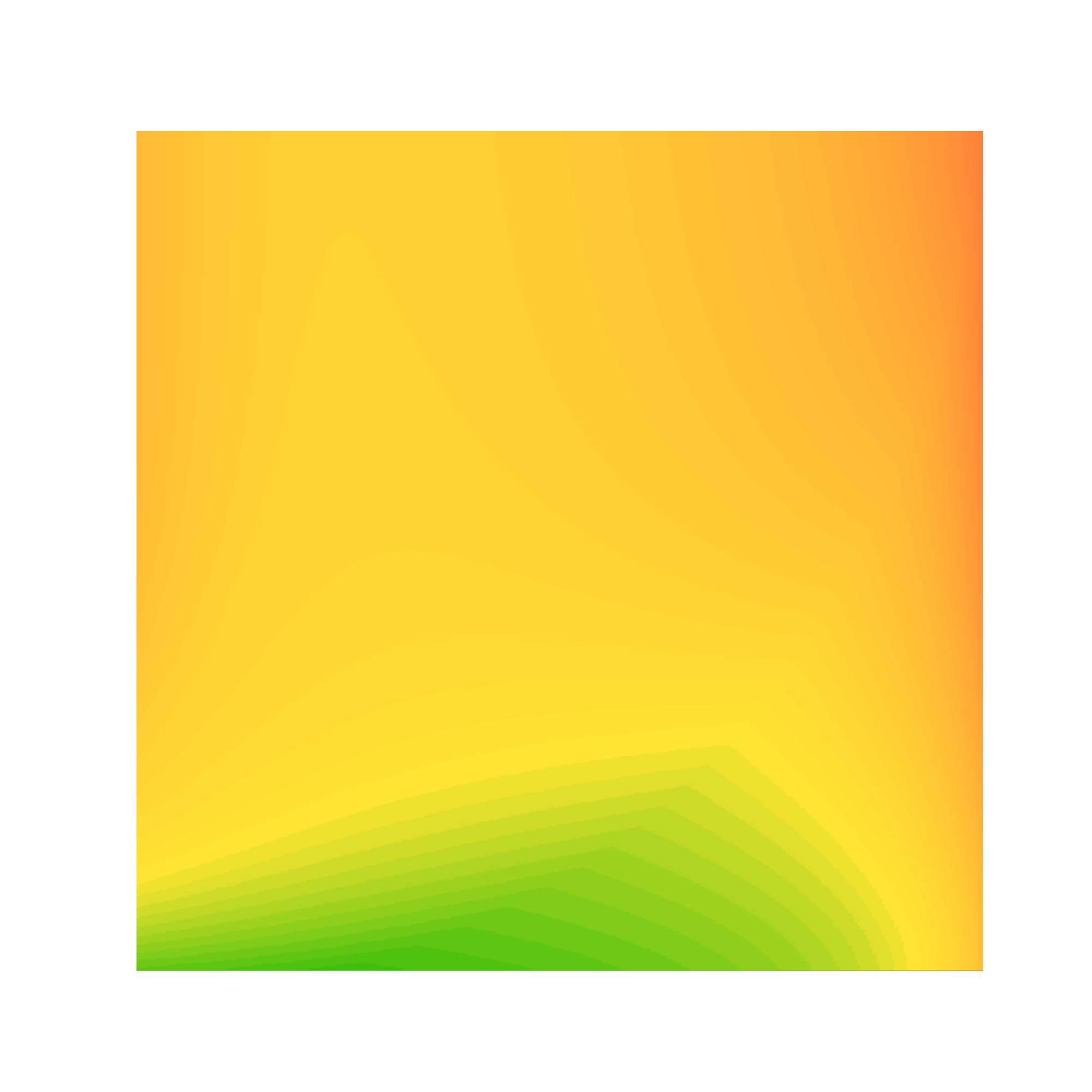}{1}{1}{$\dfrac{\VAcss}{\VAcst}=0.75$}} &
\resizebox{.24\textwidth}{!}{\maakmooieticks{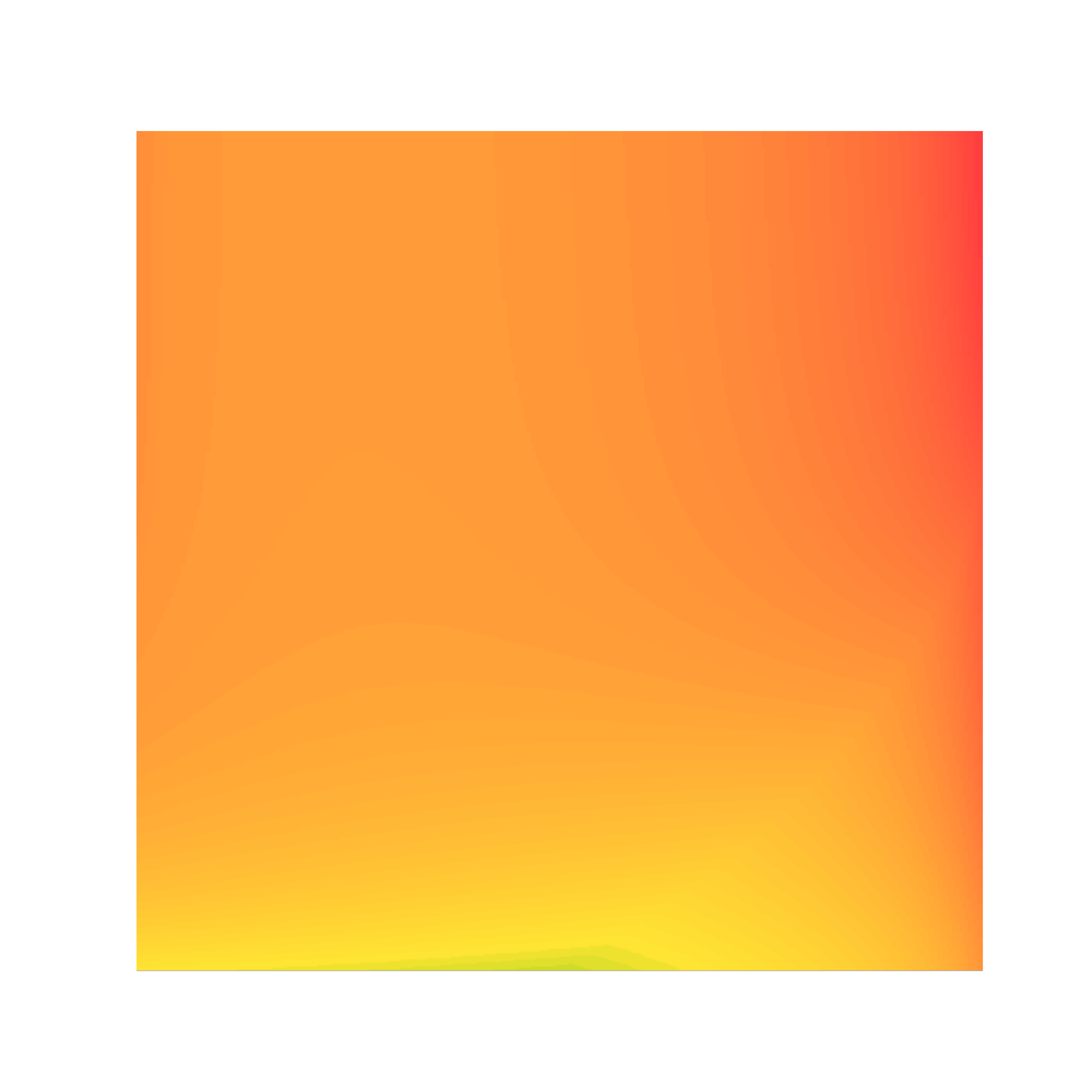}{0}{1}{$\dfrac{\VAcss}{\VAcst}=0.5$}} &
\resizebox{.24\textwidth}{!}{\maakmooieticks{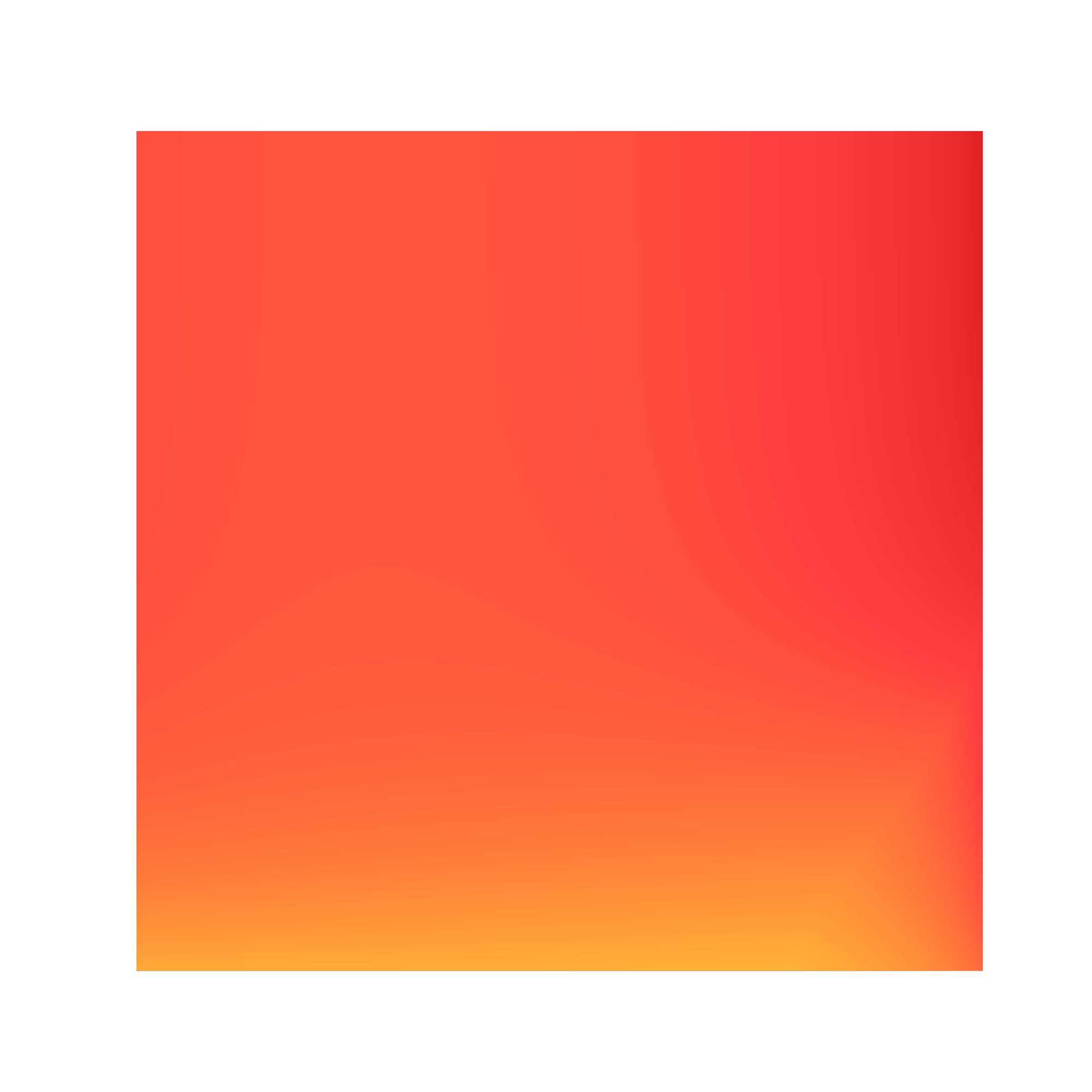}{0}{1}{$\dfrac{\VAcss}{\VAcst}=0.25$}} &
\resizebox{!}{.24\textwidth}{\maakmooietickscb{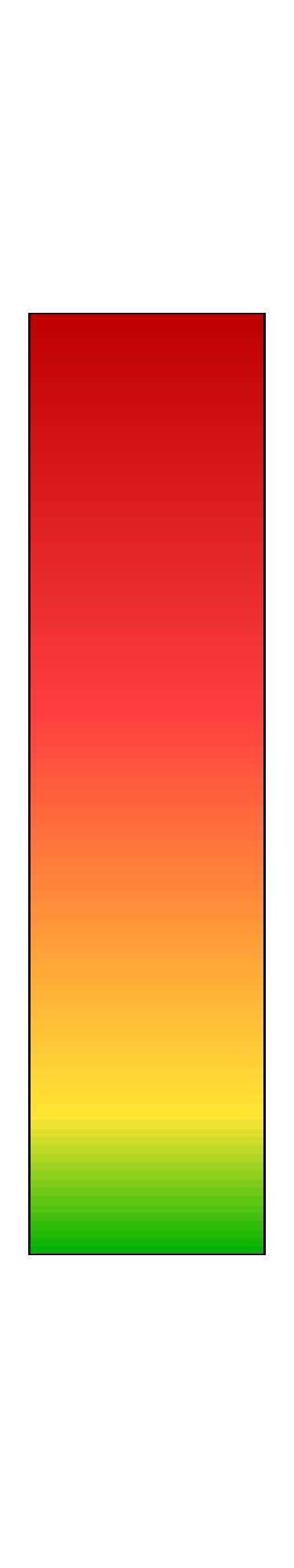}{}}
\end{tabular}
\caption{The factor by which the standard deviation increases when, in a 1D0D setting, the standard mass source estimation procedure choice, \texttt{a\_tl}, is used instead of the best estimation procedure as depicted in Figure~\ref{fig:best_1D0D_mass_stdv_bestest}.}
\label{fig:best_1D0D_mass_stdv_gain}
\end{figure}

In the EIRENE code, the default mass source estimation procedure is the \texttt{a\_tl} estimation procedure, which is not one of the three competitive estimation procedures found from the comparison based on statistical error in the 1D0D case. In Figure~\ref{fig:best_1D0D_mass_stdv_gain}, we present the factor by which the standard deviation of the \texttt{a\_tl} estimation procedure is higher than the lowest possible over all estimation procedures. As can be seen in the figure, except for a very low collisionality, or when the survival rate is high and the collisionality not high, the potential gain by using the optimal estimation procedure is very large. Before we draw conclusions, we must first note that these results are preliminary, since they apply to the simplified 1D0D case and currently only look at the statistical error. In a fusion setting,  survival probability and collisionality are typically very high, meaning that the top part of the left-most figure from Figure~\ref{fig:best_1D0D_mass_stdv_gain} should be considered. There, these preliminary results indicate that the current \texttt{a\_tl} choice is not optimal when considering the statistical error in a 1D0D setting.

\subsubsection{Best estimation procedure for the mass source in the 1D0D model based on computational cost\label{subsubsec:iif_bestest_1D0D_bestest_cost}}

A more relevant measure of performance than the statistical error, is the computational cost for a given statistical error. The number of particles required to obtain a certain statistical error is proportional to the variance of a single particle contribution, $\sigma^2$ and we approximate the simulation cost as being proportional to the number of scattering collisions. Combined, this means that the total computational cost for a given statistical error is proportional to
\begin{equation}
\sigma^2\mathbb{E}[\text{collisions per path}]\,.
\end{equation}
The number of collisions per path is independent of the chosen estimator (\texttt{a\_abs}, \texttt{a\_sc}, \texttt{c}, \texttt{tl}, \texttt{ne}), but solely depends on the simulation type (\texttt{a}, \texttt{nac}, \texttt{natl}).

\begin{figure}[H]\centering
\begin{tabular}{m{.25\textwidth}m{.25\textwidth}m{.25\textwidth}m{.2\textwidth}}
\resizebox{.24\textwidth}{!}{\maakmooieticks{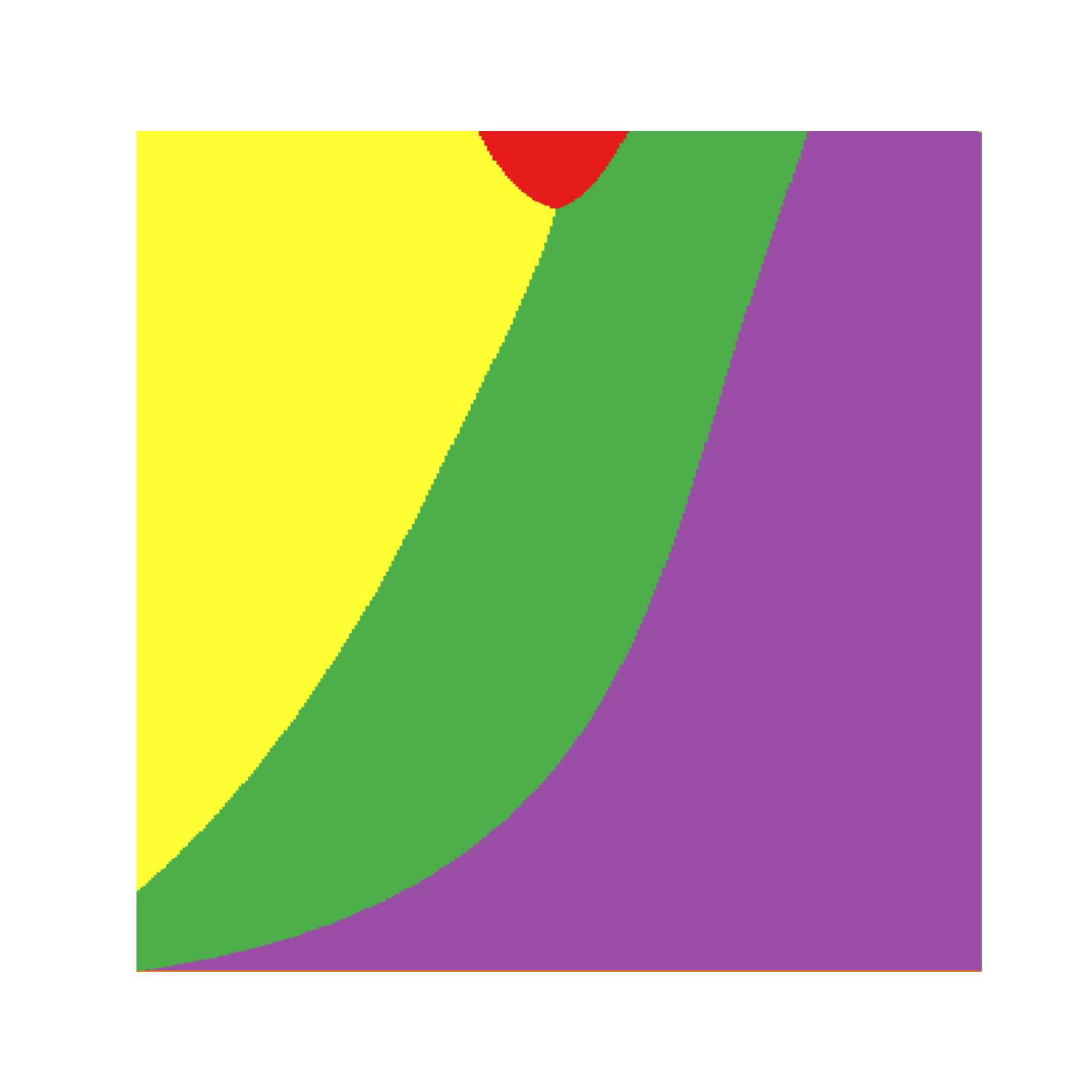}{1}{0}{$\dfrac{\VAcss}{\VAcst}=0.75$}} &
\resizebox{.24\textwidth}{!}{\maakmooieticks{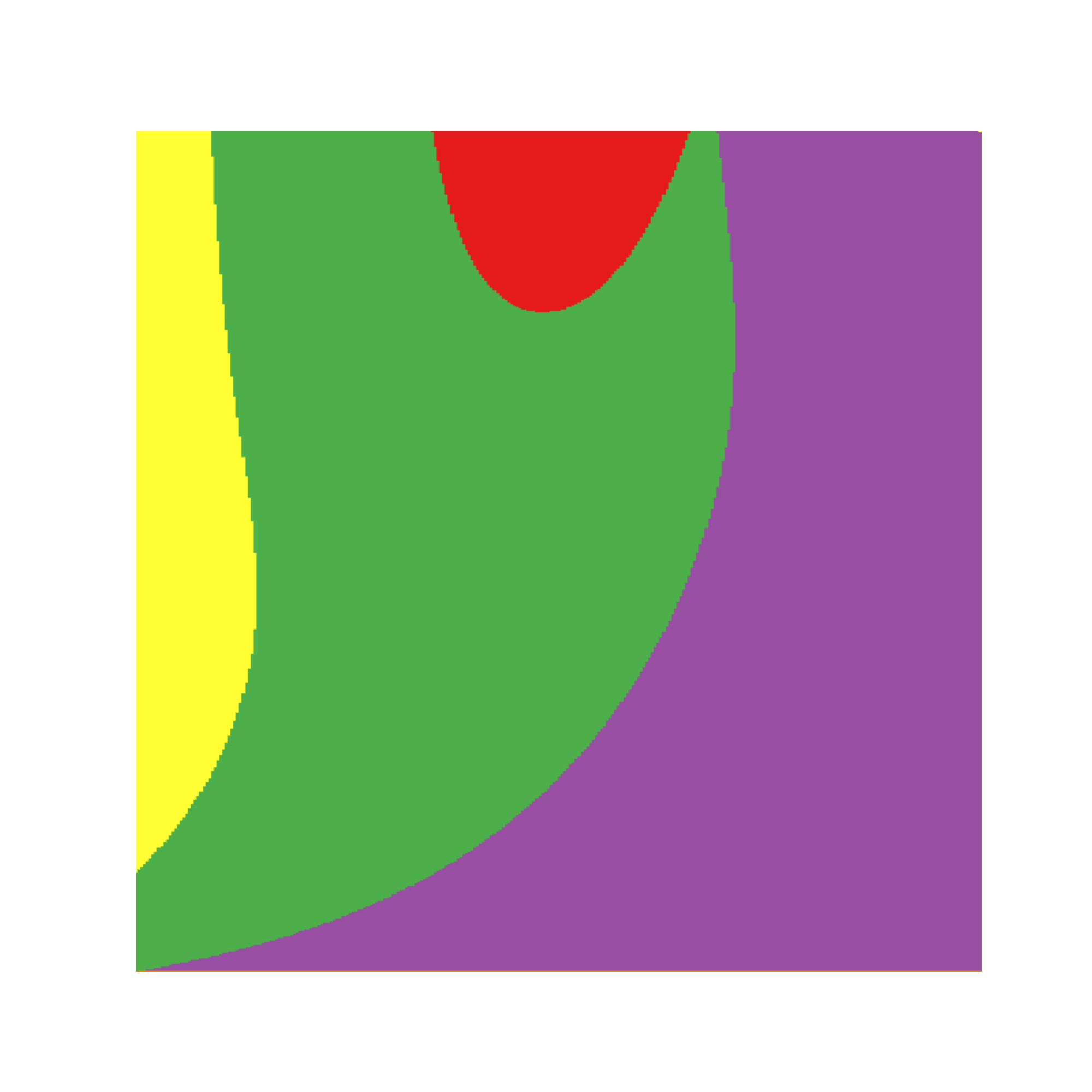}{0}{0}{$\dfrac{\VAcss}{\VAcst}=0.5$}} &
\resizebox{.24\textwidth}{!}{\maakmooieticks{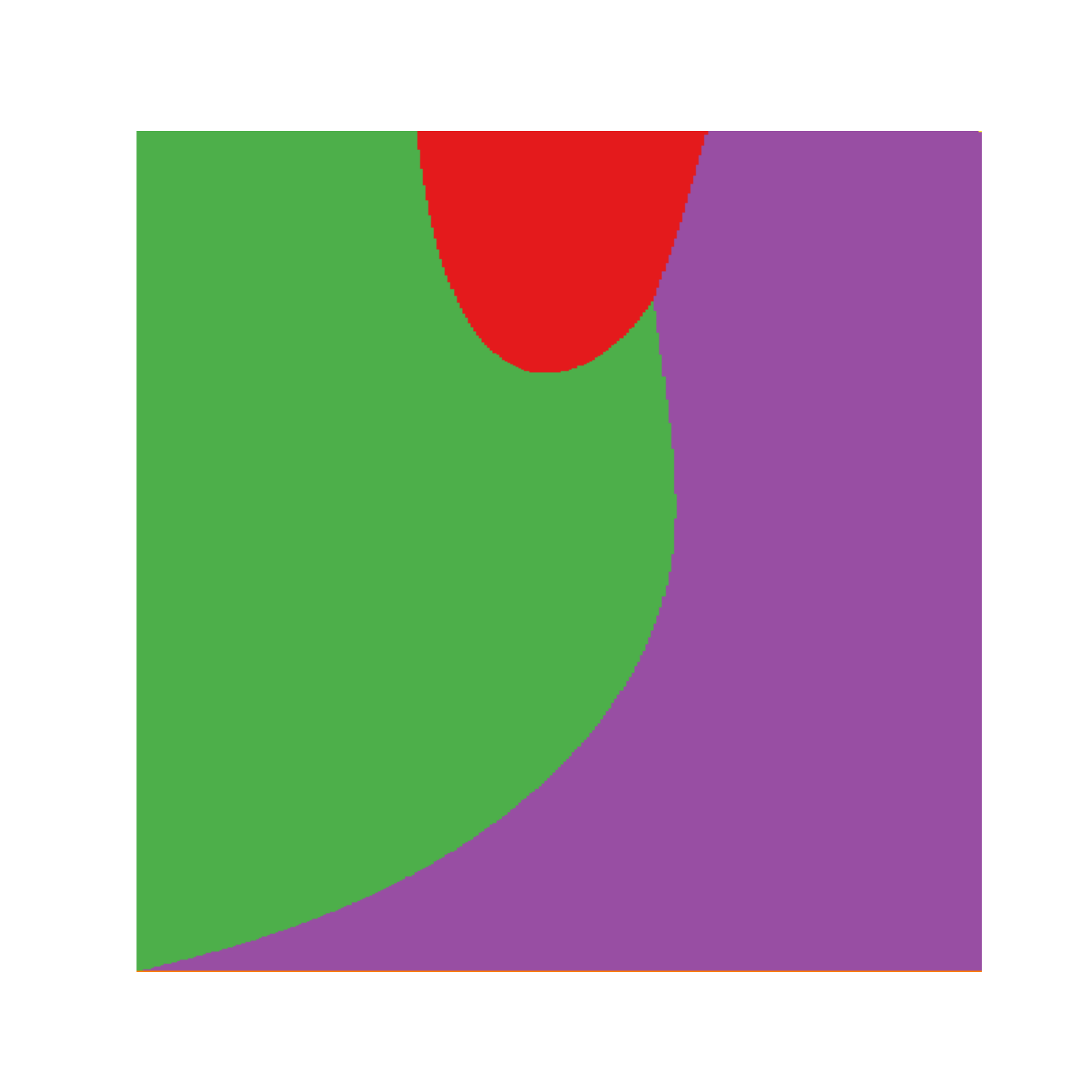}{0}{0}{$\dfrac{\VAcss}{\VAcst}=0.25$}} &
\localcolorlegendoneB
\end{tabular}\\
\begin{tabular}{b{.25\textwidth}b{.25\textwidth}b{.25\textwidth}b{.2\textwidth}}
\vspace{-.5cm}\resizebox{.24\textwidth}{!}{\maakmooieticks{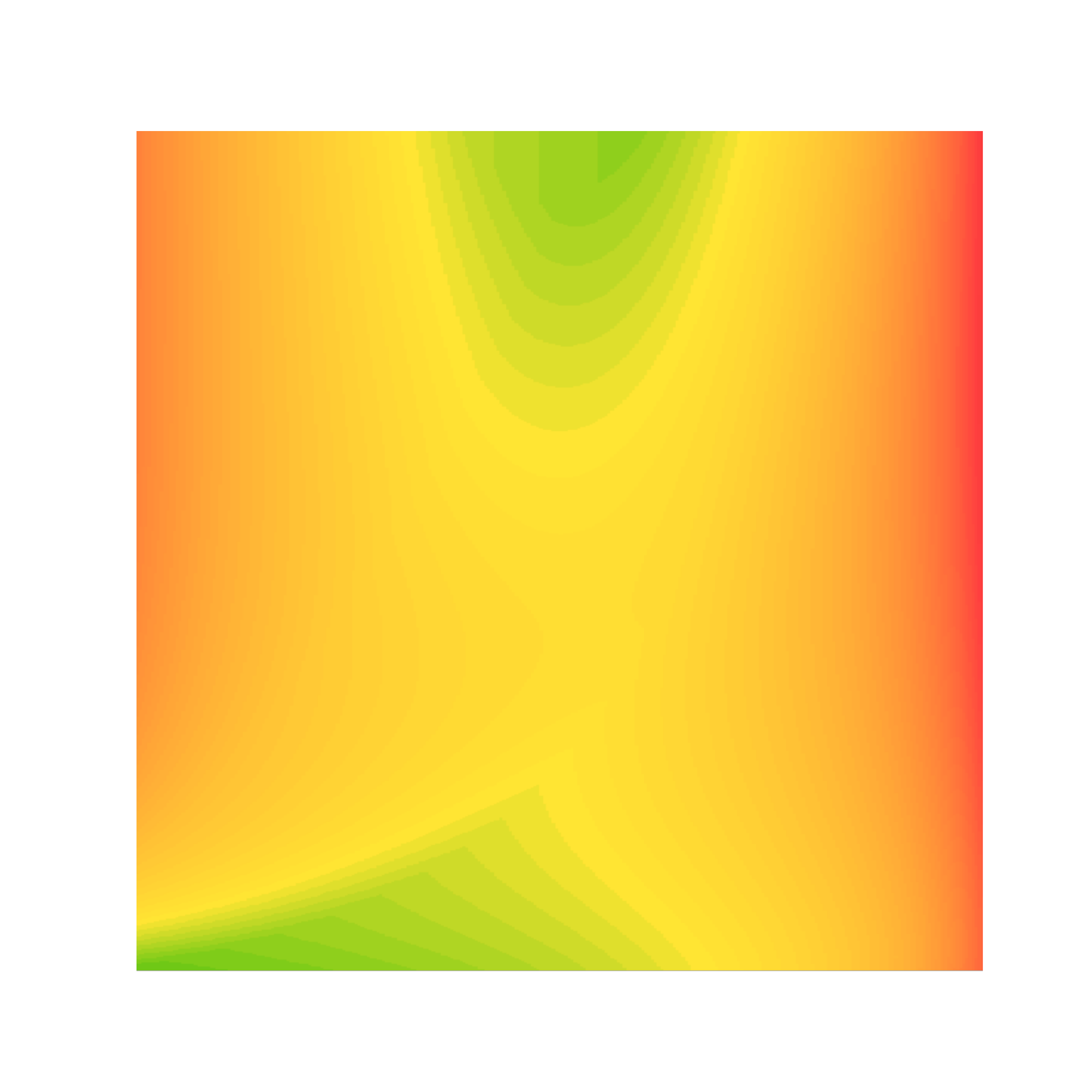}{1}{1}{}} &
\vspace{-.5cm}\resizebox{.24\textwidth}{!}{\maakmooieticks{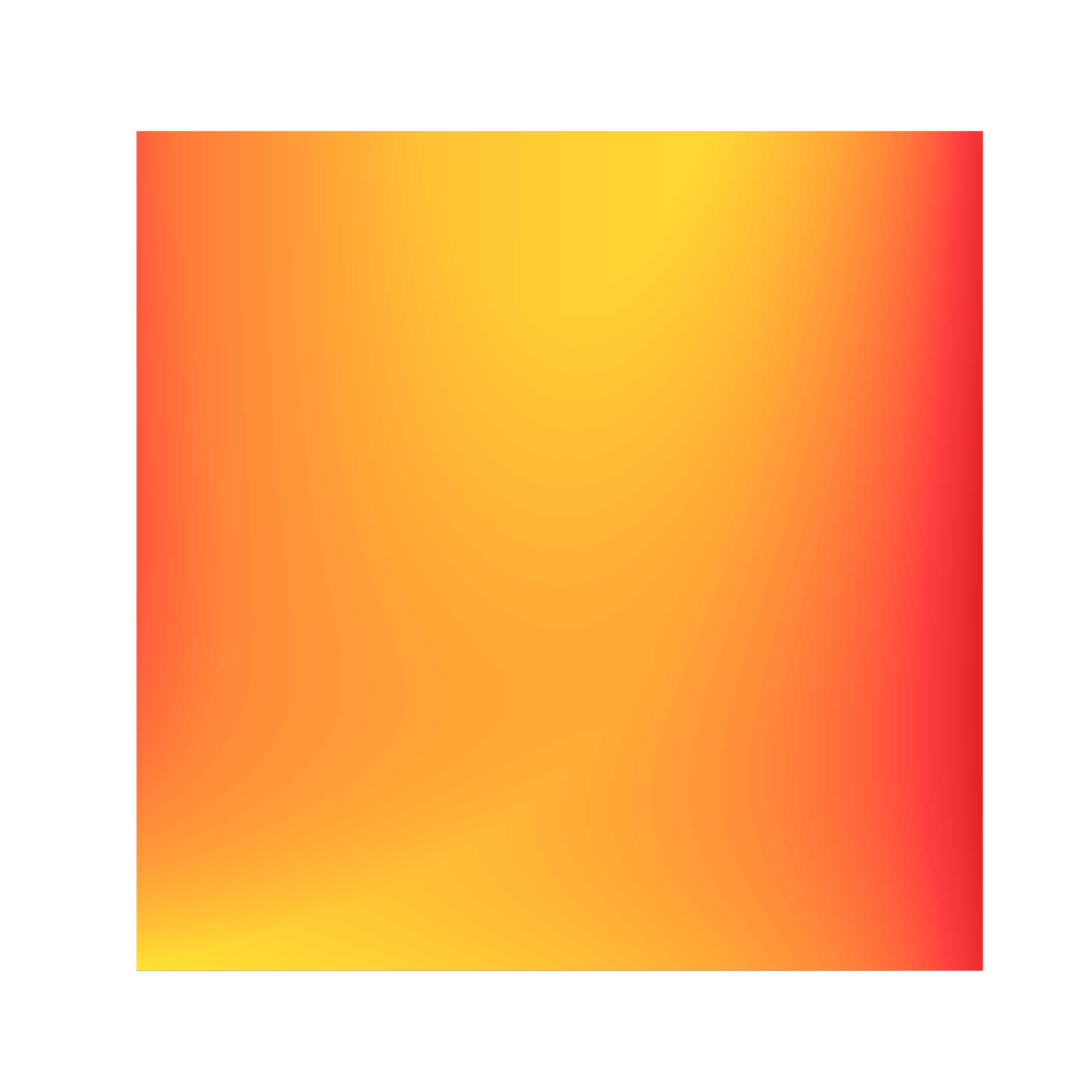}{0}{1}{}} &
\vspace{-.5cm}\resizebox{.24\textwidth}{!}{\maakmooieticks{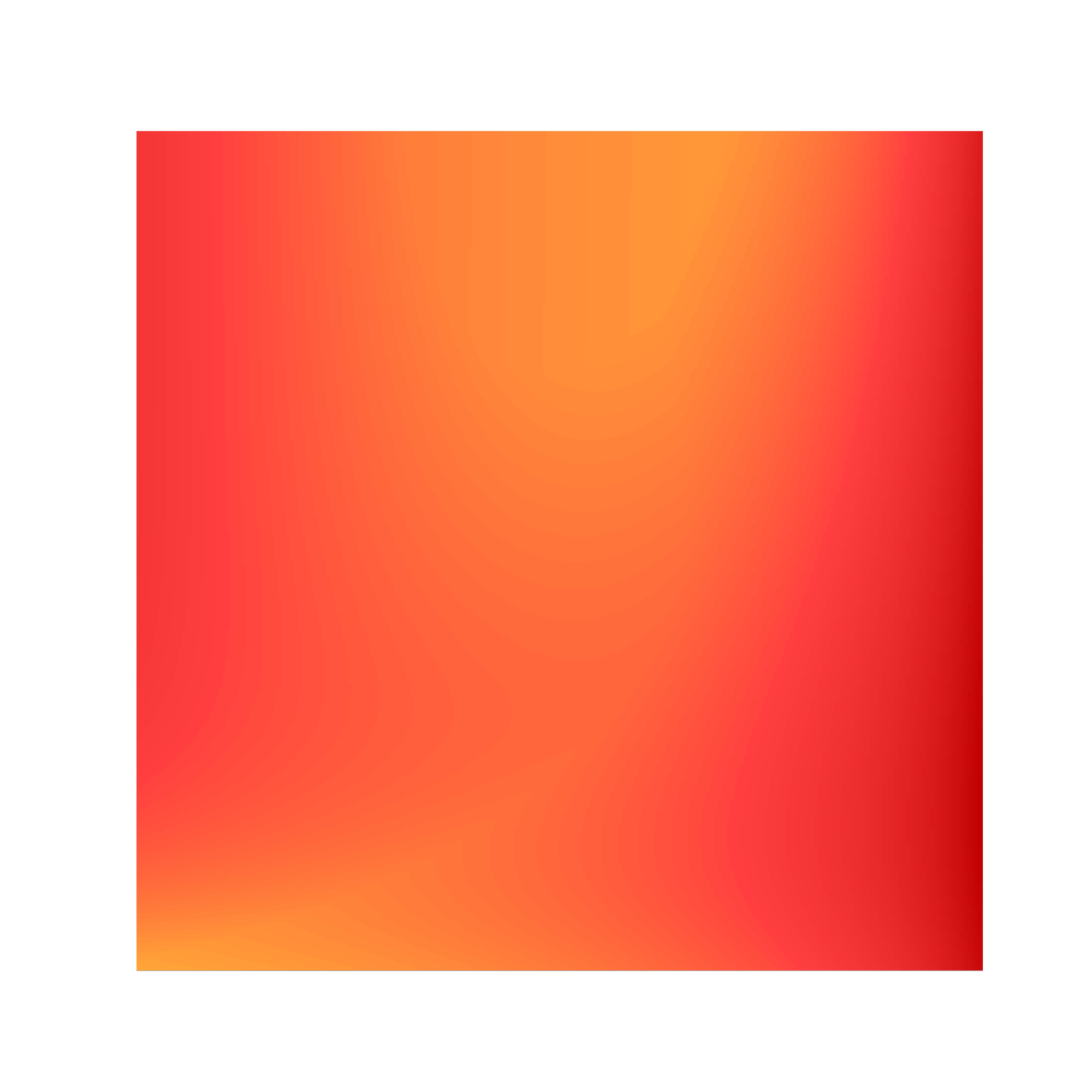}{0}{1}{}} &
\resizebox{!}{.24\textwidth}{\maakmooietickscb{figuren/improvement/improvement_colorbar.pdf}{}}
\end{tabular}
\caption{The best 1D0D mass source term estimation procedure and the potential gain when considering computational cost for a given statistical error. The first row presents a partition of the parameter domain based on the estimation procedure with the lowest statistical error. The second row presents the factor by which the standard deviation increases when the standard mass choice, \texttt{a\_tl}, is used instead of the best estimation procedure.}
\label{fig:best_1D0D_mass_cost}
\end{figure}

Figure~\ref{fig:best_1D0D_mass_cost} shows the same results as Figures~\ref{fig:best_1D0D_mass_stdv_bestest} and~\ref{fig:best_1D0D_mass_stdv_gain} but now with the computational cost as measure of performance. The connection with the relative standard deviation is that it is squared and multiplied by the simulation-dependent expected number of collisions. Hence, if there would have been a border between two estimation procedures with the same simulation type in figure~\ref{fig:best_1D0D_mass_stdv_bestest}, that border would not be different now.

The number of collisions is always lowest for \texttt{a} and highest for \texttt{nac}. These effects are visible when comparing of the first row of Figure~\ref{fig:best_1D0D_mass_cost} to Figure~\ref{fig:best_1D0D_mass_stdv_bestest}: \texttt{natl\_tl} and \texttt{natl\_ne} take over parts of the domain that went to \texttt{nac\_ne} in Figure~\ref{fig:best_1D0D_mass_stdv_bestest}. On top of that, a procedure with an analog simulation becomes competitive: \texttt{a\_ne}. For small values of the survival fraction and not too low values of $\VAcst\VAdomlength$, non-analogous particles can undergo many collisions without having a significant contribution to the estimation. Those long-living non-analogous particles constitute a computational cost, but have nearly no impact on the computed value. This is the reason analog simulations can also be the best option.

When increasing $\VAcst\VAdomlength$ now beyond 10, the effect on the expected number of scattering collisions also has to be taken into account. A simple calculation gives that this number is finite when $\VAcst\VAdomlength\rightarrow\infty$ when $\VAPr<0.5$. For values of $\VAPr$ larger than $0.5$, an analog estimation procedure, \texttt{a\_ne} for the most part, becomes optimal when $\VAcst\VAdomlength$ increases, except of course for $\VAPr=1$, since \texttt{natl\_tl} then provides zero variance and hence also zero in our computational cost measure. When $\VAPr\leq0.5$ there remains a significant fraction of the parameter domain where \texttt{natl\_ne} and \texttt{nac\_ne} are optimal, no matter how large $\VAcst\VAdomlength$ is.

When we now compare the best estimation procedure as determined by the first row of Figure~\ref{fig:best_1D0D_mass_cost} with the default choice \texttt{a\_tl}, which is shown in the second row of Figure~\ref{fig:best_1D0D_mass_cost}, we find it is now also a reasonable choice in a high-collisional isotropic situation when the survival probability is high. This is the most important regime in fusion research, meaning \texttt{a\_tl} is a proper choice there. There will be many different parameter values throughout the domain however, and as is visible in the second row of Figure~\ref{fig:best_1D0D_mass_cost}, for many parameter values \texttt{a\_tl} is not a reasonable choice.

It might be that selecting a single one estimation procedures for all parameter values, or in the context of a fusion simulation, for all grid cells, gives decent results nonetheless. This is however not the case. No matter what estimation procedure you pick, there is always a parameter set for which another estimation procedure performs infinitely better. Indeed, when the statistical error is regarded, the \texttt{natl\_tl} procedure is the only one to attain zero relative error when $\VAPr=1$, so it would be infinitely better than all others. On the other hand, the \texttt{nac\_ne} procedure outperforms it in a large part of the domain, and it can outperform it unboundedly. Namely, when the survival probability is very low, a \texttt{nac\_ne} contribution will be nearly deterministic, resulting in zero standard deviation and zero cost. If furthermore $\VAPr$ is low, the variance of \texttt{natl\_tl} can become relatively high. Furthermore, as discussed above, when increasing $\VAcst\VAdomlength$ beyond 10, \texttt{a\_ne} becomes infinitely better in terms of cost than any estimation procedure with a non-analog simulation type, when $\VAPr\ge 0.5$ and $\VAcst\VAdomlength$ becomes large.

Since, depending on the parameter set, there can be an unbounded difference in performance between estimation procedures, it can be beneficial to use different estimation procedures in different parts of the domain. Doing so is allowed, since each of estimators is unbiased and each of the simulation types corresponds to the same model.

In Section~\ref{sec:ii_res_1D1D}, we will see how well these mass source term estimation results hold up in a 1D1D setting, and we will consider momentum source term estimation.

\section{Numerical extensions to a realistic one-dimensional grid cell and momentum estimation\label{sec:ii_res_1D1D}}

In this section, we extend our analytical results in two ways: by considering continuous velocities and by considering momentum source estimation. In a setting with continuous velocities, the invariant imbedding method used in Section~\ref{sec:iif_bestest_1D0D} becomes intractable. We thus resort to numerical experiments in which an MC estimation is repeated multiple times for each estimation procedure and each parameter value and the statistical properties were computed. Because these results are obtained via this numerical method, the different parameter values are on a much coarser mesh in the parameter space than in Section~\ref{sec:iif_bestest_1D0D}. The extension to momentum is also only numerically available, both for 1D0D and 1D1D, and shows the importance of using different estimation procedures for the different quantities in most settings.

In Section~\ref{subsec:iif_1D1D_1D1D}, we explain the more realistic 1D1D simulation and present a mapping from 1D1D to 1D0D to allow a comparison of results in both settings. In Section~\ref{subsec:iif_1D1D_mass} we then present the numerical 1D1D mass source estimation results and in Section~\ref{subsec:iif_1D1D_mom}, we present the numerical extension to momentum source estimation.

\subsection{A setting with continuous velocities: 1D1D simulation\label{subsec:iif_1D1D_1D1D}}

The 1D1D setting which we experiment on, is identical to the 1D0D setting described in Section~\ref{subsec:iif_bestest_1D0D_1D0D}, except for the post-collisional velocity distribution. In a general setting, the post-collisional velocity distribution can take any form, but in most applications, a Maxwellian distribution is a proper first approximation. We will consequently focus on a Maxwellian with mean $\mu_v$ and variance $\sigma_v^2$:
\begin{equation}
\VApostcolveldistr^\text{1D1D}(v)=\frac{1}{\sqrt{2\pi}\sigma_v}e^{\frac{(v-\mu_v)^2}{2\sigma_v^2}}\,.
\end{equation}
We further retain the plasma background homogeneity (constant $\mu_v$, $\sigma_v$, $\VArates$, and $\VAratea$), fully absorbent walls, and $(x_0,v_0)=(0,1)$. Now the velocity can have different sizes, meaning the cross-section ($\VAcs=\VArate/|v|$) is no longer constant. We will however still present our results in terms of the 1D0D parameters $\VAcst\VAdomlength$, $\VAcss/\VAcst$, and $\VAPr$, by using an appropriate mapping.

\begin{figure}[H]
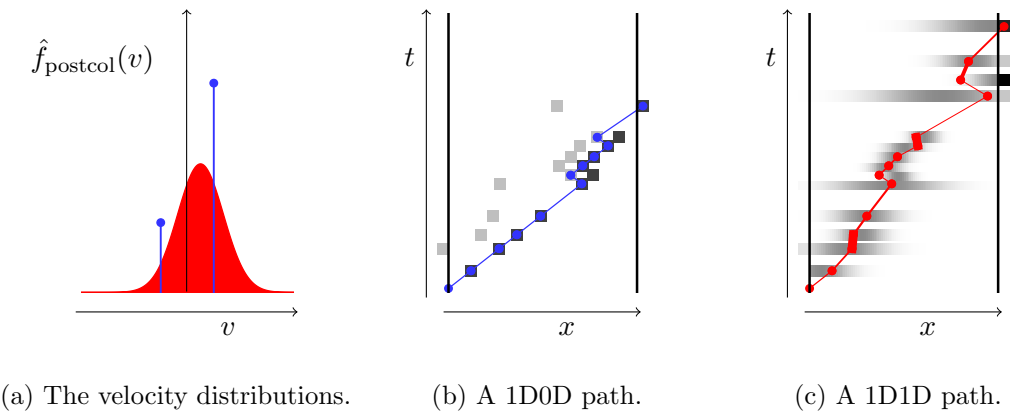

\centering
\begin{subfigure}{.28\columnwidth}
	\centering
	\oneDandzeroDpostcolveldistrplot
	\caption{The velocity distributions.}
	\label{fig:ii_1D0D1D1D_sketch_veldistr}
\end{subfigure}
\begin{subfigure}{.28\columnwidth}
	\centering
	\zeroDpathplot
	\caption{A 1D0D path.}
	\label{fig:ii_1D0D1D1D_sketch_1D0Dpath}
\end{subfigure}
\begin{subfigure}{.28\columnwidth}
	\centering
	\oneDpathplot
	\caption{A 1D1D path.}
	\label{fig:ii_1D0D1D1D_sketch_1D1Dpath}
\end{subfigure}
\caption{A sketch illustrating the difference between forward-backward scattering (blue) and more realistic scattering (red) by plotting two correlated paths.}
\label{fig:ii_1D0D1D1D_sketch}
\end{figure}

The mapping we will use from the 1D1D to the 1D0D parameter space is such that the first and second moment of the post-collisional velocity distribution match. This is achieved by choosing the parameter $\VAPr$ appropriately and by rescaling the velocity, which is implemented by rescaling the cross-sections to obtain an equivalent result. The survival probability is not modified. Such a mapping is illustrated in Figure~\ref{fig:ii_1D0D1D1D_sketch_veldistr}, where a 1D0D and a 1D1D post-collisional velocity distribution are shown with identical first and second moments.

By changing the post-collisional velocity distribution, the neutral particle behaviour changes. Figures~\ref{fig:ii_1D0D1D1D_sketch_1D0Dpath} and~\ref{fig:ii_1D0D1D1D_sketch_1D1Dpath} illustrate some of the effects by showing a neutral path with forward-backward scattering in Figure~\ref{fig:ii_1D0D1D1D_sketch_1D0Dpath} and a neutral path with scattering according to a normal distribution with equal first two moments in Figure~\ref{fig:ii_1D0D1D1D_sketch_1D1Dpath}. Besides matching of the first two moments of the distributions, the paths have also been sampled in a correlated manner, resulting in equal collision times and similar but not equal velocities. In these figures we have also illustrated the post-collisional velocity distribution by the greyscale squares, respectively lines. These show the probability distribution of the next collision, given the collision time and the previous collision position are known. This visualisation lays bare the most prevalent difference in the neutral behaviour between the two cases, namely the additional stochasticity on the travelled length. One derivative effect of this is the interaction with the boundaries: in a 1D1D simulation, crossing the boundary as a next event is always possible (indicated by the (sometimes extremely faintly) small squares outside of the domain), which is not the case for a 1D0D simulation, where, with given collision times, the probability of leaving the domain is only non-zero close to the boundaries. 

Another aspect that can impact the variance associated with estimation procedures comes from the dependence of the velocity in many scoring procedures. Slow parts of the path, may not contribute much to the distance travelled, but have a much larger impact on the score and weight than their slower counterparts. The velocity of the particle is clear in Figure~\ref{fig:ii_1D0D1D1D_sketch_1D1Dpath} by the slope of the lines, but we have emphasized the importance of the slow parts by having the thickness of the lines be proportional to $1/v_\VAeventno$. This effect was crucial in the prevalence of \texttt{natl\_tl} in Section~\ref{subsubsec:iif_bestest_1D0D_bestest_var}.

These new effects can impact the different estimation procedures profoundly, as we will investigate in the next section.

\subsection{Mass source estimation in 1D1D\label{subsec:iif_1D1D_mass}}

\begin{figure}[H]\centering
\begin{tabular}{m{.25\textwidth}m{.25\textwidth}m{.25\textwidth}m{.2\textwidth}}
\resizebox{.24\textwidth}{!}{\maakmooieticks{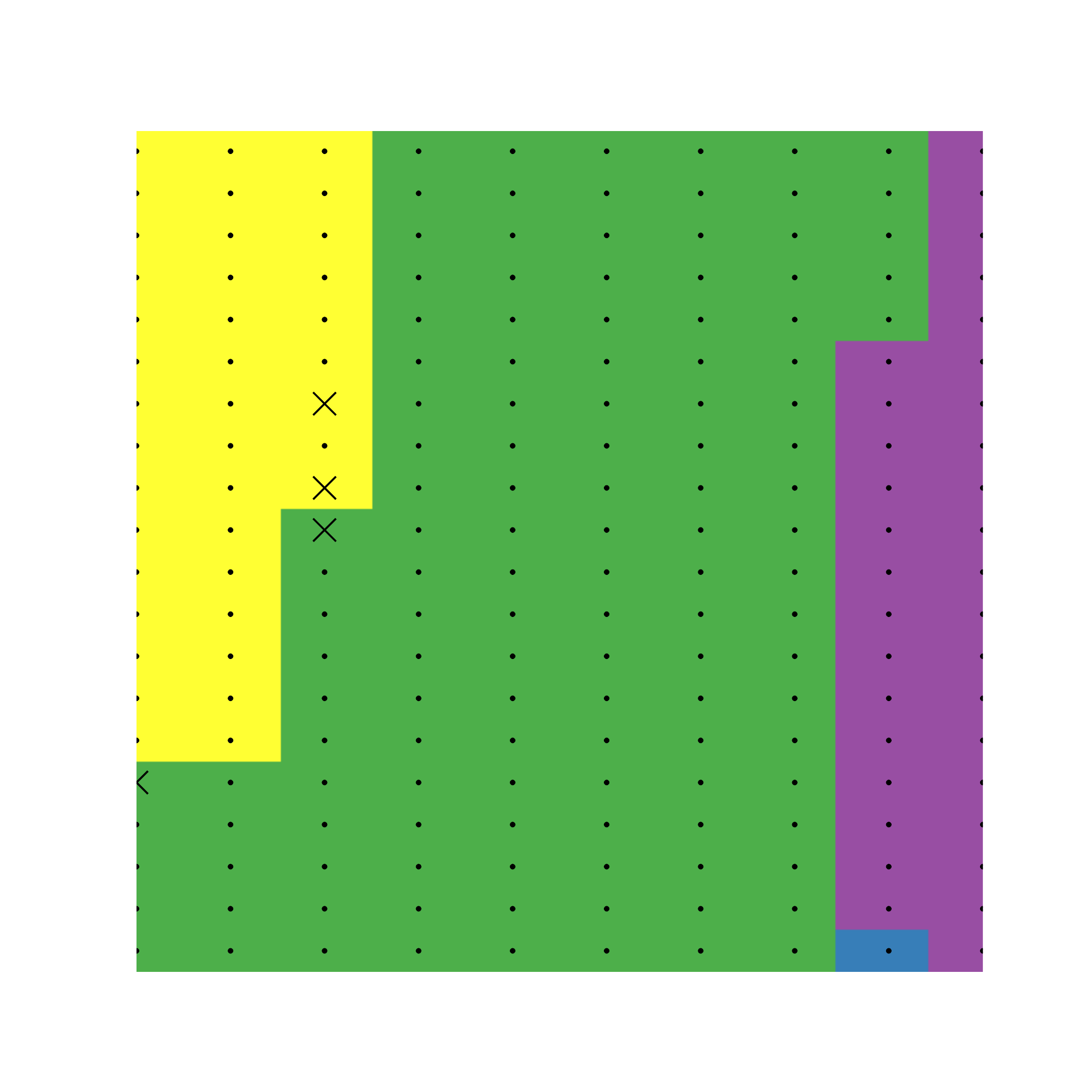}{1}{0}{$\dfrac{\VAcss}{\VAcst}=0.75$}} &
\resizebox{.24\textwidth}{!}{\maakmooieticks{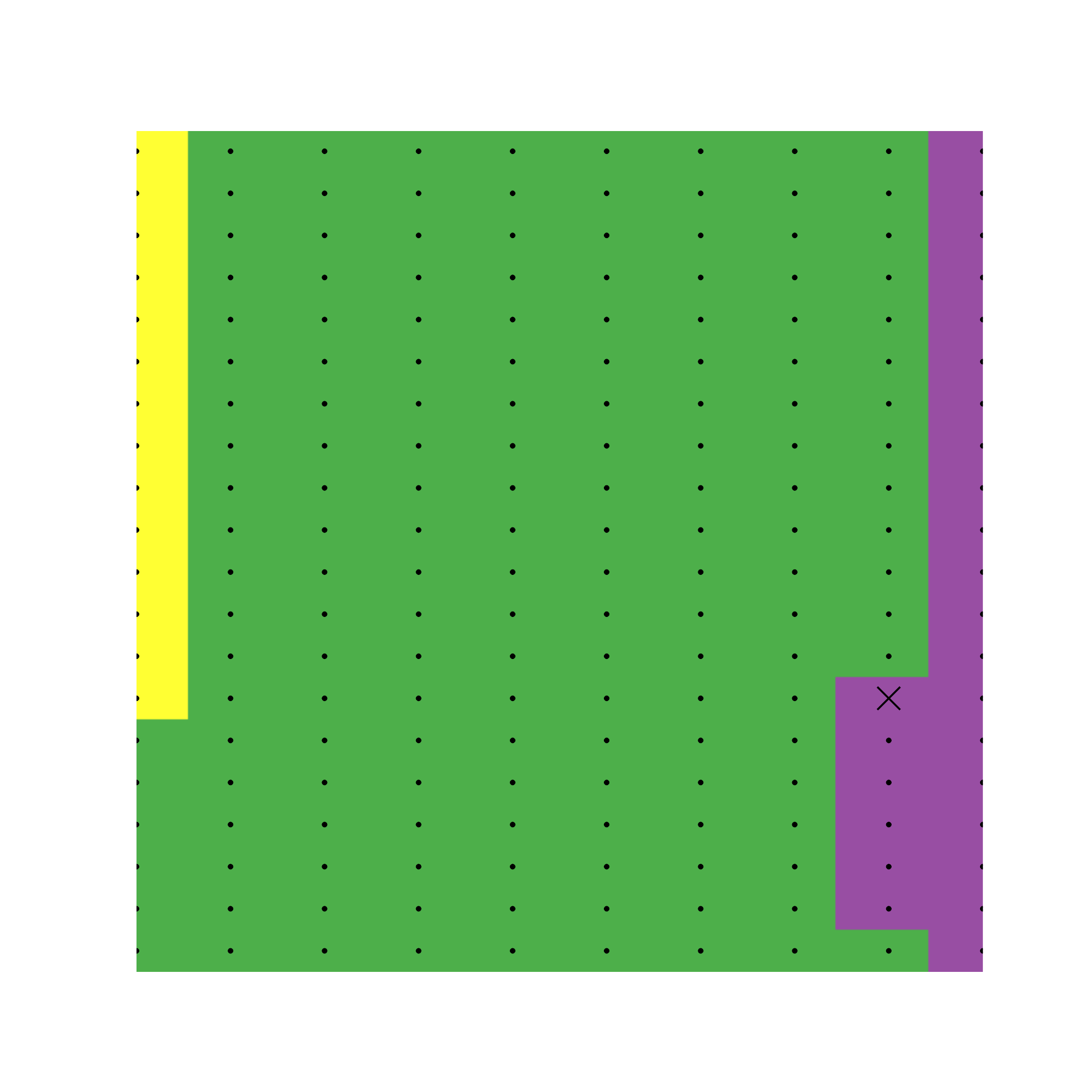}{0}{0}{$\dfrac{\VAcss}{\VAcst}=0.5$}} &
\resizebox{.24\textwidth}{!}{\maakmooieticks{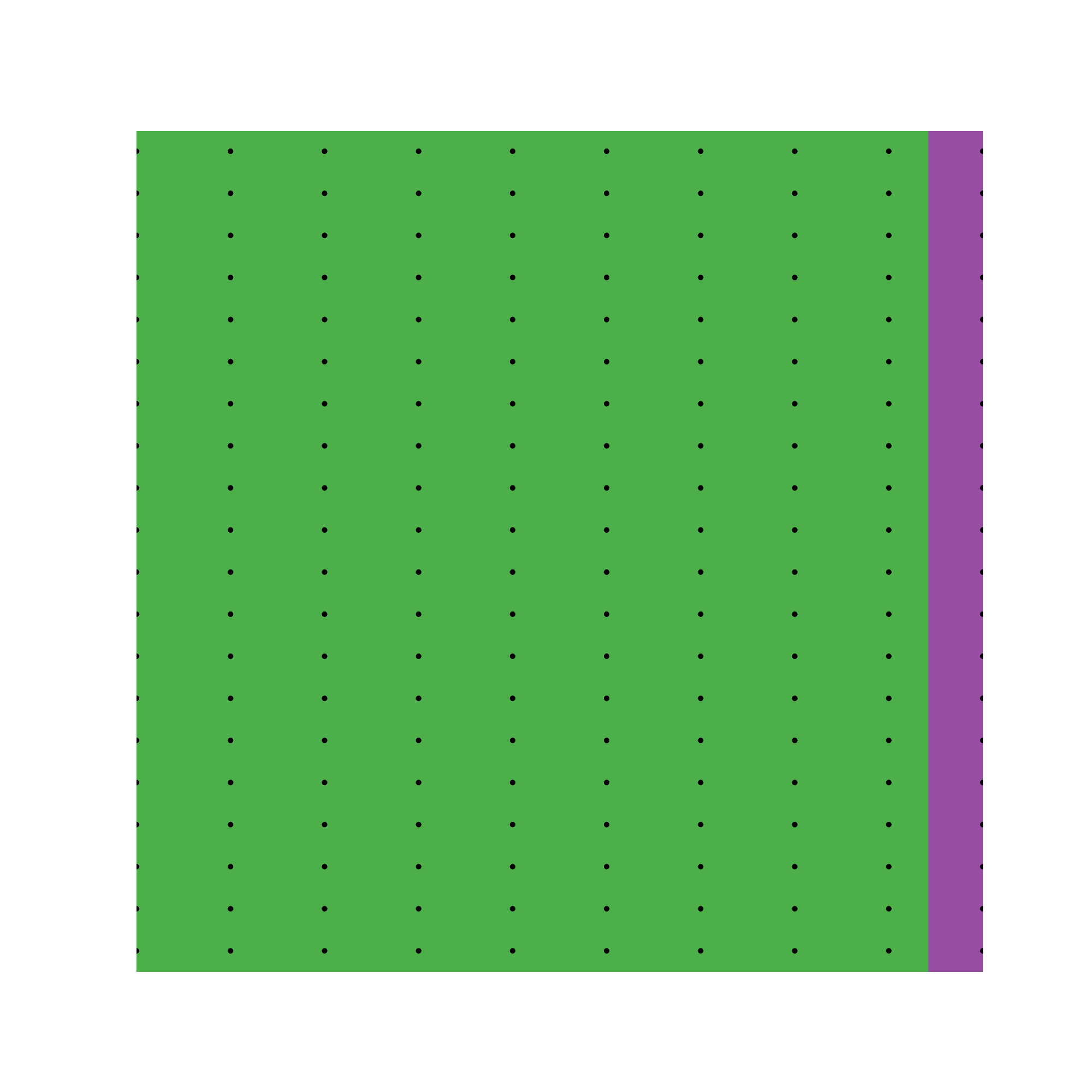}{0}{0}{$\dfrac{\VAcss}{\VAcst}=0.25$}} &
\localcolorlegendoneA
\end{tabular}\\
\begin{tabular}{b{.25\textwidth}b{.25\textwidth}b{.25\textwidth}b{.2\textwidth}}
\vspace{-.5cm}\resizebox{.24\textwidth}{!}{\maakmooieticks{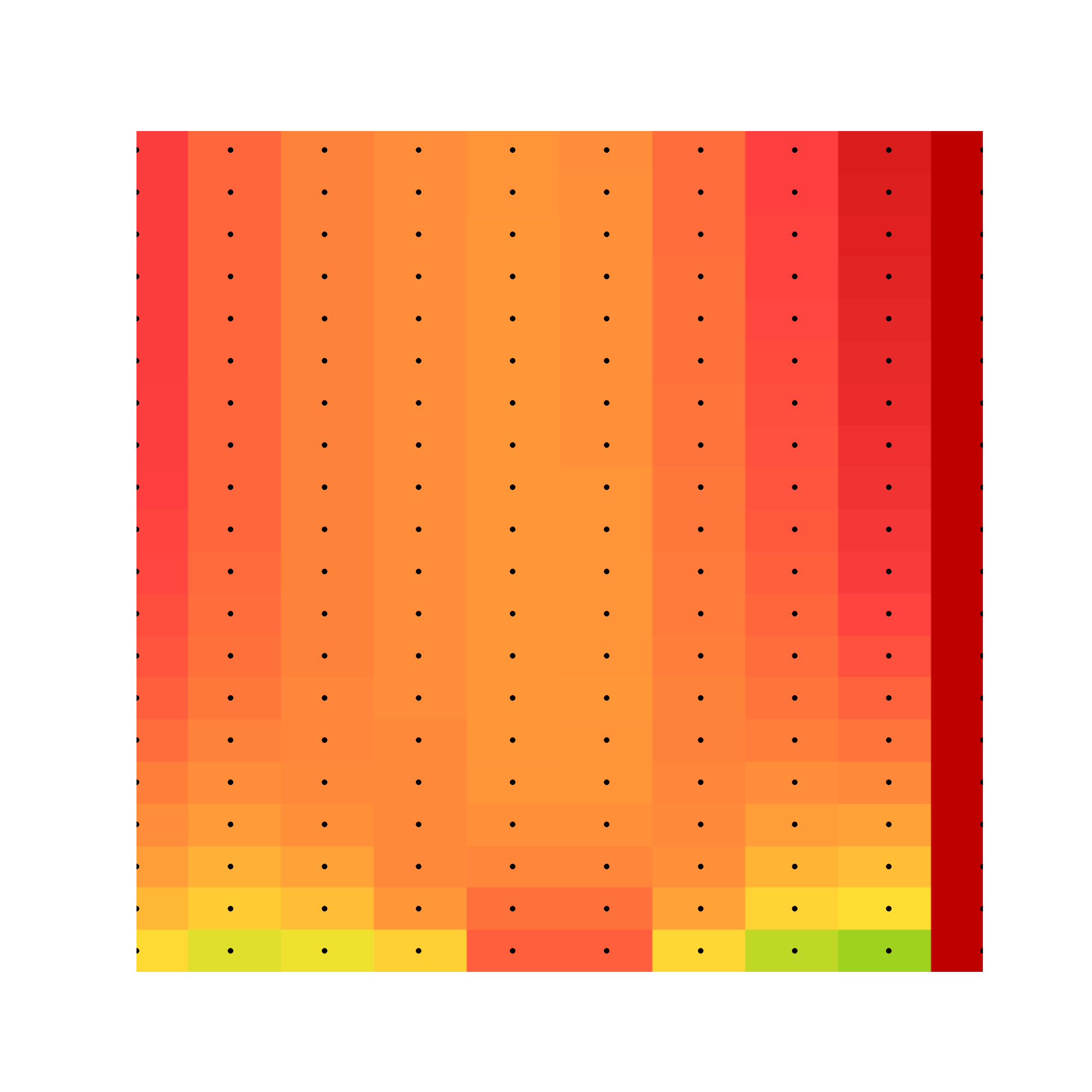}{1}{1}{}} &
\vspace{-.5cm}\resizebox{.24\textwidth}{!}{\maakmooieticks{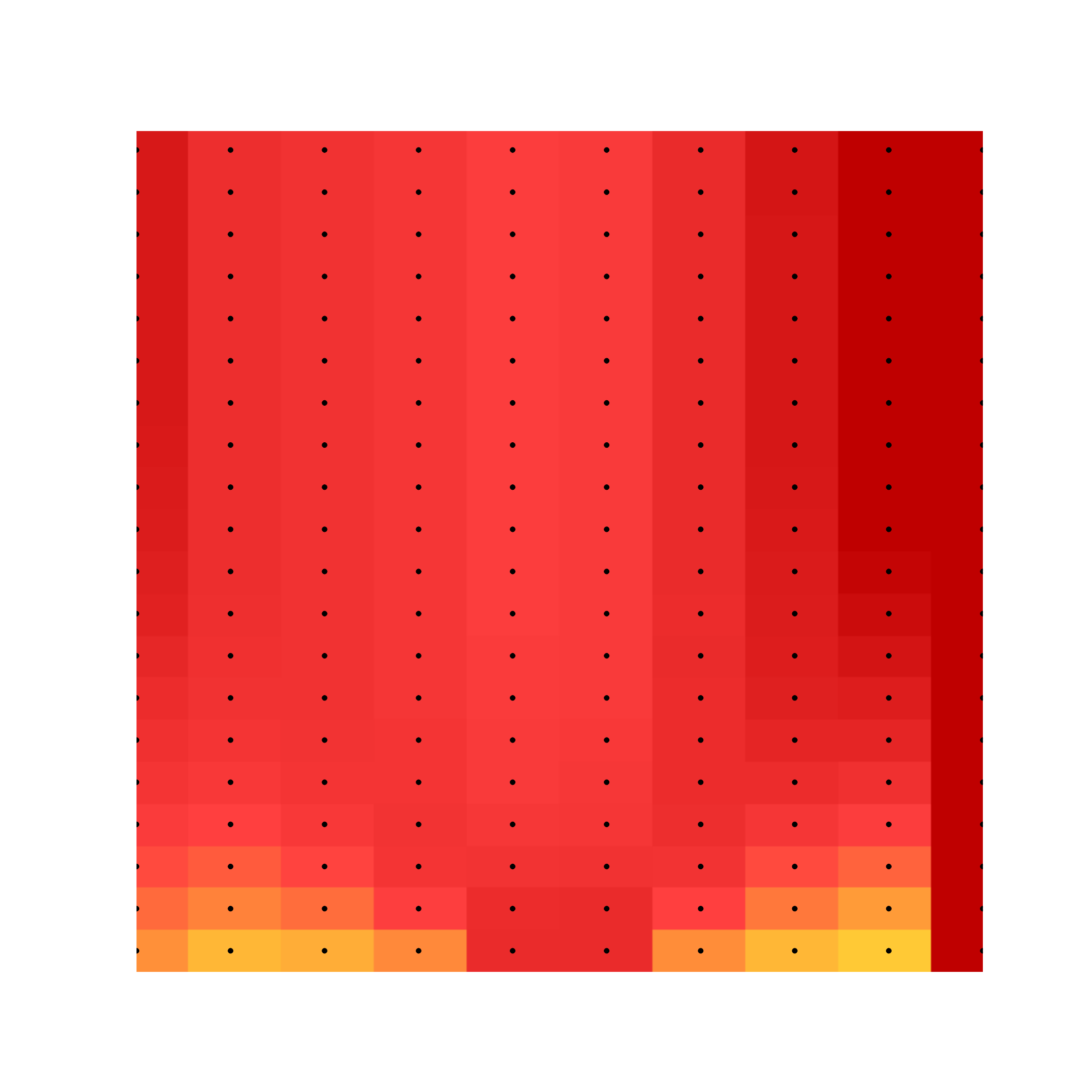}{0}{1}{}} &
\vspace{-.5cm}\resizebox{.24\textwidth}{!}{\maakmooieticks{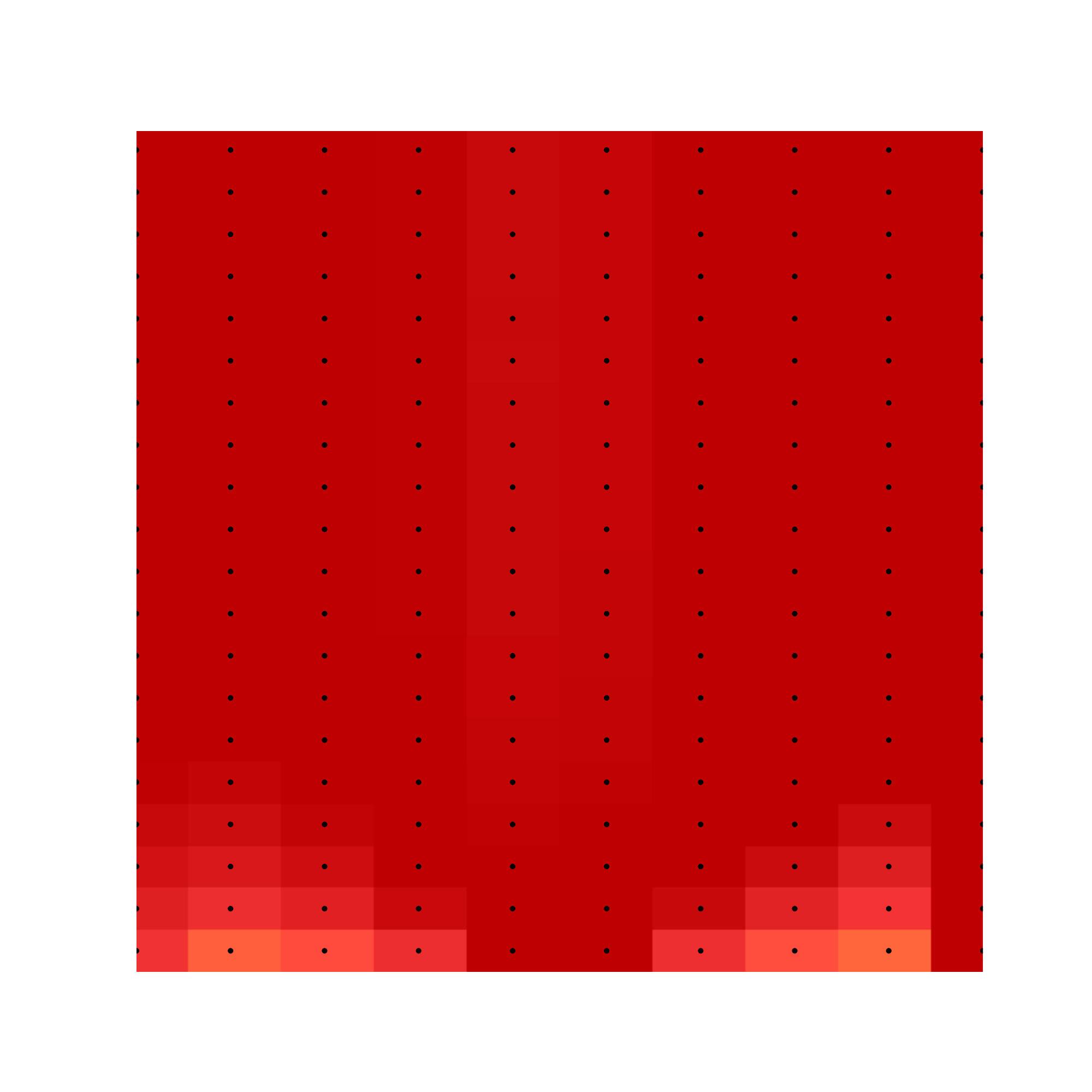}{0}{1}{}} &
\resizebox{!}{.24\textwidth}{\maakmooietickscb{figuren/improvement/improvement_colorbar.pdf}{}}
\end{tabular}
\caption{The best 1D1D mass source term estimation procedure and the potential gain when considering the statistical error. The first row presents a partition of the parameter domain based on the estimation procedure with the lowest statistical error. The dots represent the exact position in parameter space where the simulation occurred, crosses indicate inconclusive results. The second row presents the factor by which the standard deviation increases when the standard mass choice, \texttt{a\_tl}, is used instead of the best estimation procedure.}
\label{fig:ii_best_1D1D_mass_var}
\end{figure}

The effect of using continuous velocities (1D1D) instead of forward-backward scattering (1D0D) is clear when we compare the best mass source estimation procedure throughout the parameter domain based on statistical error for 1D1D in the first row of Figure~\ref{fig:ii_best_1D1D_mass_var} with the results for 1D0D from Figure~\ref{fig:best_1D0D_mass_stdv_bestest}. The \texttt{natl\_tl} estimation procedure, which benefited from the connection between travelled distance and score, loses terrain compared to other estimation procedures. This is most notable in the region close to $\VAPr=0.5$, which corresponds to high variance of the post-collisional velocity, and consequently additional variation on the travelled distance compared to the 1D0D case. Of special interest is the appearance of another estimation procedure, \texttt{nac\_tl}, under select circumstances. Further, the boundary between \texttt{natl\_ne} and \texttt{nac\_ne} shifts marginally in the advantage of \texttt{natl\_ne}.

The second row of Figure~\ref{fig:ii_best_1D1D_mass_var} shows the potential loss if the default estimation procedure, \texttt{a\_tl} is used instead of the best estimation procedure, and is to be compared with Figure~\ref{fig:best_1D0D_mass_stdv_gain} for 1D0D. Both sets of figures are very similar, with an exception to the isotropic situation, where the variance on the post-collisional velocity is much larger. There, \texttt{a\_tl} has deteriorated compared to the 1D0D situation. This deterioration can be understood from the same context as the deterioration of \texttt{natl\_tl}, since both \texttt{a\_tl} and \texttt{natl\_tl} benefit from the, now diminished, connection between travelled distance and scoring.

\begin{figure}[H]\centering
\begin{tabular}{m{.25\textwidth}m{.25\textwidth}m{.25\textwidth}m{.2\textwidth}}
\resizebox{.24\textwidth}{!}{\maakmooieticks{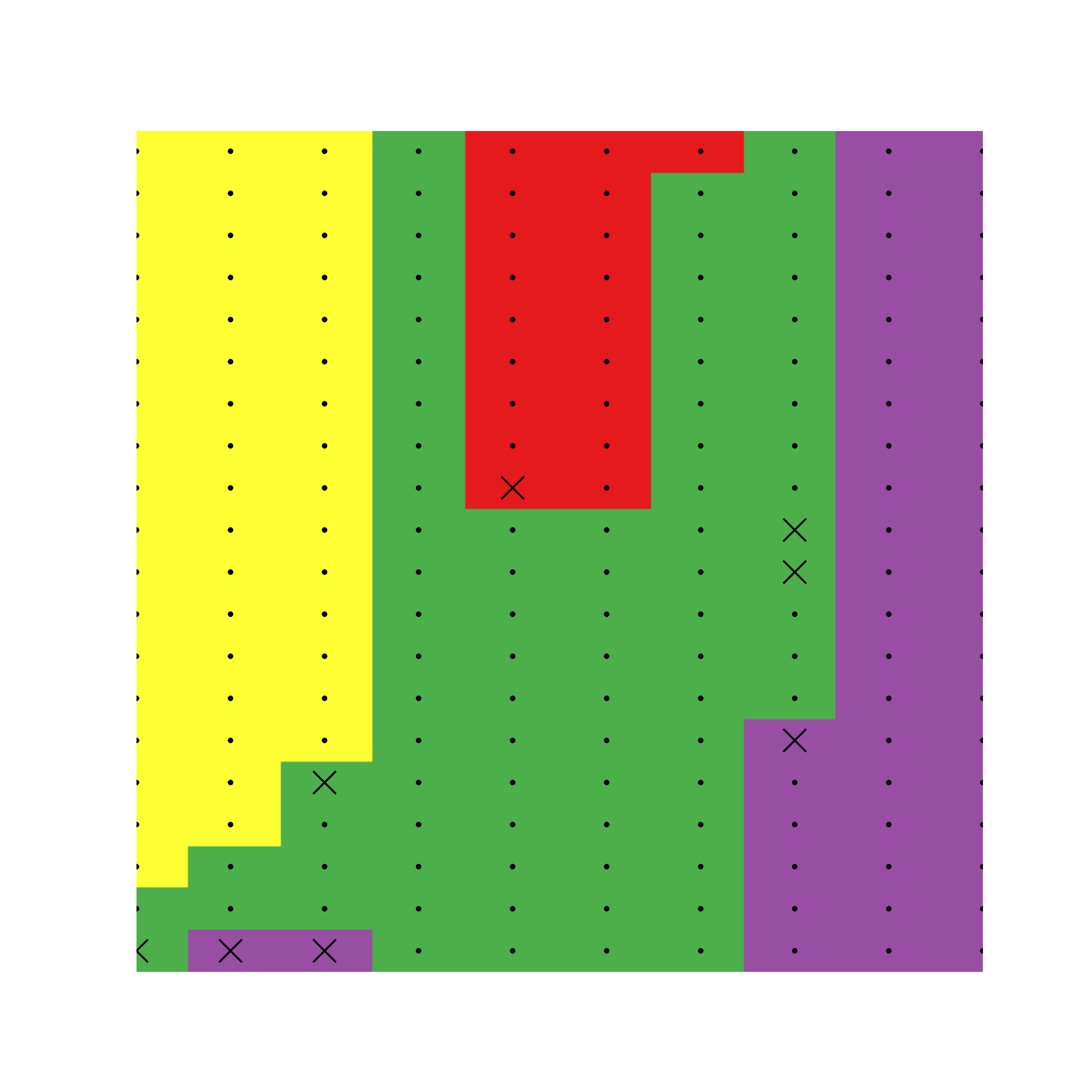}{1}{0}{$\dfrac{\VAcss}{\VAcst}=0.75$}} &
\resizebox{.24\textwidth}{!}{\maakmooieticks{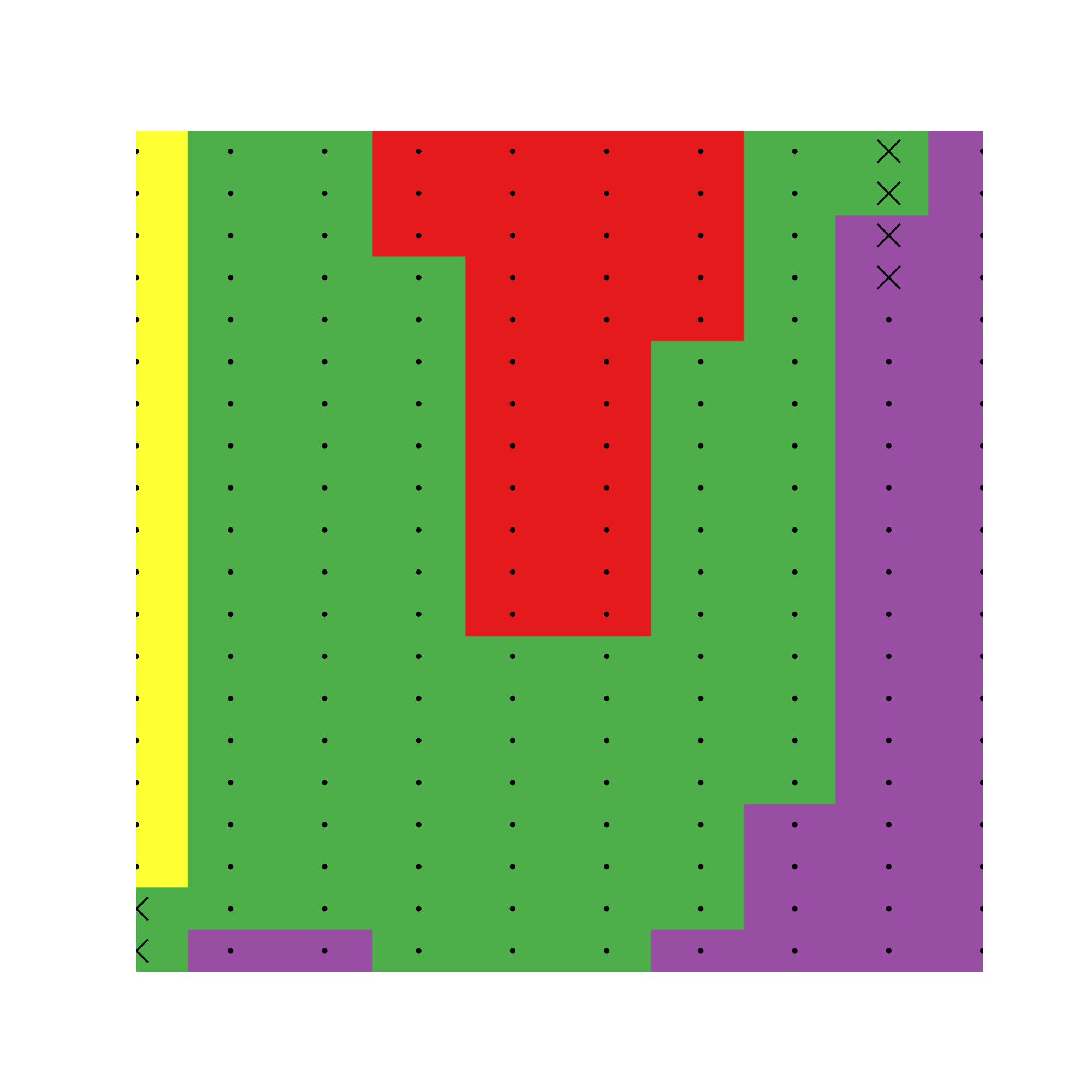}{0}{0}{$\dfrac{\VAcss}{\VAcst}=0.5$}} &
\resizebox{.24\textwidth}{!}{\maakmooieticks{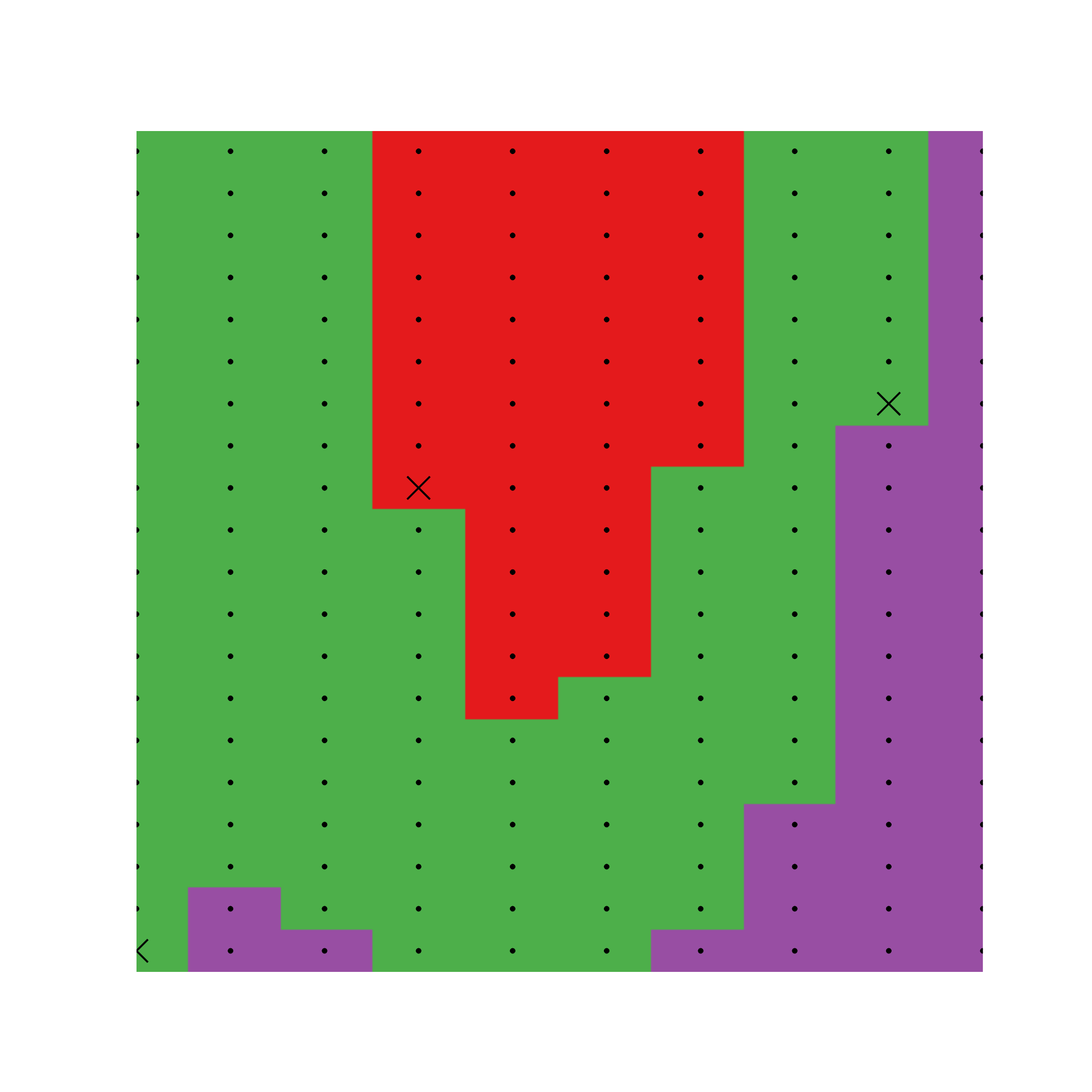}{0}{0}{$\dfrac{\VAcss}{\VAcst}=0.25$}} &
\localcolorlegendoneB
\end{tabular}\\
\begin{tabular}{b{.25\textwidth}b{.25\textwidth}b{.25\textwidth}b{.2\textwidth}}
\vspace{-.5cm}\resizebox{.24\textwidth}{!}{\maakmooieticks{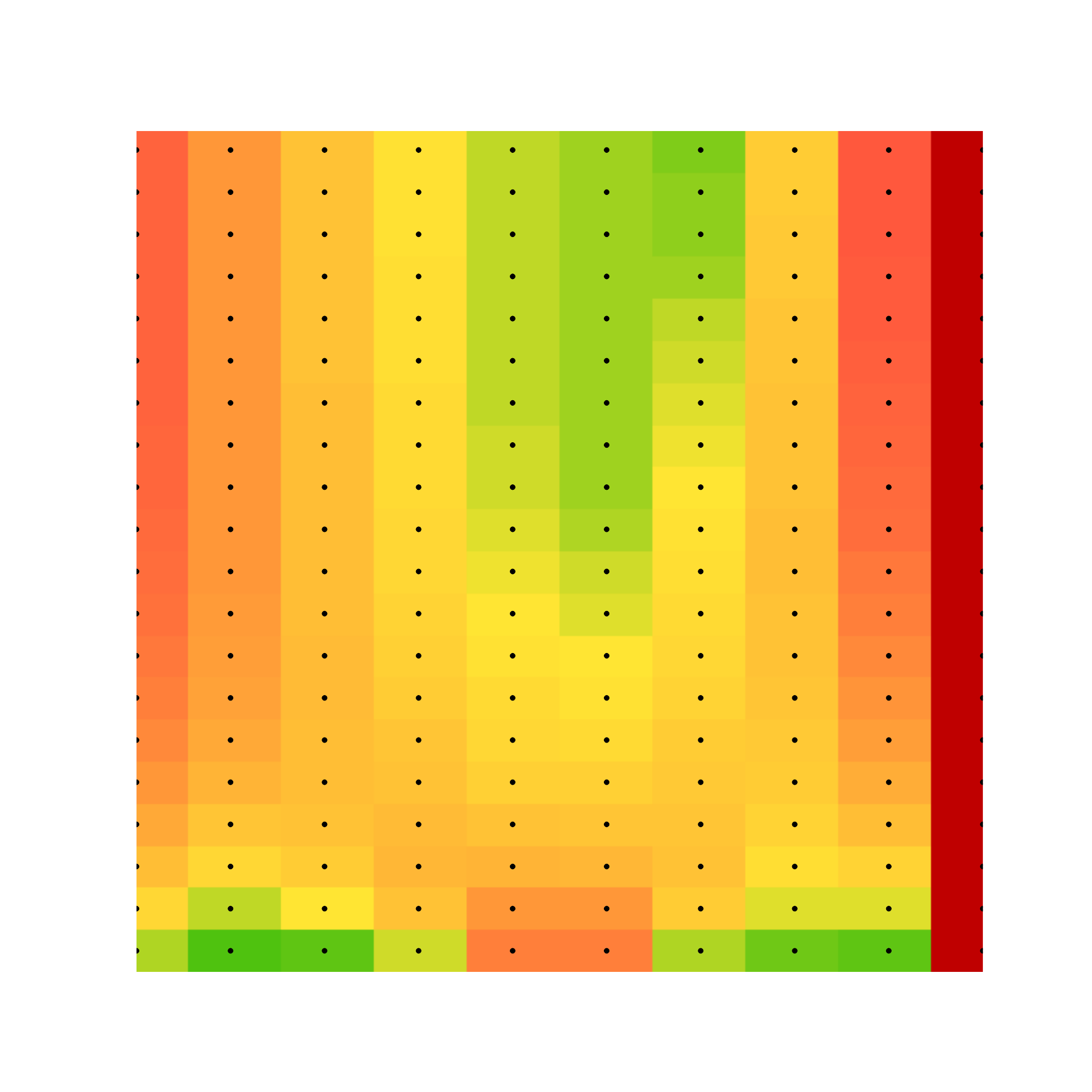}{1}{1}{}} &
\vspace{-.5cm}\resizebox{.24\textwidth}{!}{\maakmooieticks{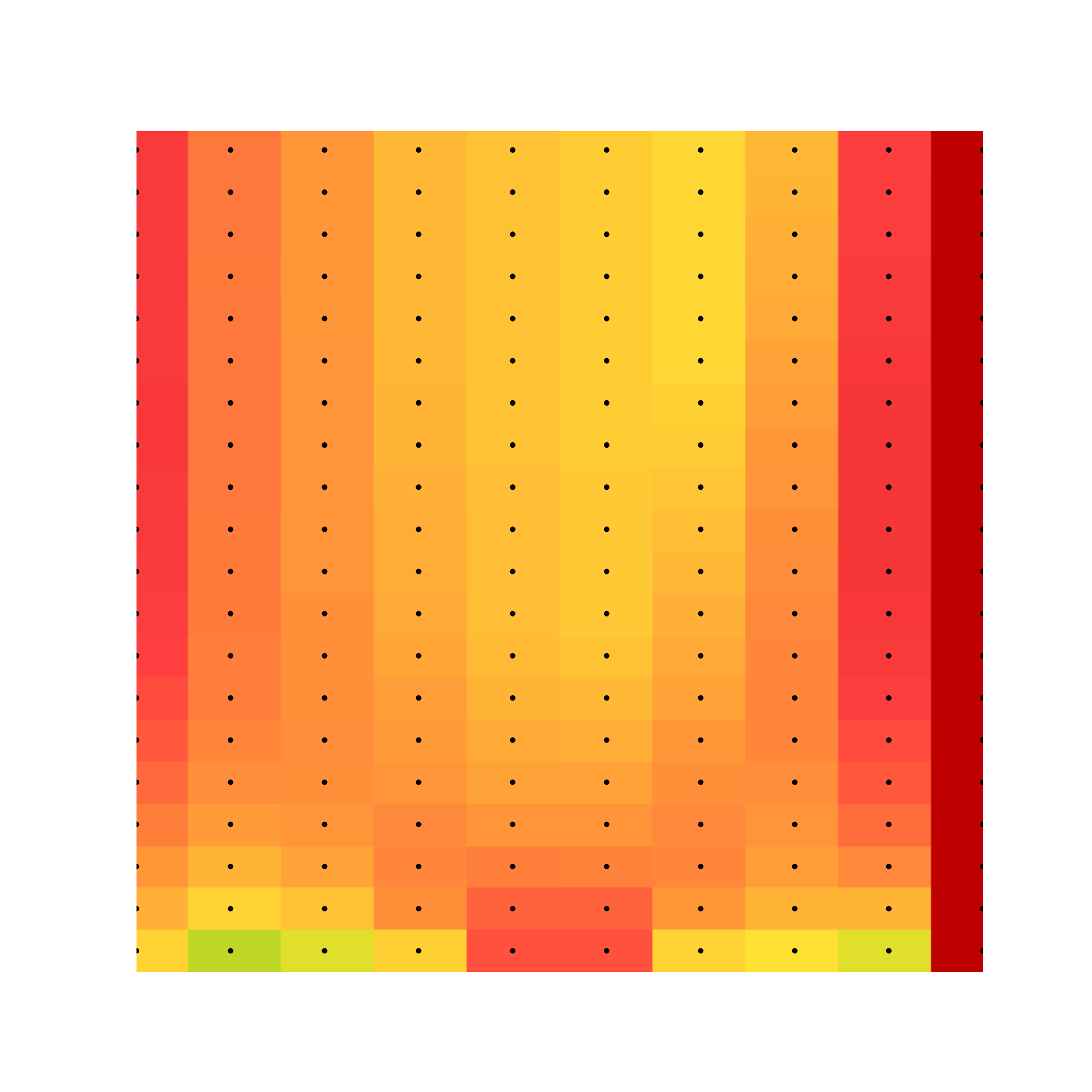}{0}{1}{}} &
\vspace{-.5cm}\resizebox{.24\textwidth}{!}{\maakmooieticks{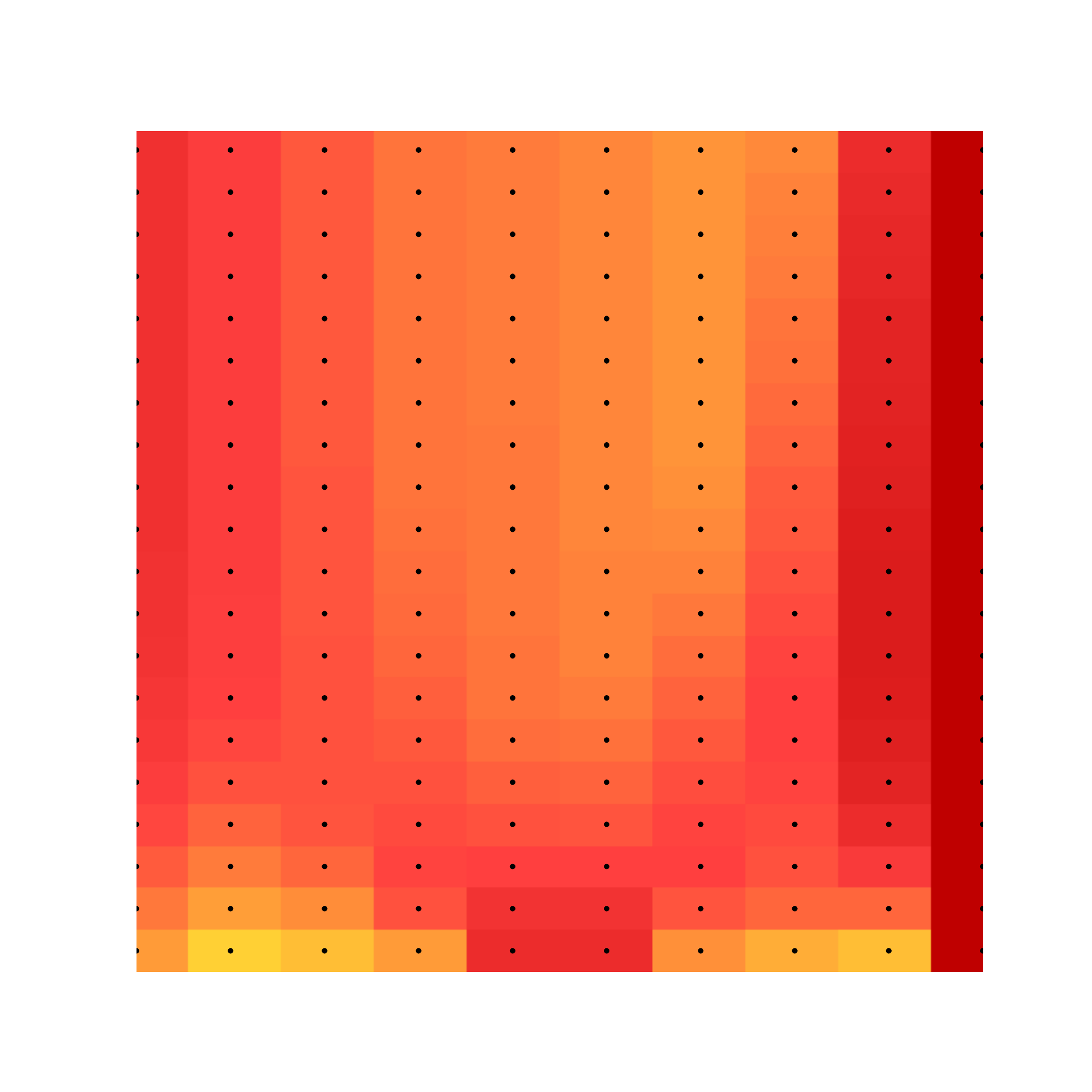}{0}{1}{}} &
\resizebox{!}{.24\textwidth}{\maakmooietickscb{figuren/improvement/improvement_colorbar.pdf}{}}
\end{tabular}
\caption{The best 1D1D mass source term estimation procedure and the potential gain when considering computational cost for a given statistical error. The first row presents a partition of the parameter domain based on the estimation procedure with the lowest computational cost error. The dots represent the exact position in parameter space where the simulation occurred, crosses indicate inconclusive results. The second row presents the factor by which the computational cost increases when the standard mass choice, \texttt{a\_tl}, is used instead of the best estimation procedure.}
\label{fig:iif_1D1D_mass_cost}
\end{figure}

The results when considering cost (variance times number of collisions) are shown in Figure~\ref{fig:iif_1D1D_mass_cost}. As was the case for 1D0D in Figure~\ref{fig:best_1D0D_mass_cost}, going from variance as a measure of performance to cost is advantageous for \texttt{a} and \texttt{natl} simulations. The most notable effect is the increased domain of \texttt{a\_ne} in the region of large variance on the post-collisional velocity (around $\VAPr$). When the variance on the post-collisional velocity is larger, the positive effect of more scoring events on the variance is mitigated somehow, so executing an absorption is less adverse.

The important conclusion from the second row of Figure~\ref{fig:best_1D0D_mass_cost}, that the \texttt{a\_tl} estimation procedure is sufficiently good in the relevant high-collisional isotropic case when the survival probability is high, is retained in the results of Figure~\ref{fig:iif_1D1D_mass_cost}. The only large difference with the 1D0D results is the deterioration of the \texttt{a\_tl} estimation procedure in the low-collisional isotropic case. This follows from the similar detoriation of the standard deviation visible when comparing Figure~\ref{fig:best_1D0D_mass_stdv_gain} and the second row of Figure~\ref{fig:ii_best_1D1D_mass_var}.

An important additional conclusion from this section is the large parallel between the results for the 1D0D case and for the 1D1D case. The conclusion that the \texttt{a\_tl} estimation procedure functions well in fusion-relevant case, is retained, as well as most of the trade-off between \texttt{a\_tl} and the best estimation procedure. Consequently also the conclusion that in regions of the domain with deviating parameter values, different methods should be used is retained.

\subsection{Momentum source estimation in 1D1D\label{subsec:iif_1D1D_mom}}

\begin{figure}\centering
\begin{tabular}{m{.25\textwidth}m{.25\textwidth}m{.25\textwidth}m{.2\textwidth}}
\resizebox{.24\textwidth}{!}{\maakmooieticks{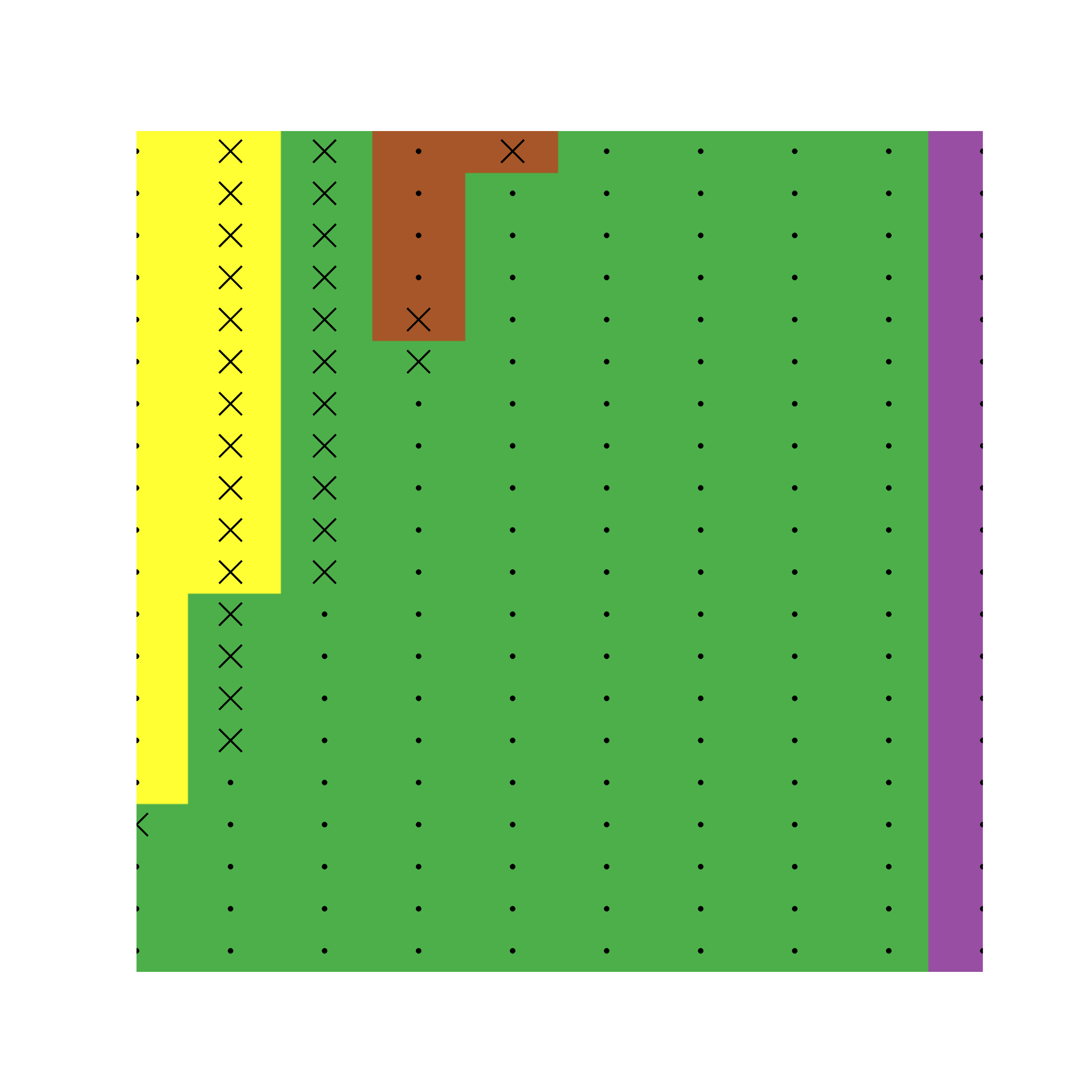}{1}{0}{$\dfrac{\VAcss}{\VAcst}=0.98$}} &
\resizebox{.24\textwidth}{!}{\maakmooieticks{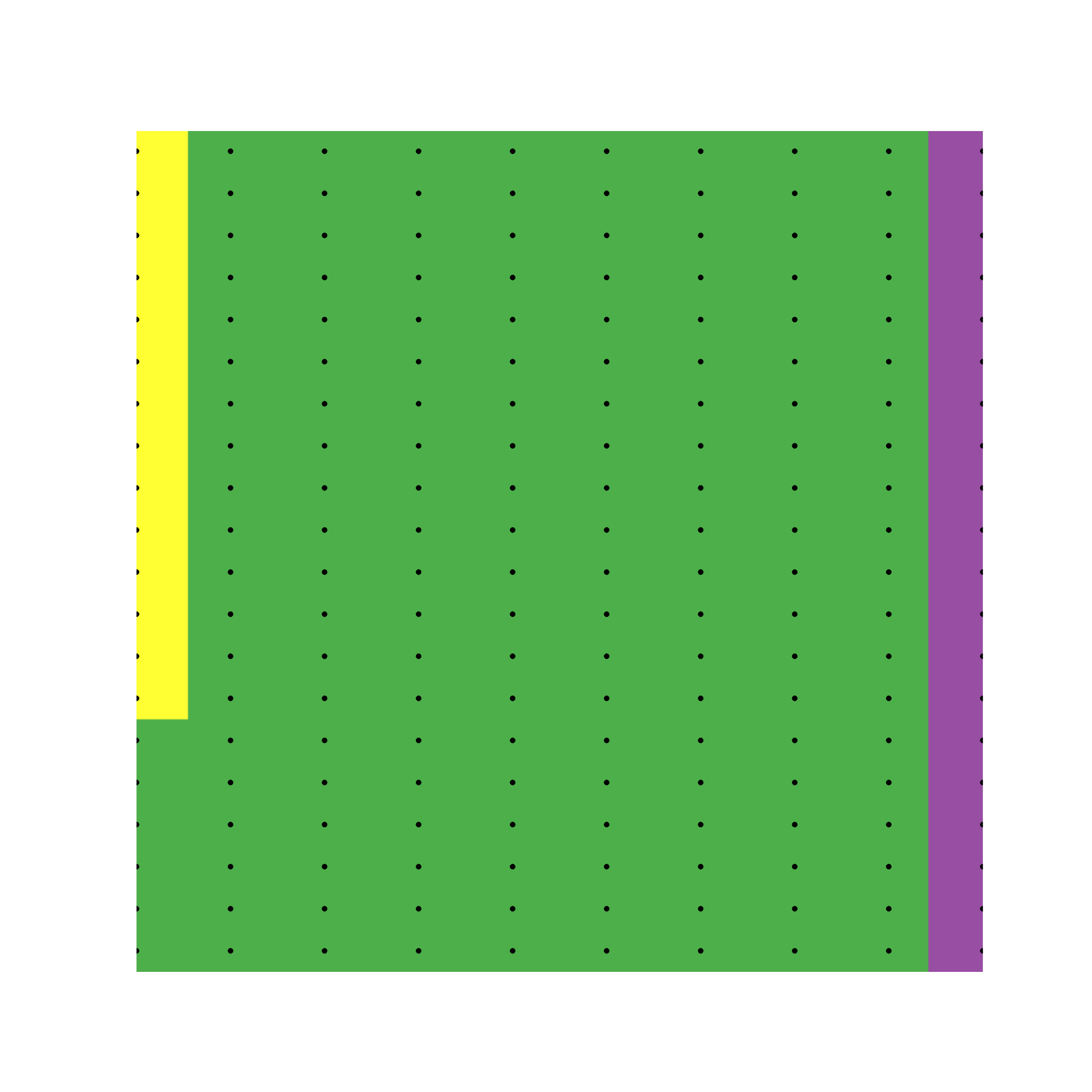}{0}{0}{$\dfrac{\VAcss}{\VAcst}=0.5$}} &
\resizebox{.24\textwidth}{!}{\maakmooieticks{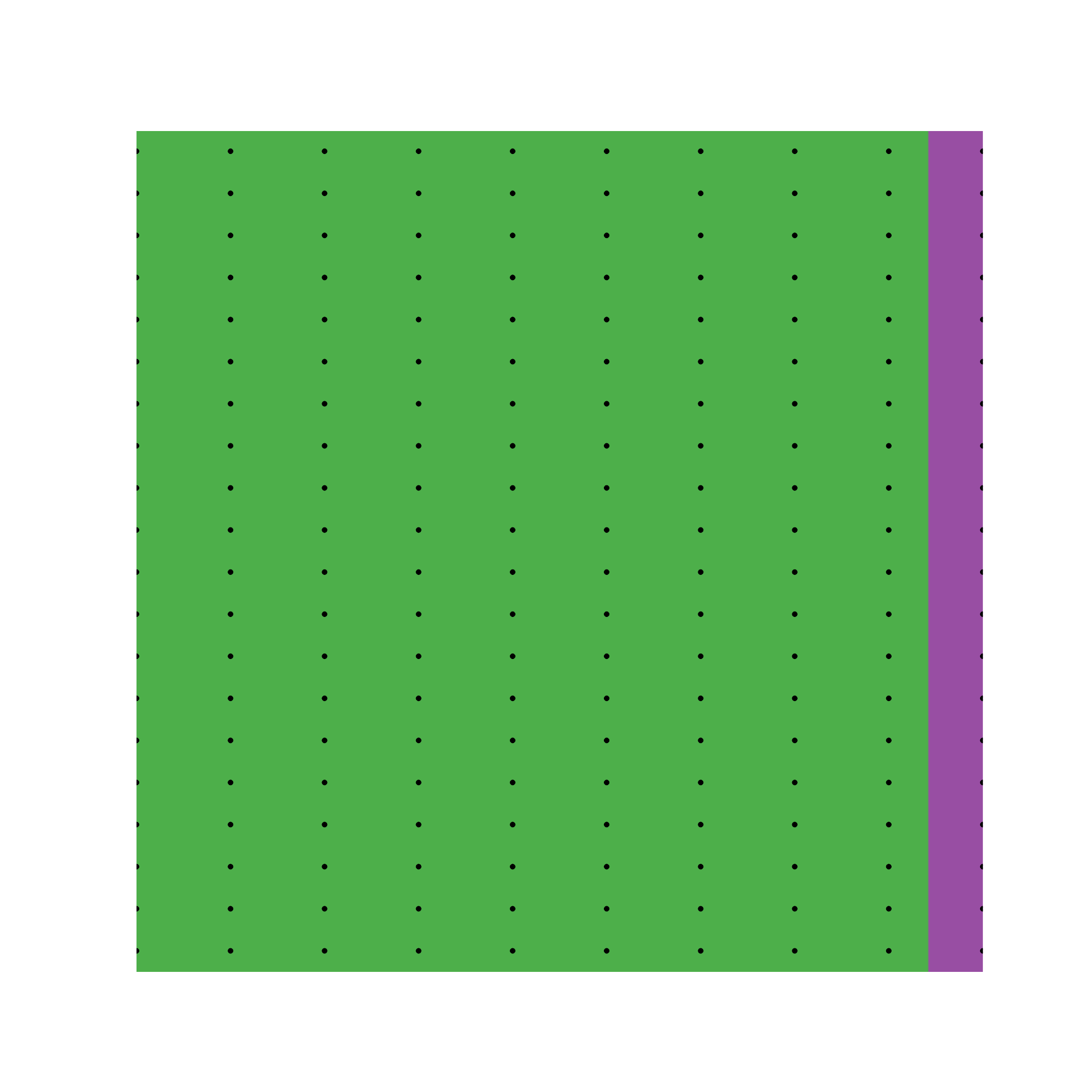}{0}{0}{$\dfrac{\VAcss}{\VAcst}=0.25$}} &
\localcolorlegendoneC\\
\end{tabular}\\
\begin{tabular}{b{.25\textwidth}b{.25\textwidth}b{.25\textwidth}b{.2\textwidth}}
\vspace{-0.5cm}\resizebox{.24\textwidth}{!}{\maakmooieticks{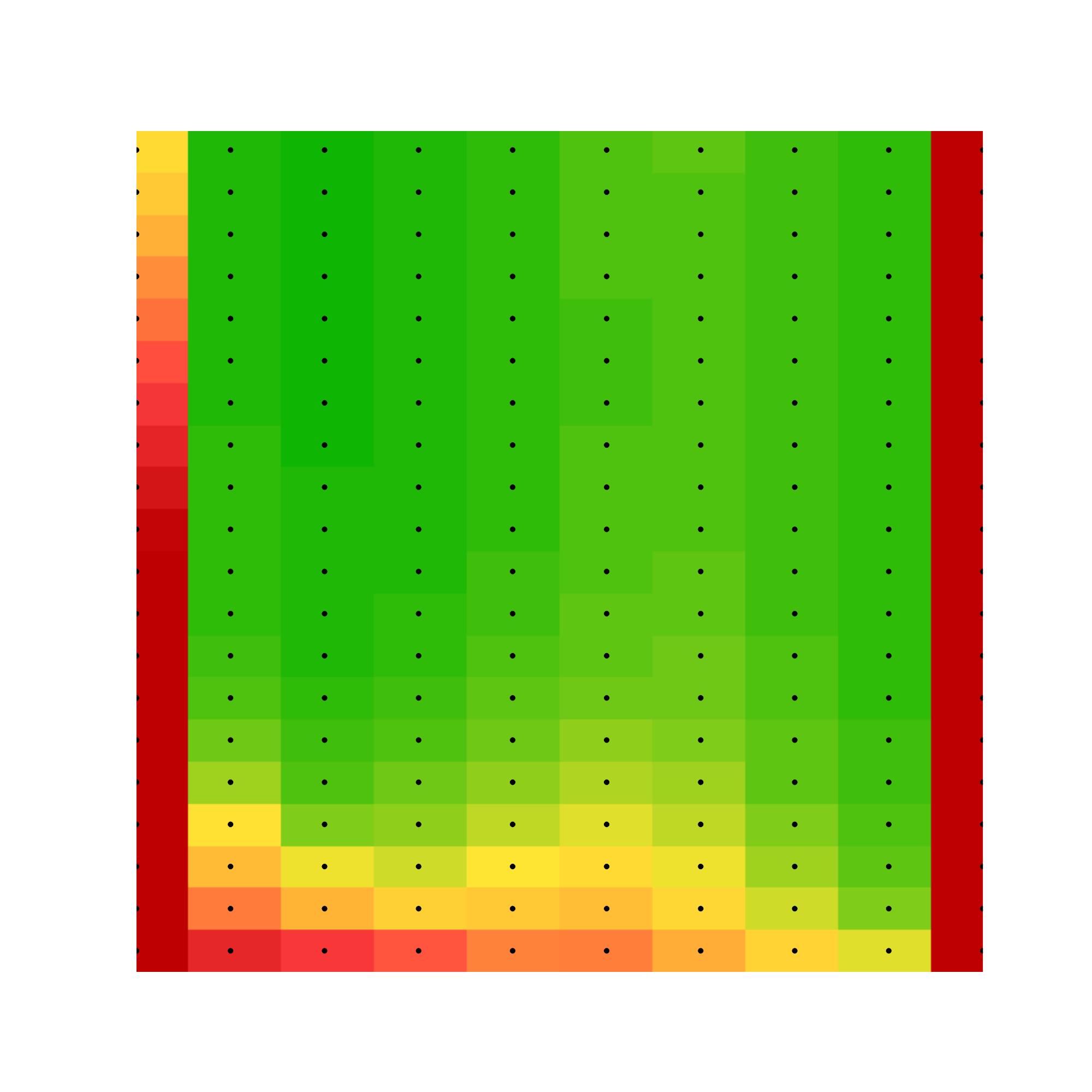}{1}{1}{}} &
\vspace{-0.5cm}\resizebox{.24\textwidth}{!}{\maakmooieticks{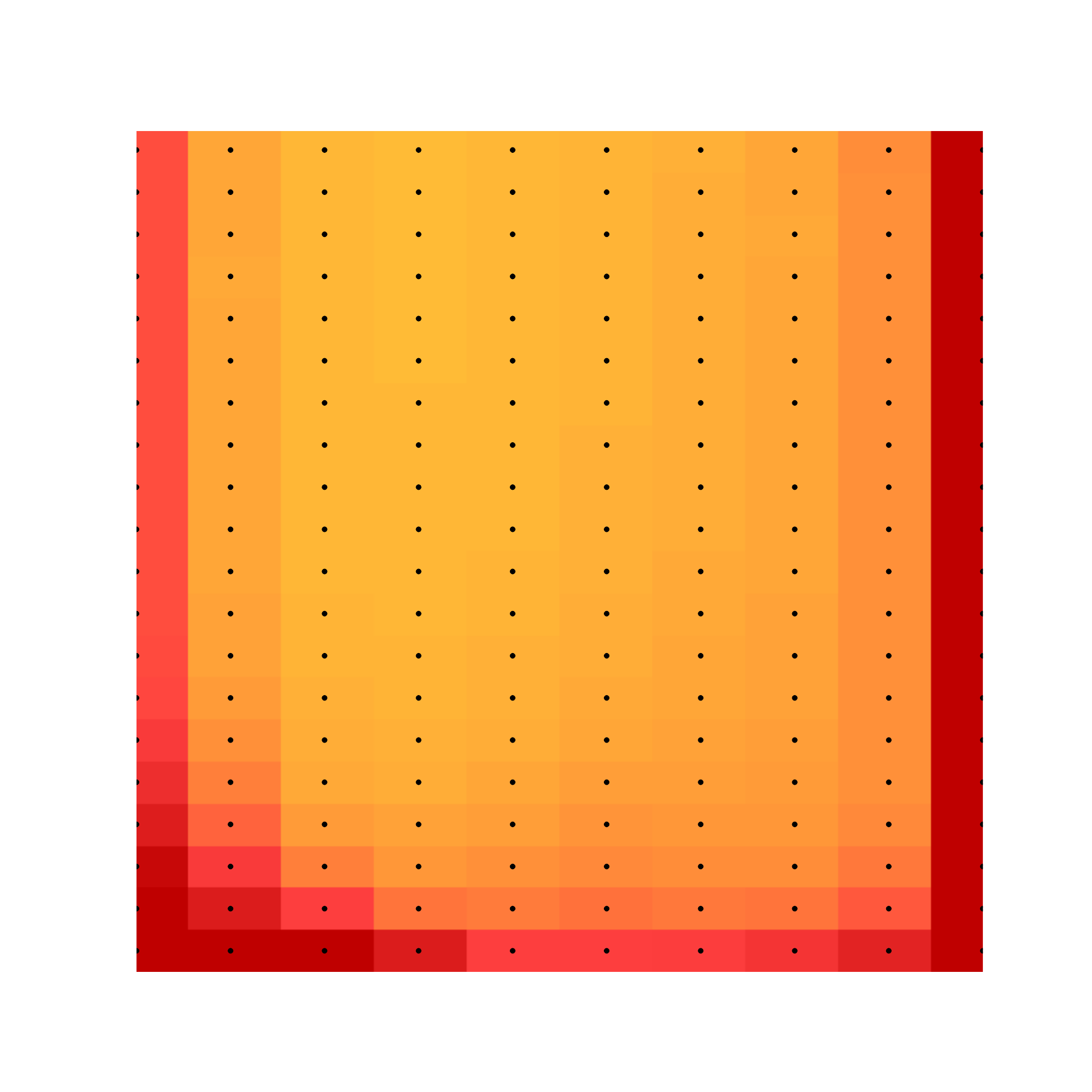}{0}{1}{}} &
\vspace{-0.5cm}\resizebox{.24\textwidth}{!}{\maakmooieticks{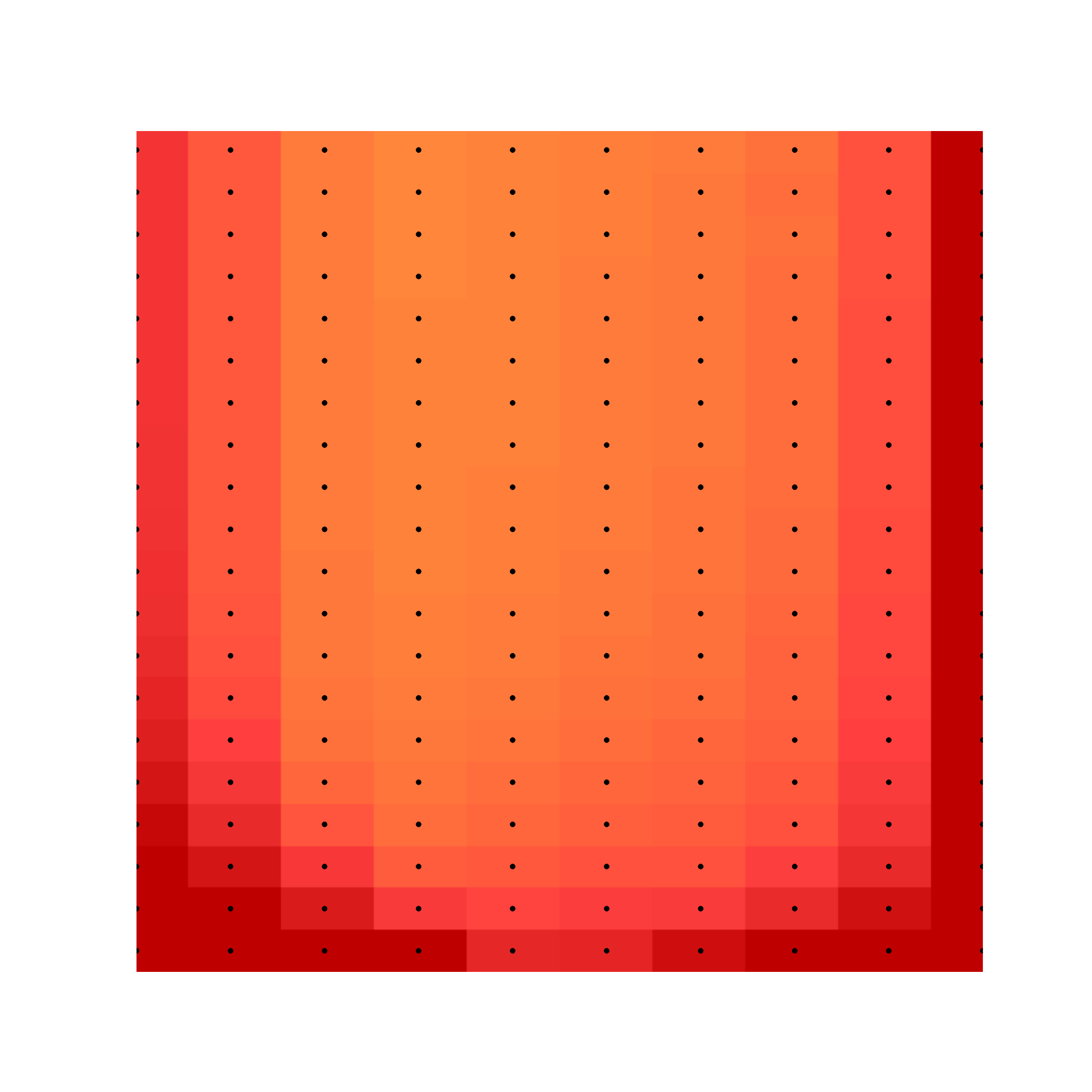}{0}{1}{}} &
\resizebox{!}{.24\textwidth}{\maakmooietickscb{figuren/improvement/improvement_colorbar.pdf}{}}
\end{tabular}
\caption{The best 1D1D momentum source term estimation procedure and the potential gain when considering the statistical error. The first row presents a partition of the parameter domain based on the estimation procedure with the lowest statistical error. The dots represent the exact position in parameter space where the simulation occurred, crosses indicate inconclusive results. The second row presents the factor by which the standard deviation increases when the standard mass choice, \texttt{a\_tl}, is used instead of the best estimation procedure.}
\label{fig:iif_1D1D_mom_var}
\end{figure}

\begin{figure}\centering
\begin{tabular}{m{.25\textwidth}m{.25\textwidth}m{.25\textwidth}m{.2\textwidth}}
\resizebox{.24\textwidth}{!}{\maakmooieticks{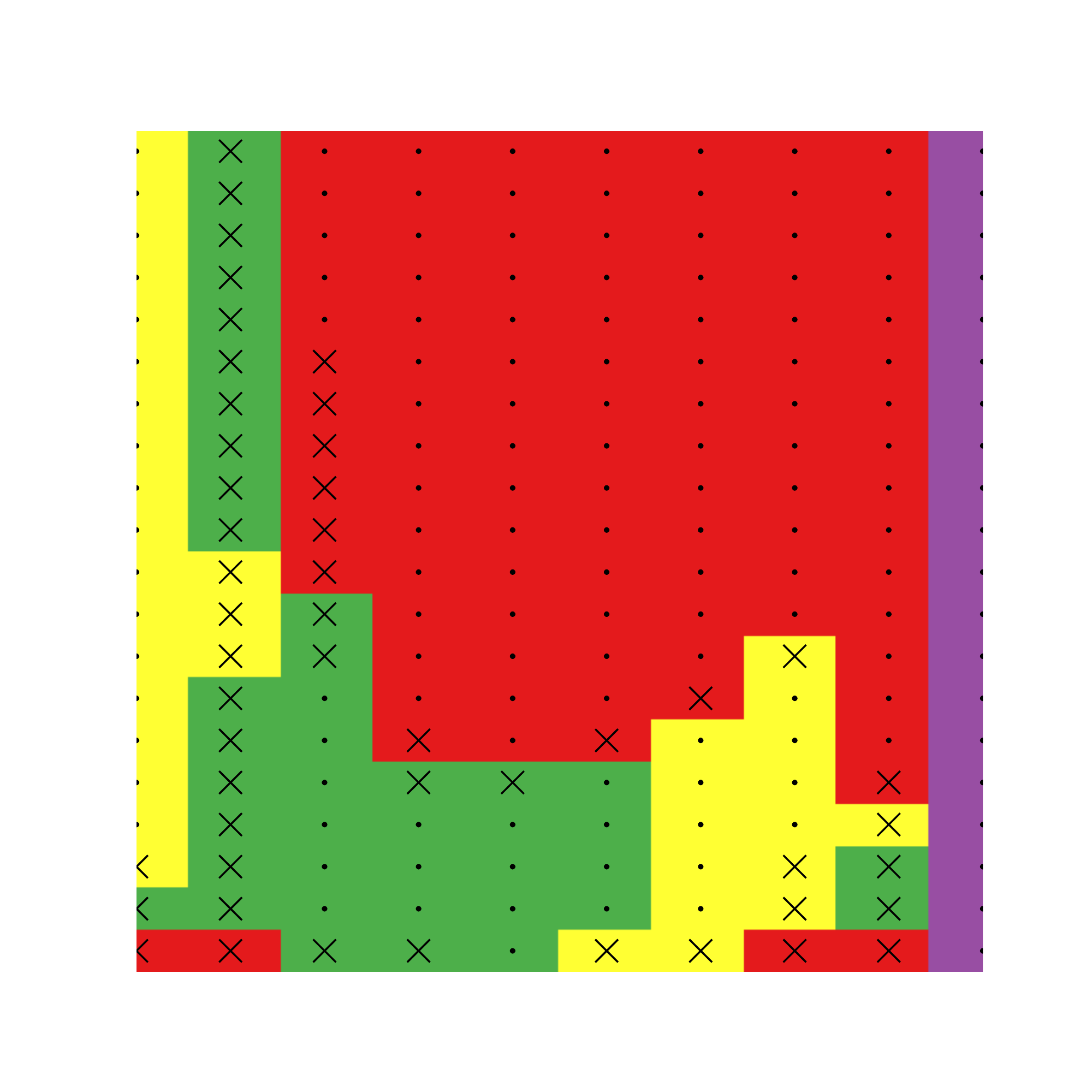}{1}{0}{$\dfrac{\VAcss}{\VAcst}=0.94$}} &
\resizebox{.24\textwidth}{!}{\maakmooieticks{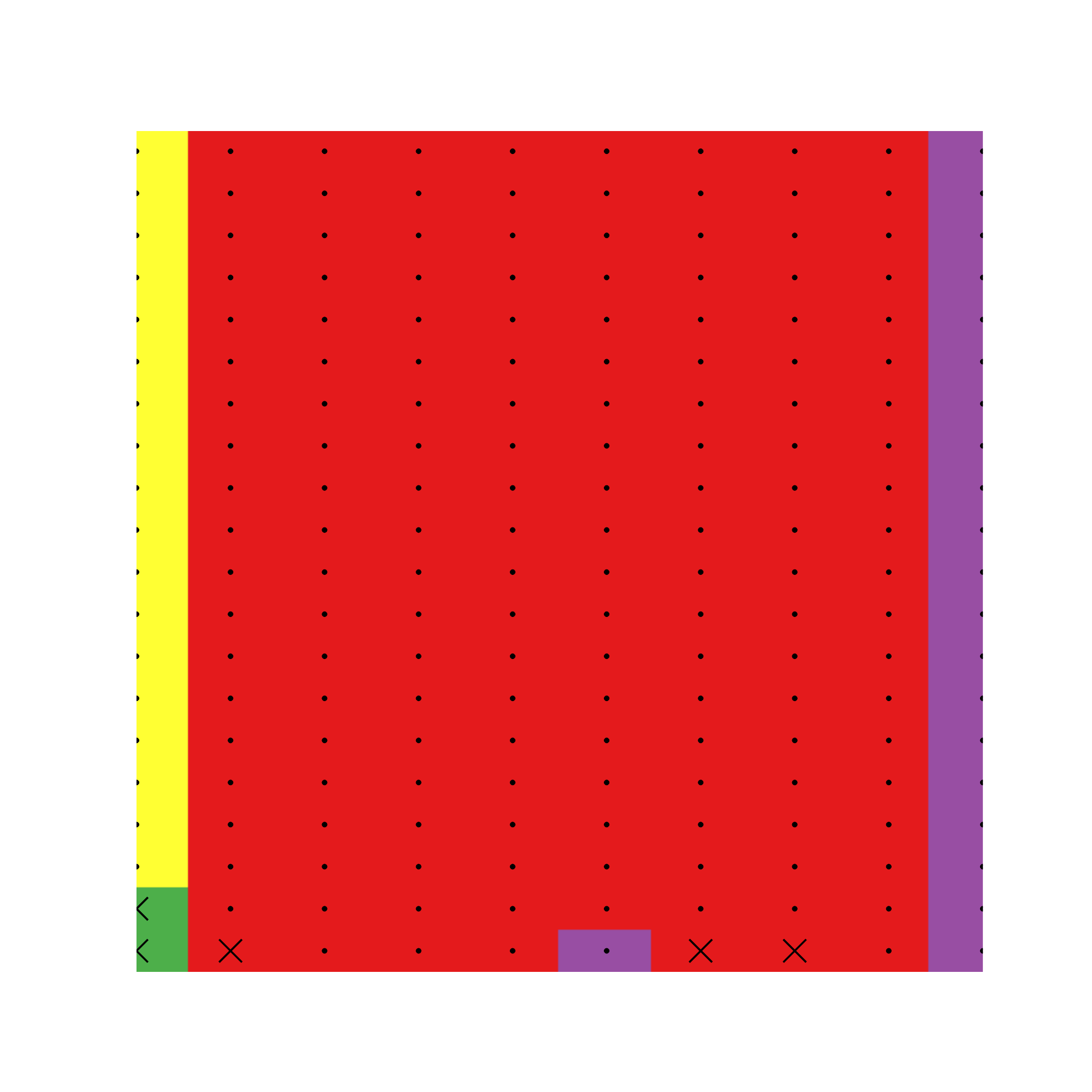}{0}{0}{$\dfrac{\VAcss}{\VAcst}=0.5$}} &
\resizebox{.24\textwidth}{!}{\maakmooieticks{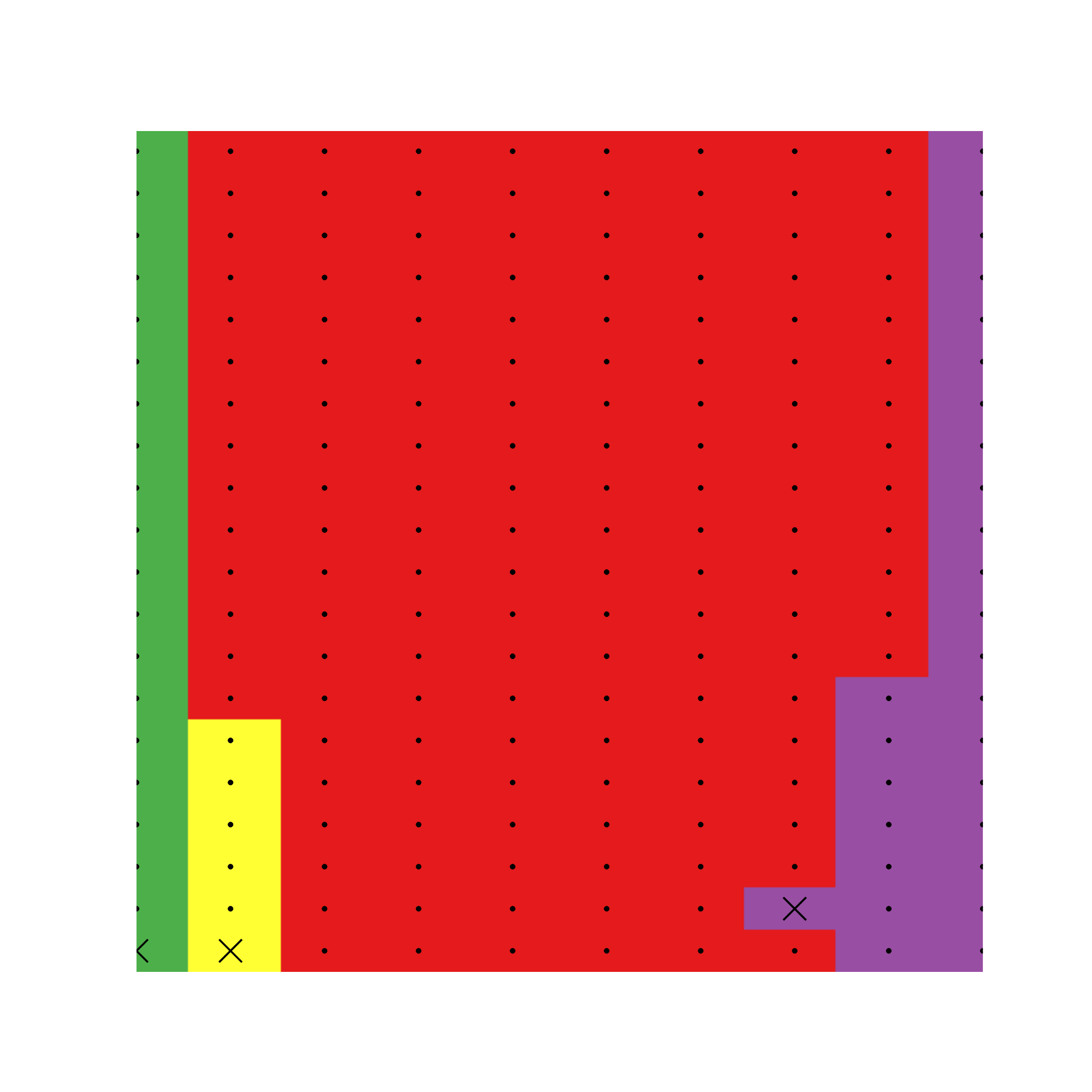}{0}{0}{$\dfrac{\VAcss}{\VAcst}=0.25$}} &
\localcolorlegendoneB\\
\end{tabular}\\
\begin{tabular}{b{.25\textwidth}b{.25\textwidth}b{.25\textwidth}b{.2\textwidth}}
\vspace{-0.5cm}\resizebox{.24\textwidth}{!}{\maakmooieticks{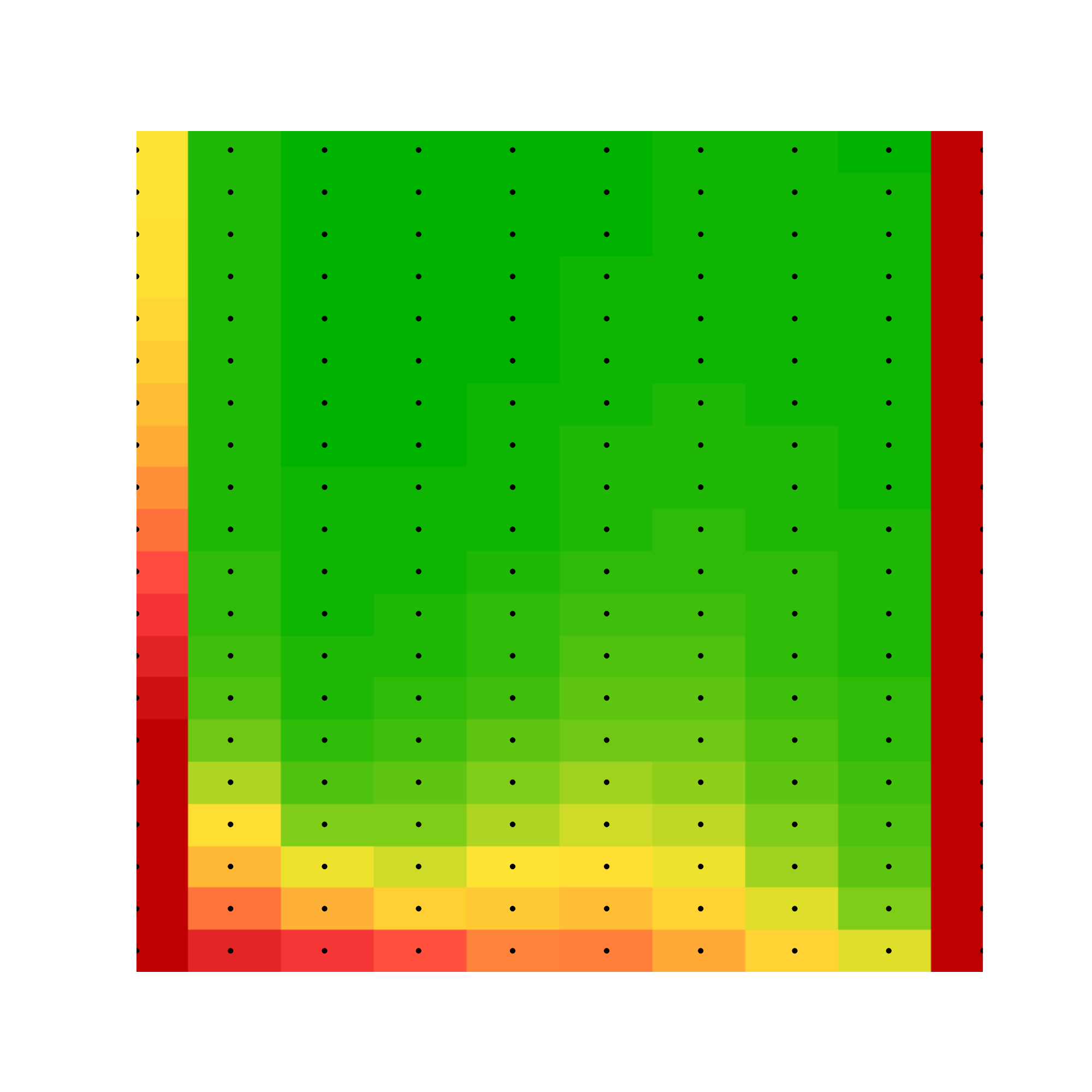}{1}{1}{}} &
\vspace{-0.5cm}\resizebox{.24\textwidth}{!}{\maakmooieticks{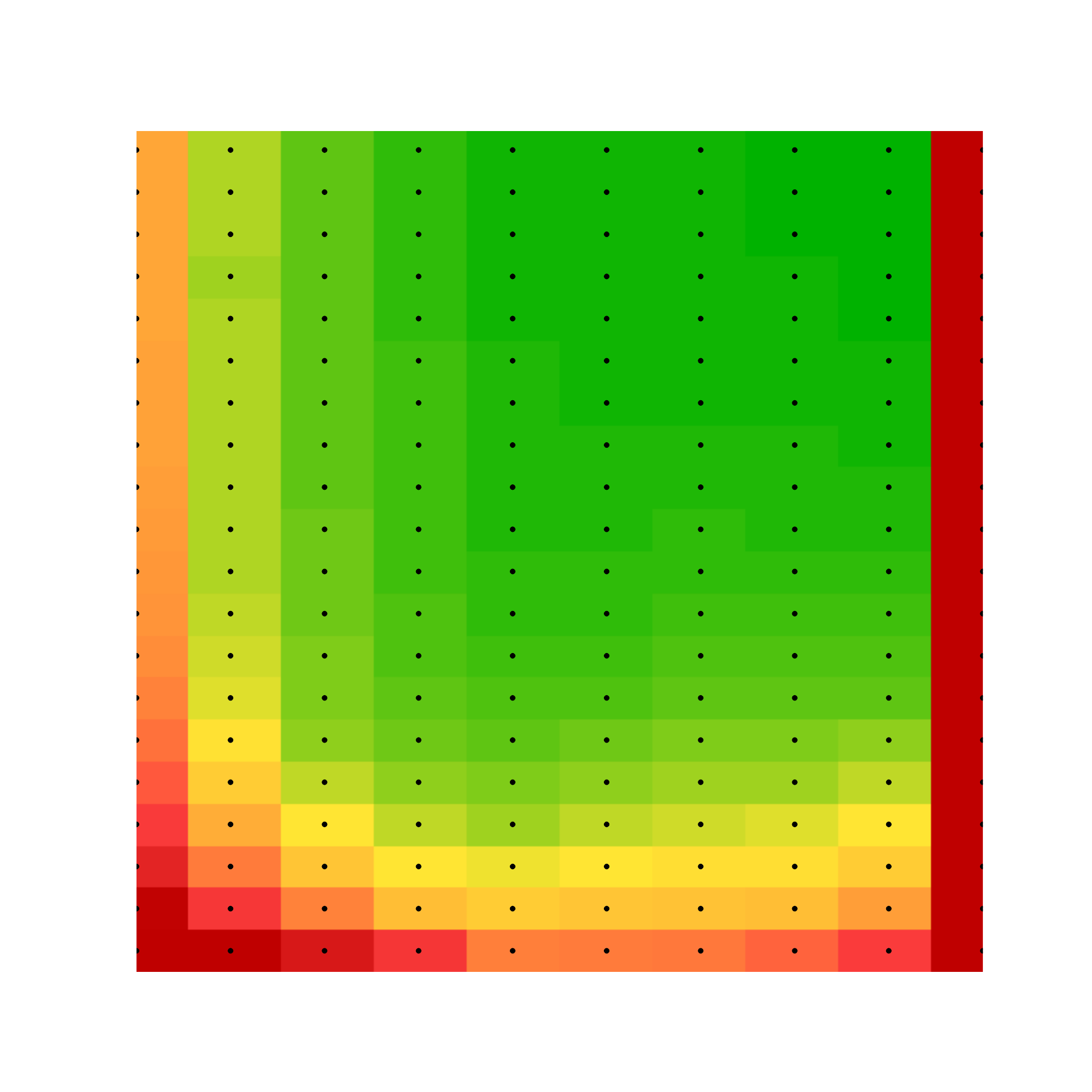}{0}{1}{}} &
\vspace{-0.5cm}\resizebox{.24\textwidth}{!}{\maakmooieticks{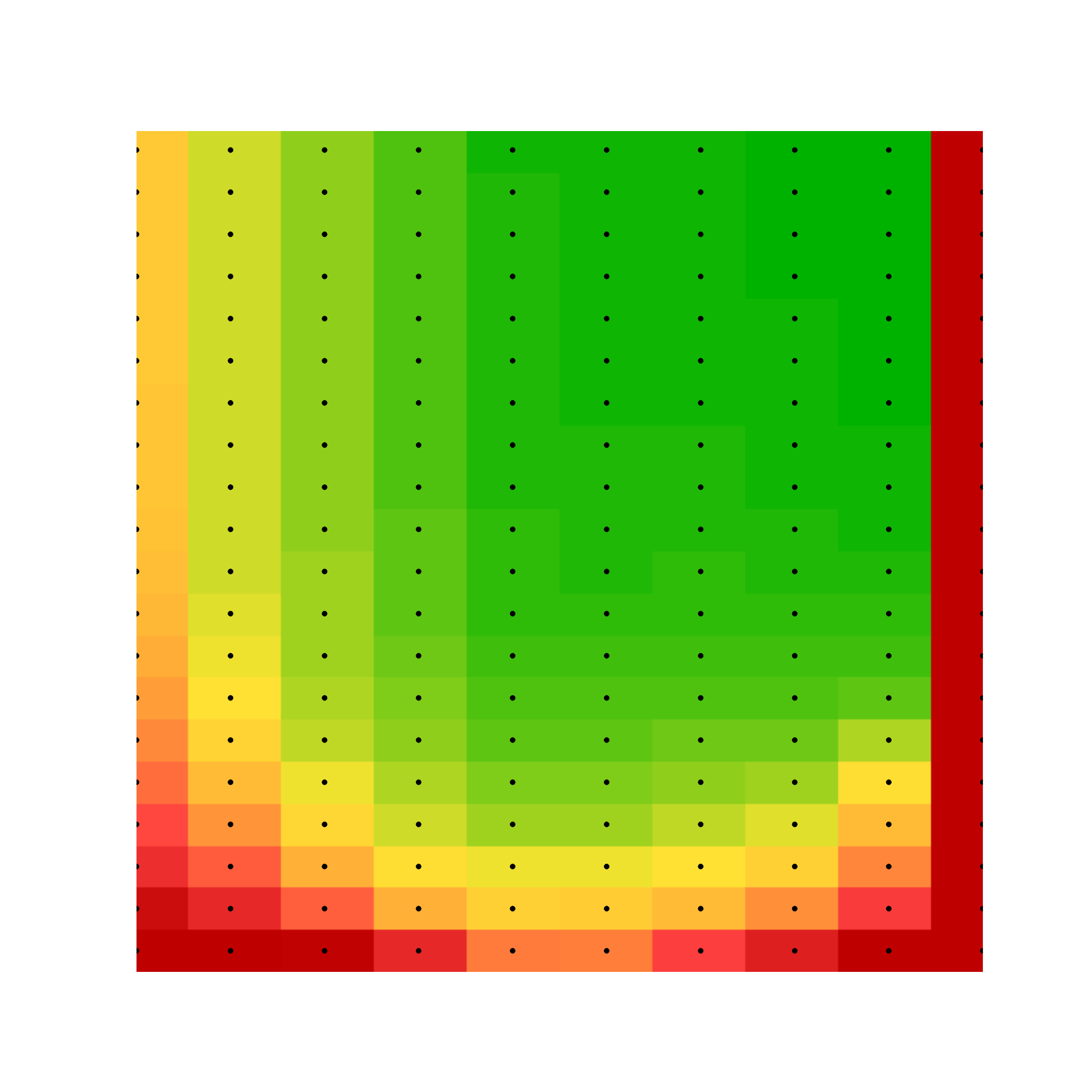}{0}{1}{}} &
\resizebox{!}{.24\textwidth}{\maakmooietickscb{figuren/improvement/improvement_colorbar.pdf}{}}
\end{tabular}
\caption{The best 1D1D momentum source term estimation procedure and the potential gain when considering computational cost for a given statistical error. The first row presents a partition of the parameter domain based on the estimation procedure with the lowest computational cost error. The dots represent the exact position in parameter space where the simulation occurred, crosses indicate inconclusive results. The second row presents the factor by which the computational cost increases when the standard mass choice, \texttt{a\_tl}, is used instead of the best estimation procedure.}
\label{fig:iif_1D1D_mom_cost}
\end{figure}

We have numerically extended our analysis to momentum source estimation. The results presenting the best momentum estimation procedure for the 1D1D setting are shown in Figure~\ref{fig:iif_1D1D_mom_var} for statistical error and in Figure~\ref{fig:iif_1D1D_mom_cost} for computational cost. The default momentum source estimation procedure in the EIRENE code is currently \texttt{a\_c}, which was selected via trial and error.

When we compare the best momentum estimation procedure to the best mass estimation procedure, we note the presence of a collision estimator \texttt{nac\_c}, and the much larger prevalence of the \texttt{nac\_ne} estimator when considering statistical error in the first row of Figure~\ref{fig:iif_1D1D_mom_var} and the \texttt{a\_ne} estimator when considering computational cost in the first row of Figure~\ref{fig:iif_1D1D_mom_cost}. As can be seen in the second row of Figure~\ref{fig:iif_1D1D_mom_cost}, the default momentum estimation procedure, \texttt{a\_c}, now functions very well for the largest part of the parameter domain when considering the computational cost, with notable exceptions when the background is highly anisotropic ($\VAPr\approx0$ or $\VAPr\approx1$) or for very low collisionality.

\begin{figure}\centering
\begin{tabular}{m{.25\textwidth}m{.25\textwidth}m{.25\textwidth}m{.2\textwidth}}
\resizebox{.24\textwidth}{!}{\maakmooieticks{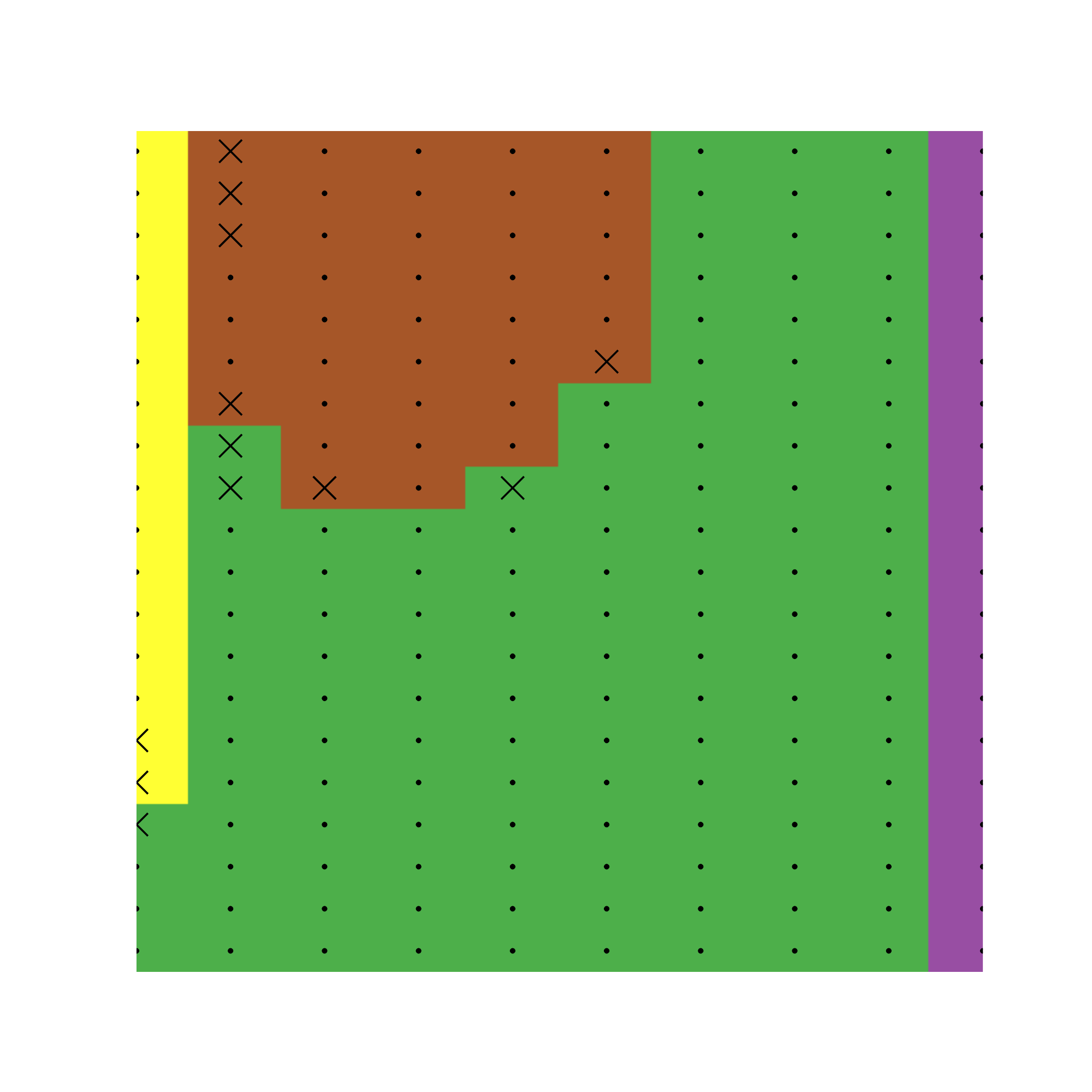}{1}{0}{$\dfrac{\VAcss}{\VAcst}=0.98$}} &
\resizebox{.24\textwidth}{!}{\maakmooieticks{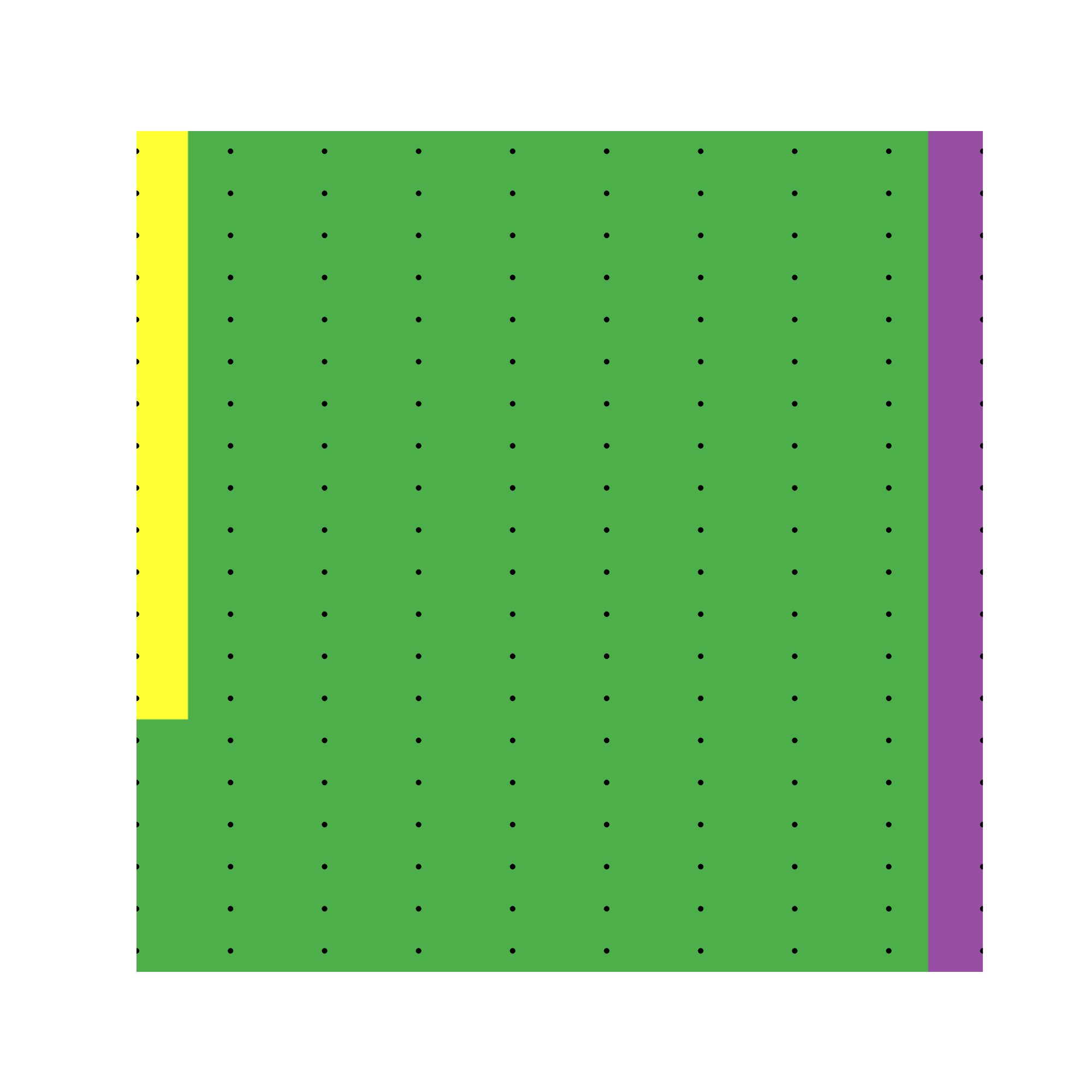}{0}{0}{$\dfrac{\VAcss}{\VAcst}=0.5$}} &
\resizebox{.24\textwidth}{!}{\maakmooieticks{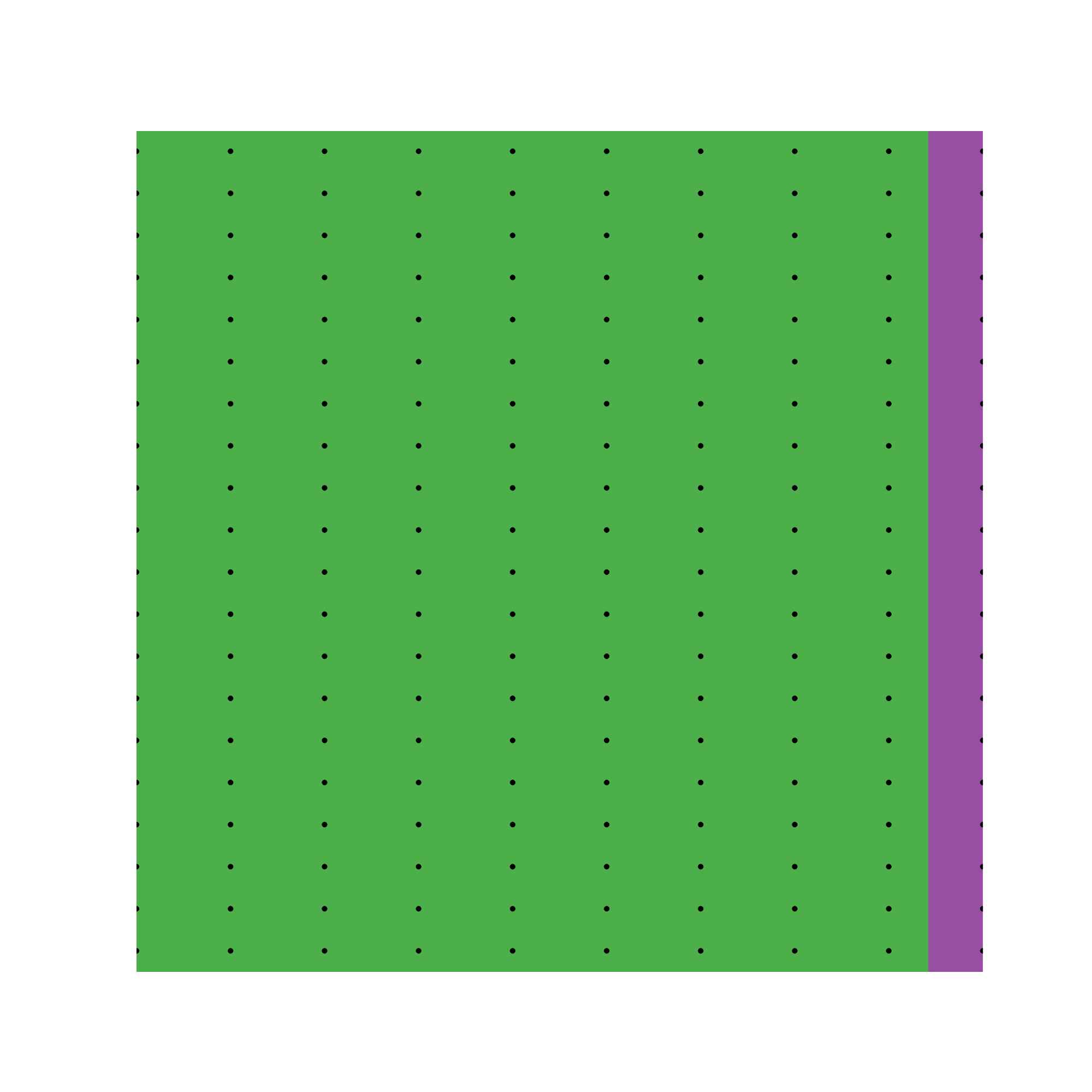}{0}{0}{$\dfrac{\VAcss}{\VAcst}=0.25$}} &
\localcolorlegendoneC\\
\end{tabular}\\
\begin{tabular}{b{.25\textwidth}b{.25\textwidth}b{.25\textwidth}b{.2\textwidth}}
\vspace{-0.5cm}\resizebox{.24\textwidth}{!}{\maakmooieticks{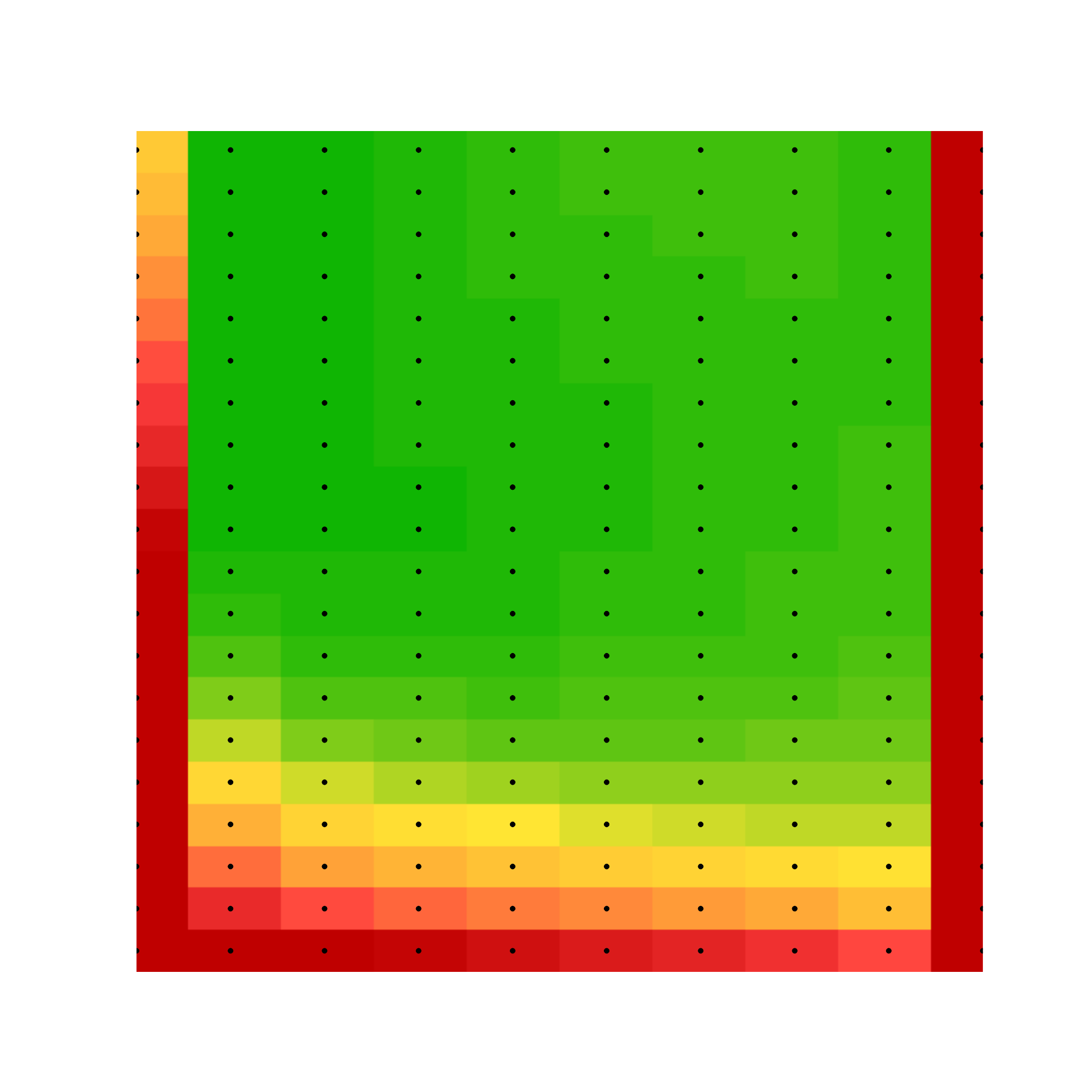}{1}{1}{}} &
\vspace{-0.5cm}\resizebox{.24\textwidth}{!}{\maakmooieticks{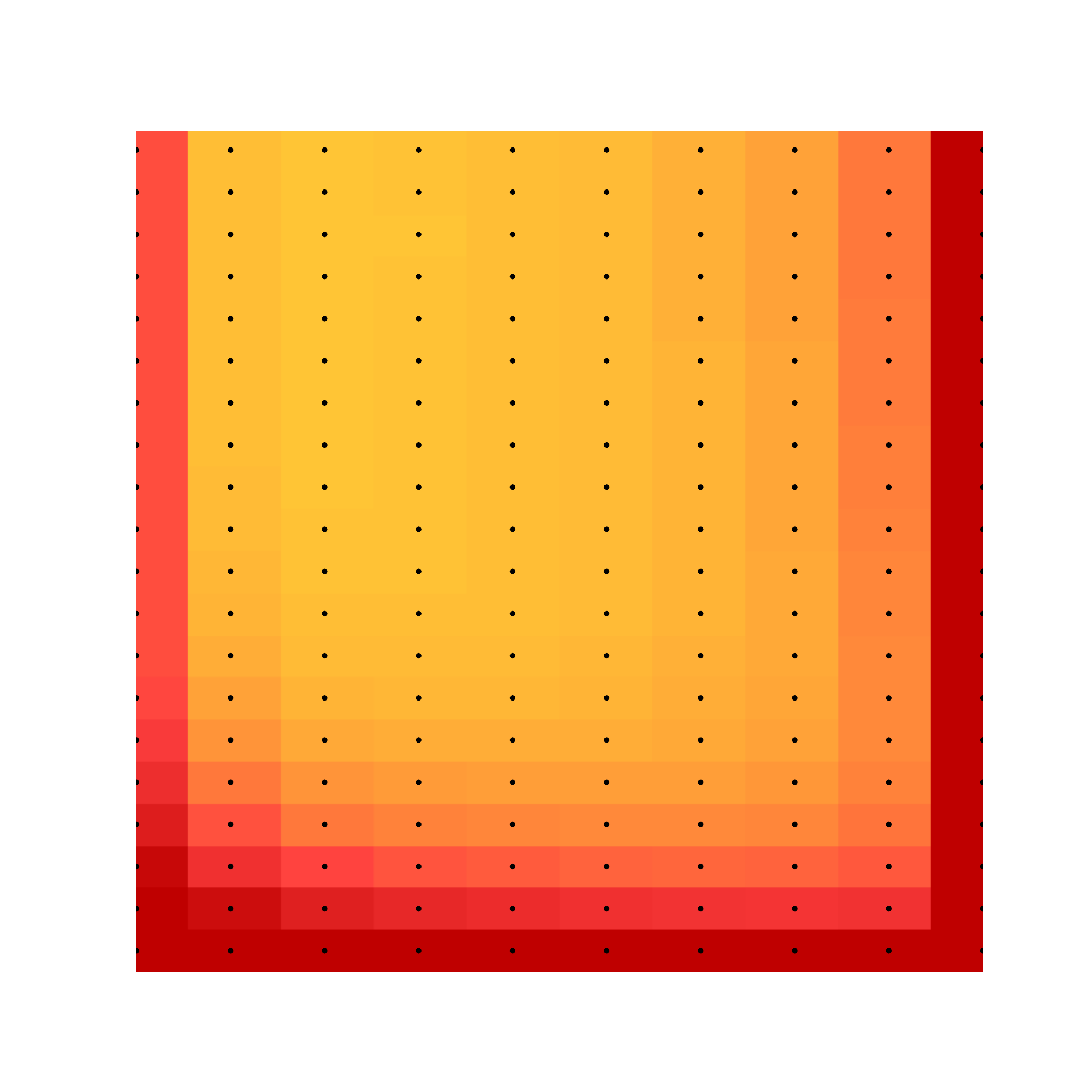}{0}{1}{}} &
\vspace{-0.5cm}\resizebox{.24\textwidth}{!}{\maakmooieticks{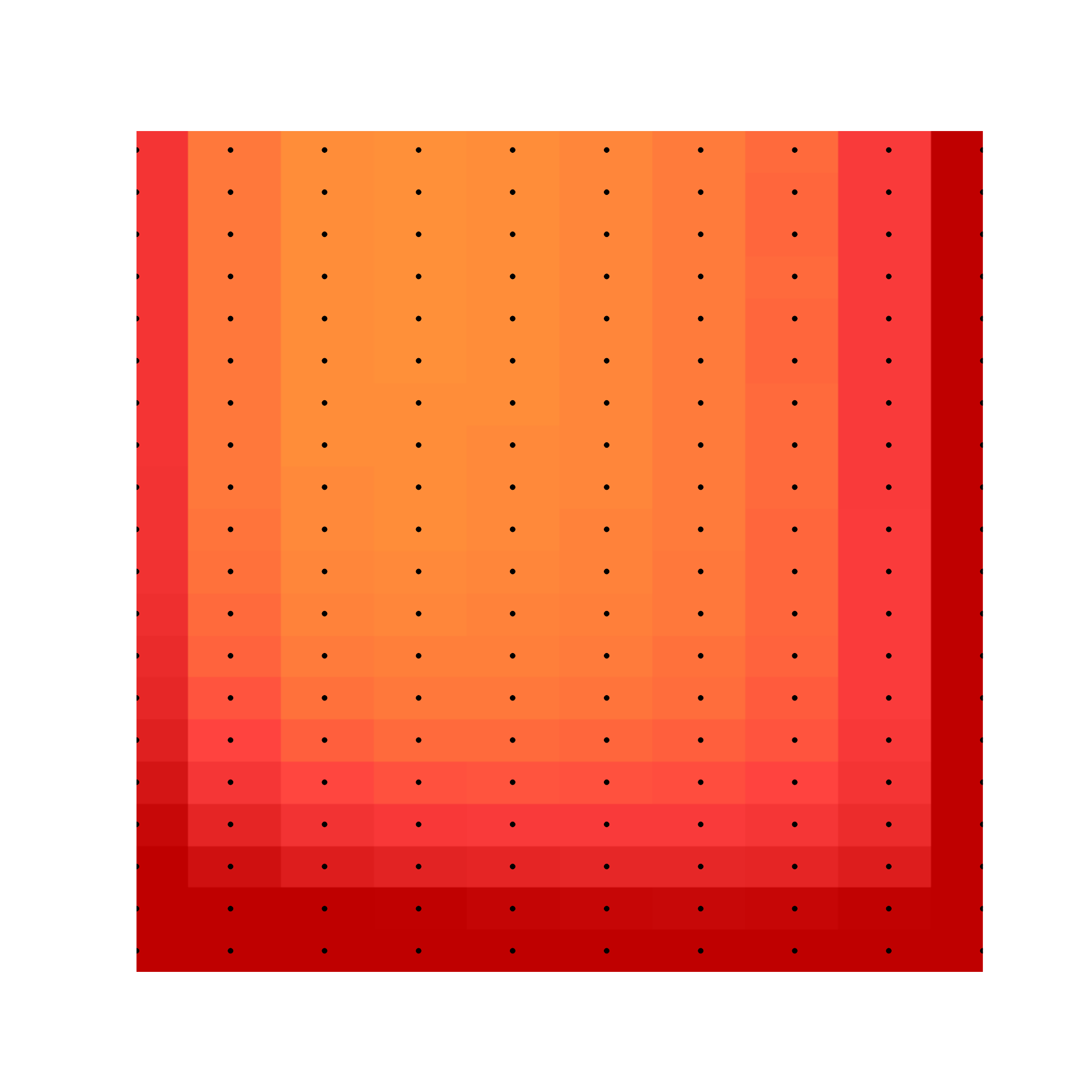}{0}{1}{}} &
\resizebox{!}{.24\textwidth}{\maakmooietickscb{figuren/improvement/improvement_colorbar.pdf}{}}
\end{tabular}
\caption{The best 1D0D momentum source term estimation procedure and the potential gain when considering the statistical error. The first row presents a partition of the parameter domain based on the estimation procedure with the lowest statistical error. The dots represent the exact position in parameter space where the simulation occurred, crosses indicate inconclusive results. The second row presents the factor by which the standard deviation increases when the standard mass choice, \texttt{a\_tl}, is used instead of the best estimation procedure.}
\label{fig:iif_1D0D_mom_var}
\end{figure}

\begin{figure}\centering
\begin{tabular}{m{.25\textwidth}m{.25\textwidth}m{.25\textwidth}m{.2\textwidth}}
\resizebox{.24\textwidth}{!}{\maakmooieticks{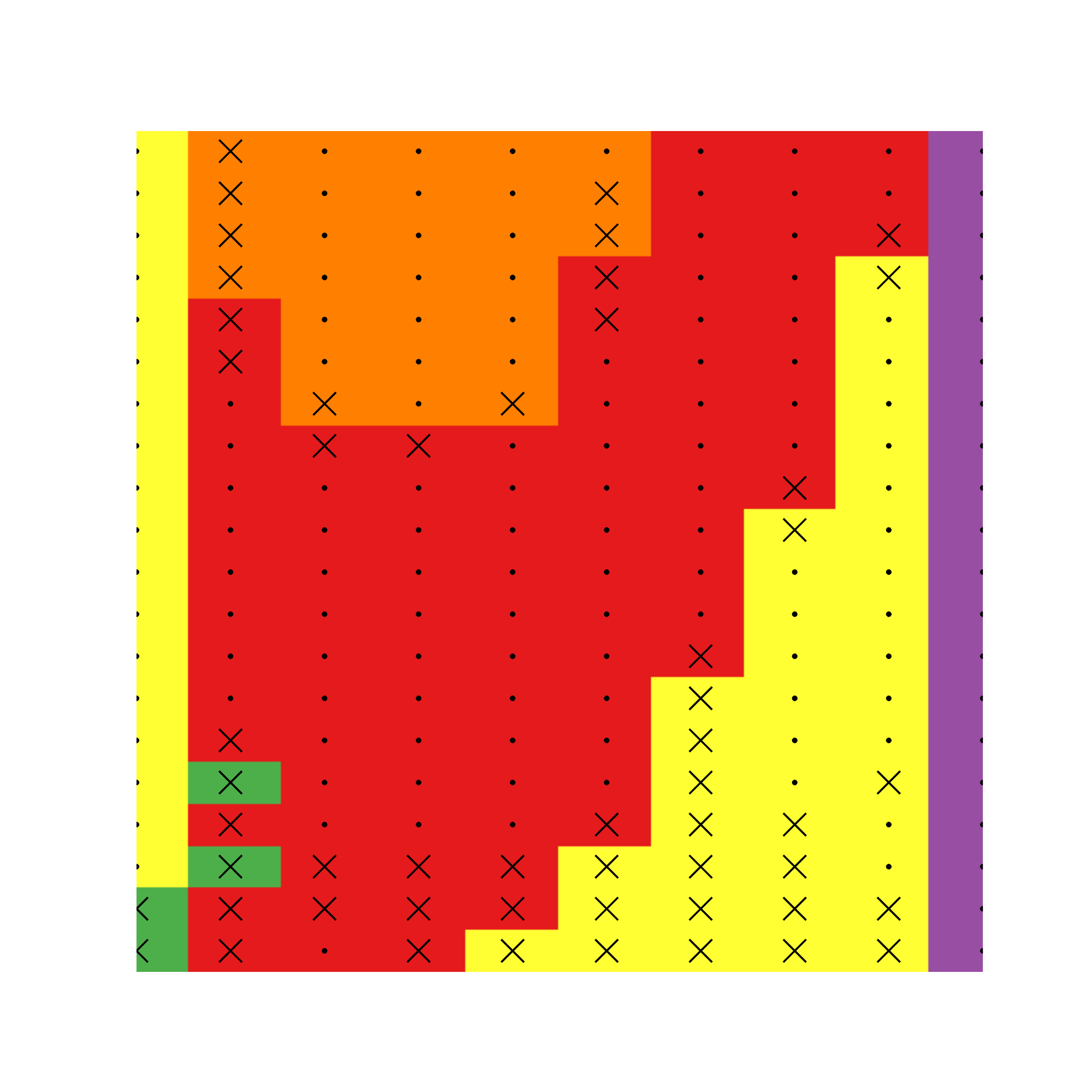}{1}{0}{$\dfrac{\VAcss}{\VAcst}=0.94$}} &
\resizebox{.24\textwidth}{!}{\maakmooieticks{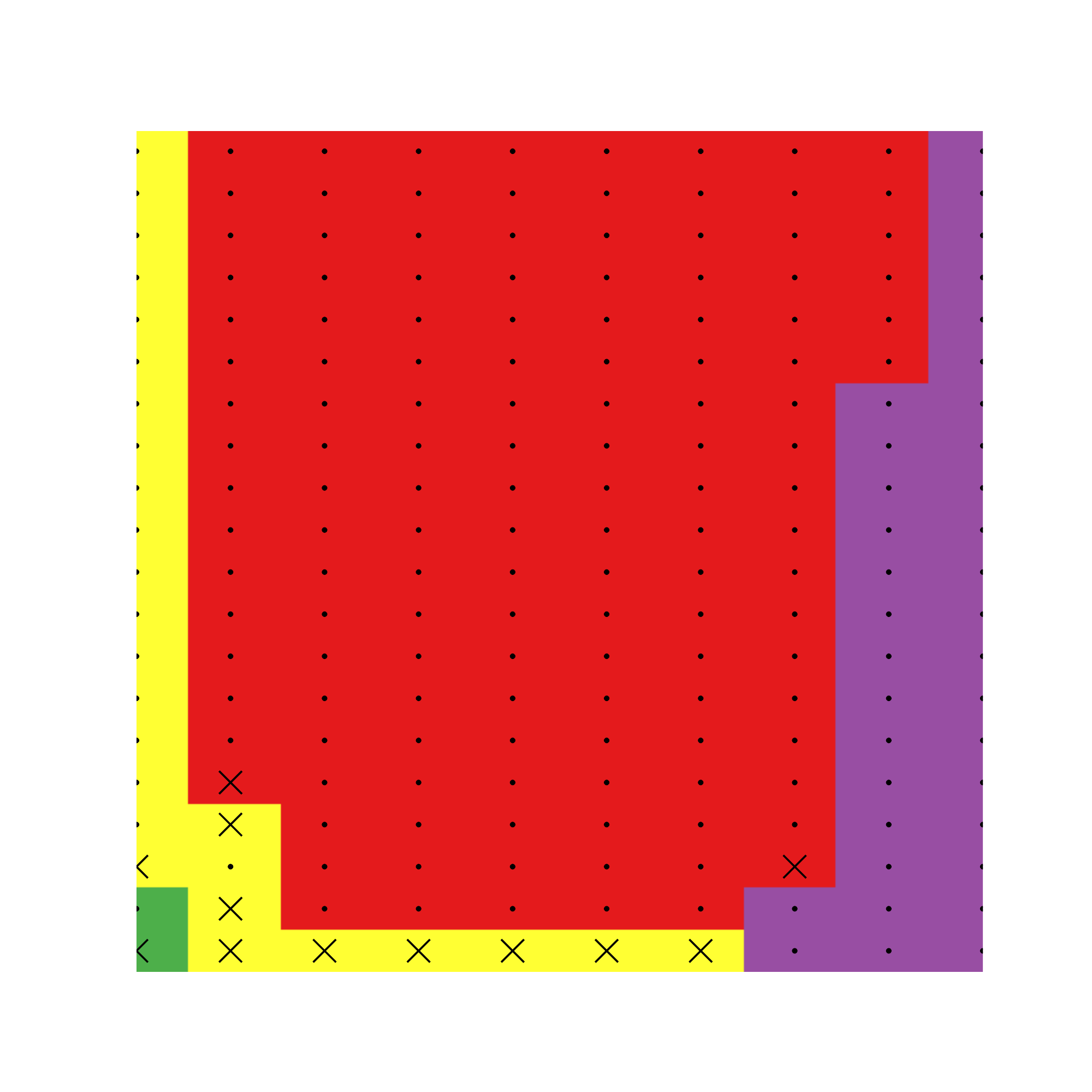}{0}{0}{$\dfrac{\VAcss}{\VAcst}=0.5$}} &
\resizebox{.24\textwidth}{!}{\maakmooieticks{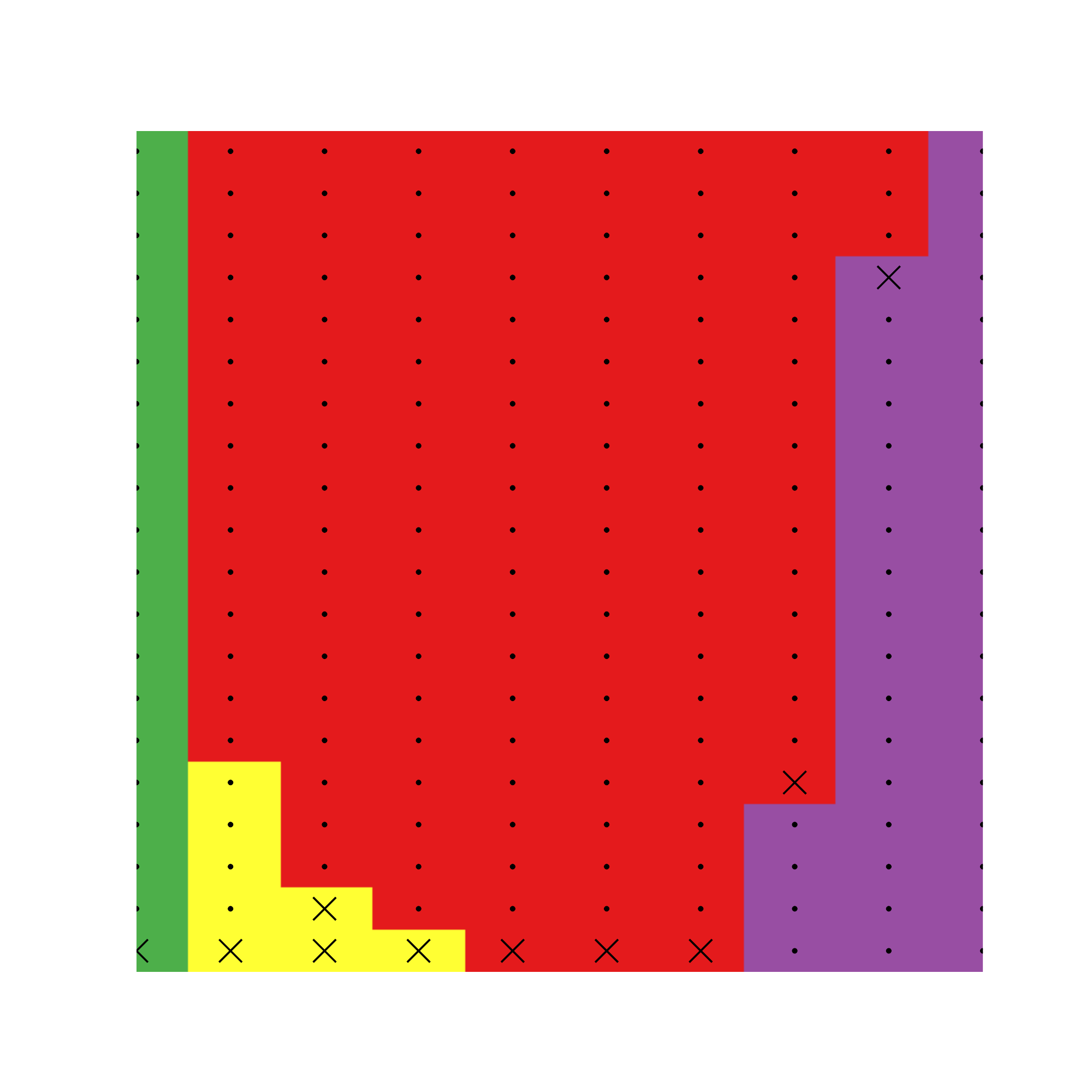}{0}{0}{$\dfrac{\VAcss}{\VAcst}=0.25$}} &
\localcolorlegendoneD\\
\end{tabular}\\
\begin{tabular}{b{.25\textwidth}b{.25\textwidth}b{.25\textwidth}b{.2\textwidth}}
\vspace{-0.5cm}\resizebox{.24\textwidth}{!}{\maakmooieticks{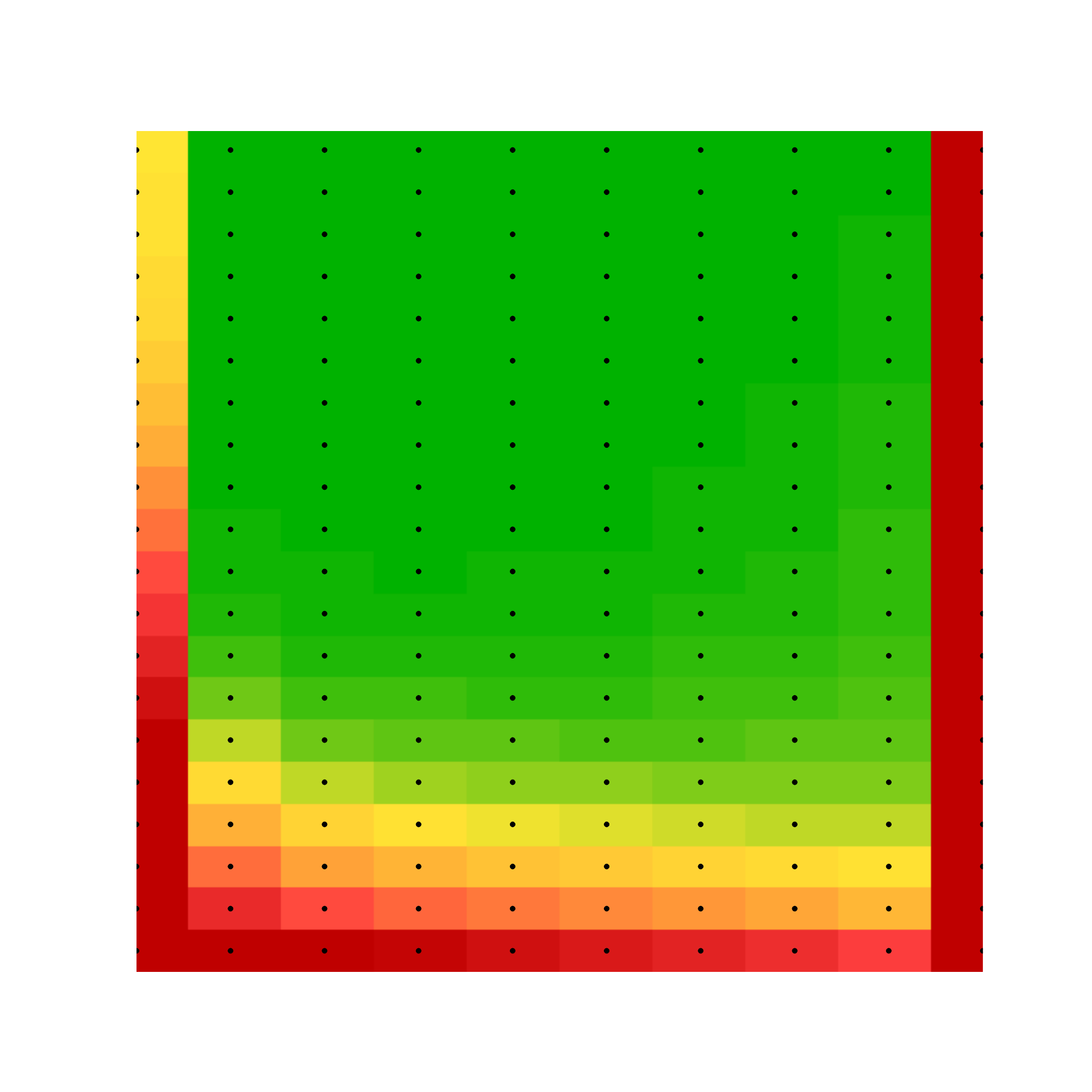}{1}{1}{}} &
\vspace{-0.5cm}\resizebox{.24\textwidth}{!}{\maakmooieticks{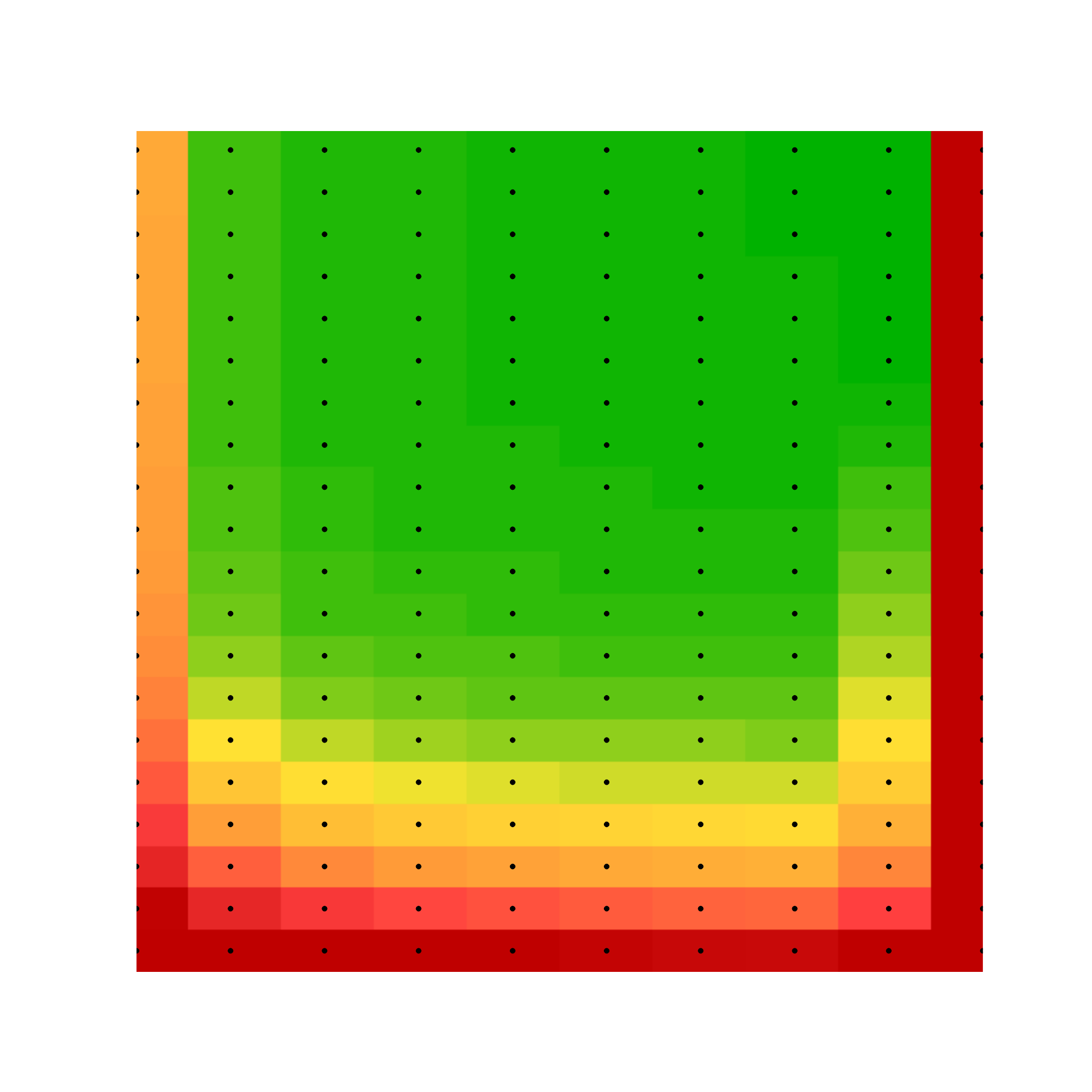}{0}{1}{}} &
\vspace{-0.5cm}\resizebox{.24\textwidth}{!}{\maakmooieticks{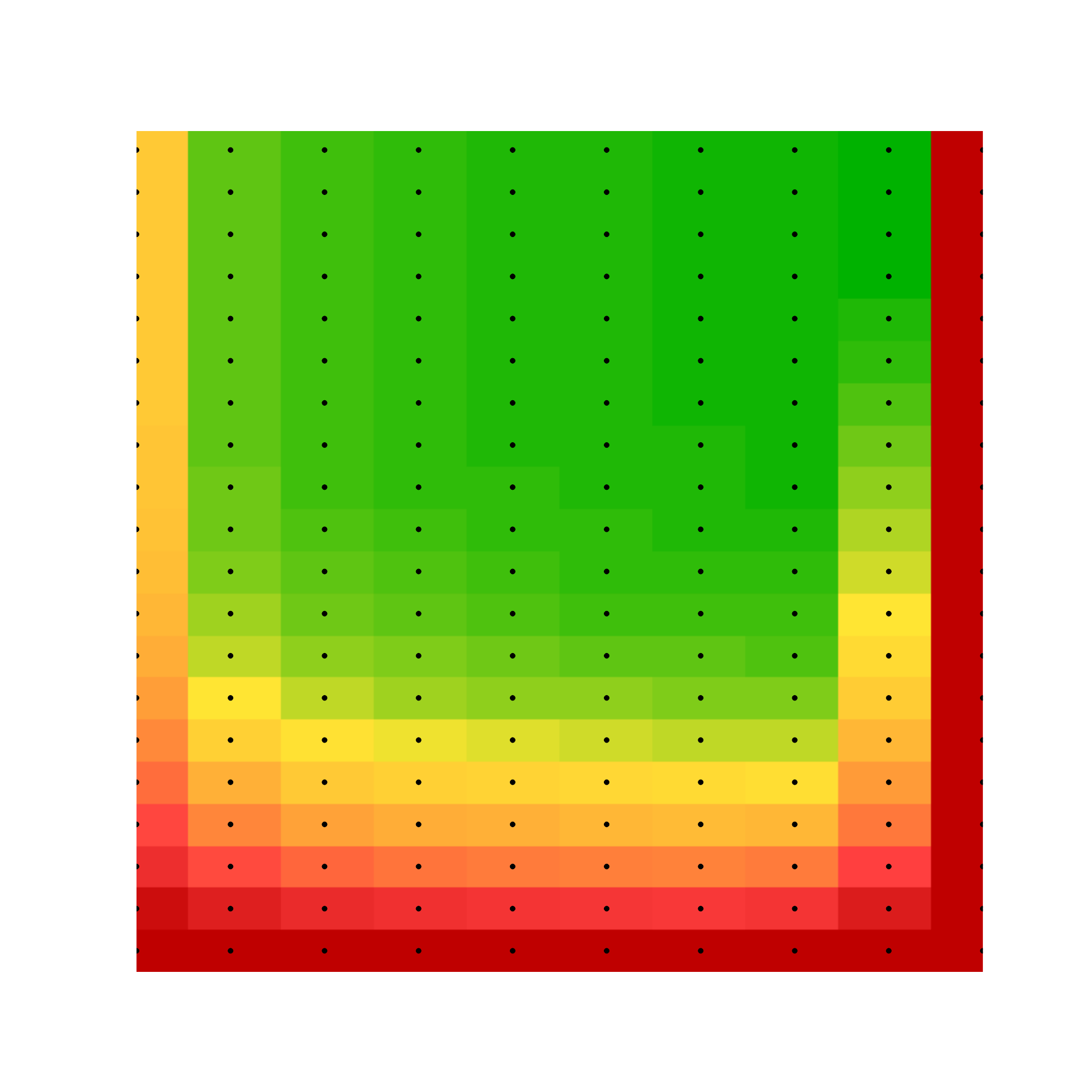}{0}{1}{}} &
\resizebox{!}{.24\textwidth}{\maakmooietickscb{figuren/improvement/improvement_colorbar.pdf}{}}
\end{tabular}
\caption{The best 1D0D momentum source term estimation procedure and the potential gain when considering computational cost for a given statistical error. The first row presents a partition of the parameter domain based on the estimation procedure with the lowest computational cost error. The dots represent the exact position in parameter space where the simulation occurred, crosses indicate inconclusive results. The second row presents the factor by which the computational cost increases when the standard mass choice, \texttt{a\_tl}, is used instead of the best estimation procedure.}
\label{fig:iif_1D0D_mom_cost}
\end{figure}

To evaluate the representativeness of the 1D0D results as predictors for the 1D1D results, we have performed the same experiments for the 1D0D case, the results of which are shown in Figures~\ref{fig:iif_1D0D_mom_var} and~\ref{fig:iif_1D0D_mom_cost}. When we compare Figures~\ref{fig:iif_1D0D_mom_var} and~\ref{fig:iif_1D0D_mom_cost} with Figures~\ref{fig:iif_1D1D_mom_var} and~\ref{fig:iif_1D1D_mom_cost}, we note remarkable agreement, but with several shifts in which estimator performs best, and the addition of \texttt{a\_c} as a potentially competitive estimator in the 1D0D case.

The appearance of collision estimators (\texttt{nac\_c} and \texttt{a\_c}) as best estimation procedures is to be contrasted to the proven fact by Lux~\cite{lux1991MCPT} that collision estimators can never outperform next-event estimators for mass source estimation. This is now clearly shown to not hold for momentum source estimation.

In this section we have shown that different estimation procedures are optimal for momentum source estimation than for mass source estimation. Typically, we estimate both with the same simulation. Since the noise on the higher moments is typically larger, the selection of the momentum estimation procedure can have a higher priority. The simulation type for the mass estimation procedure is then fixed, but the estimator can still be chosen freely.

\section{Conclusion and future prospects\label{sec:iif_conclusion}}

We compared the performance of a coherent set of source term estimation procedures for a particle tracing Monte Carlo method for the Boltzmann-BGK kinetic equation. We first considered the best estimation procedure with the analytical results for mass source estimation in a simplified 1D0D setting with forward-backward scattering, which were obtained via the invariant imbedding methodology in~\cite{mortier2020iiappendix}. Our results show three estimation procedures are competitive for mass source estimation in a 1D0D setting when variance is considered (\texttt{natl\_tl}, \texttt{nac\_ne}, and \texttt{natl\_ne}) and four when cost is considered (\texttt{a\_ne}, \texttt{natl\_tl}, \texttt{nac\_ne}, and \texttt{natl\_ne}). We have also shown that the potential profit of using the optimal estimation procedure instead of the standard choice and discussed the unboundedness of the potential gain when the background parameters are unknown. Numerical experimentation showed that our analysis for the mass source term in a 1D0D setting remains relevant in a multi-velocity setting (1D1D): the same estimation procedures remain optimal, but there are some shifts in the regions in which they are optimal due to the changed neutral model.

The mass source estimation results also clearly show the currently used default estimation procedure (\texttt{a\_tl}) in the EIRENE code is a valid choice in the regime that is relevant for the application. This is to be expected, since this estimation procedure has been selected by trial-and-error. A next-event estimator, or a more flexible, region-dependent choice of estimation procedure, might provide sharp improvements.

As an additional numerical extension, we compared the different estimation procedures for momentum source estimation. Compared with the results for mass estimation, the optimal estimation procedure changes more drastically, especially when cost is considered as a measure of performance. Now five or seven estimation procedures are optimal when variance, respectively cost is considered. Again, the current default momentum source estimation option, \texttt{a\_c}, is shown to be a proper choice according to our study.

Our results provide an indicative tool to select the proper estimation procedure without having to resort to trial-and-error. This enables the selection of different estimation procedures for different regions of the domain, which can be extremely relevant when the domain is highly heterogeneous or when grid refinement drastically changes the total collisionality per grid cell. Due to the unbiasedness of the estimators and the consistency of the different simulation types, there is no restriction on using different estimation procedures even conditionally on the velocity at entry in the grid cell.

Higher dimensions, or energy source term estimation are possible extensions of this work. These form future work. Since neutrals are transported mostly along lines, the one-dimensional analysis, perhaps on the quantitative level, and definitely the qualitative results are expected to hold in higher dimensions. 

Our comparison has been restricted to a single grid cell, hence another logical expansion of our results is to multiple grid cells. This will entail non-local effects due to the simulation types. In an analog simulation the probability of reaching a grid cell is restricted compared to non-analog simulations, which, in turn, experience adverse effects due to their weight distribution. Non-analog simulations entail small-weight particles that will have little impact on the estimate, but still require an equally costly simulation as particles with a higher weight. Evidently, non-analog simulation types will hold the possibility for the lowest variance, but when regarding the measure of cost, they are expected to lose performance with respect to analog estimation procedures compared to the single grid cell situation.

\section*{Acknowledgements}

The first author is funded by a PhD fellowship of the Research Foundation Flanders (FWO) under fellowship number 1189919N.\\
This work has been carried out within the framework of the EUROfusion Consortium and has received funding from the Euratom research and training programme 2014-2018 and 2019-2020 under grant agreement No 633053. The views and opinions expressed herein do not necessarily reflect those of the European Commission.\\
Parts of the work are supported by the Research Foundation Flanders (FWO) under project grant G078316N.\\
The computational resources and services used in this work were provided by the VSC (Flemish Supercomputer Center), funded by the Research Foundation Flanders (FWO) and the Flemish Government, department EWI.

\bibliographystyle{abbrv}
\bibliography{referenties}

\end{document}